\documentclass{amsart}
\usepackage{amsmath}

\usepackage{graphicx}
\newtheorem{defn}{Definition}
\newtheorem{con}{Conjecture}

\newtheorem{alg}{Algorithm}

\begin{document}

\title{Lozenge Tilings Of the Equilateral Triangle}

\author{Richard J. Mathar}
\urladdr{https://www.mpia-hd.mpg.de/homes/mathar}
\address{Max-Planck Institute of Astronomy, K\"onigstuhl 17, 69117 Heidelberg, Germany}

\date{\today}

\begin{abstract}
We consider incomplete tilings of the equilateral triangle of edge length $n$
that is subdivided into $n^2$ regular equilateral smaller unit triangles.
Pairs of the unit triangles that share a side may be converted
into lozenges, leaving some subset of the unit triangles untouched.
We count numerically these coverings by lozenges and unit triangles for edge
lengths $n\le 15$: the total and the detailed refinement as a function of the 
number of lozenges.
\end{abstract}

\maketitle

\section{Lozenge Tilings}
\subsection{Basic Geometry}

An equilateral triangle of integer side length $n$
may be divided into $n^2$ equilateral triangles of unit side length
by regular subdivision of each side into $n$ sections 
and drawing lines through these parallel to all three sides.
This creates a graph with \cite[A000217]{sloane}
\begin{equation}
(n+1)+n+(n-1)+\cdots +1 = T_{n+1}=\frac{(n+1)(n+2)}{2}
\label{eq.V}
\end{equation}
vertices (corners of the unit triangles).

\begin{defn}(Triangular Numbers)
\begin{equation}
T_n\equiv \left\{ \begin{array}{ll}
n(n+1)/2,& n\ge 0,\\
0,& n<0.
\end{array}\right.
\end{equation}
\end{defn}
The graph has \cite[A045943]{sloane} 
\begin{equation}
M_n=3T_n = 3\frac{n(n+1)}{2}
\label{eq.v}
\end{equation}
edges (edges of the unit triangles). 
Euler's Formula \cite{LakatosBJPS14} states that 
the number of faces plus the number of vertices equals the number of edges plus 1:
\begin{equation}
n^2+ T_{n+1} = M_n+1.
\end{equation}
The number of edges on the perimeter of the big triangle is three times the number
of segments, $3n$, so the number of edges internal to the big triangle is
\begin{equation}
M_n-3n= M_{n-1}.
\label{eq.inte}
\end{equation}

The number of vertices \emph{inside} the big triangle is
the number of vertices which are not on one of the sides of the big triangle;
so subtracting  $3n$, the number of vertices on the big triangle's sides,
from \eqref{eq.V} yields  the number of internal vertices:
\begin{equation}
T_{n+1}-3n=T_{n-2}.
\label{eq.intv}
\end{equation}

The triangle of side length $n$ is considered fixed with one of the three corners 
pointing up. It contains
$T_n$ unit triangles pointing up and $n^2-T_n=T_{n-1}$ unit triangles pointing down.

\subsection{Lozenge Sets} \label{sec.corr}
In conjunction with this work, a lozenge is created by removing one of the
inner edges; this merges the two unit adjacent triangles that have that edge in common. 
A lozenge tiling with $l$ non-overlapping lozenges is created by removing
$l$ of the inner edges under the constraint that no pair of removed
edges must be two edges of the same triangle---which would create
tiles that are larger than a lozenge. So the constraint
means that once an inner edge has been removed (to become the short diagonal
of a unit lozenge), none of the 4 edges of that lozenge must be removed.

If $l$ is the number of lozenges, $n^2-2l$ is the number
of free triangles. An obvious upper bound for $l$ is  the ``capacity''
\begin{equation} 
l\le n^2/2 
\end{equation}
 because
each lozenge covers 2 triangles.

\begin{defn}
$L_{n,l}$ is the number of tilings of the equilateral triangle
with edges of length $n$ with $l$ non-overlapping lozenges (and
in consequence $n^2-2l$ unit triangles not covered by lozenges).
\end{defn}

\begin{alg}\label{alg.1}
A simple strategy to count the tilings $L_{n,l}$
is to generate the set of $M_{n-1}$ inner edges, to scan all $2^{M_{n-1}}$
subsets of removing them, and to count all the subsets that meet the criterion
that no pair of removed edges is part of the same triangle. If the constraint
were absent, the number of subsets follows from the usual combinatorial selection,
so with (\ref{eq.inte}) this constitutes an upper bound
\begin{equation}
L_{n,l} \le \binom{M_{n-1}}{l}.
\label{eq.bino}
\end{equation}
\end{alg}

The lozenges have three different orientations with axes differing
by angles of 120$^\circ$. We classify them according to the removed edge being
horizontal, falling left-to-right or rising left-to-right.
If one takes a set of lozenges of a common orientation and shoves them
in closest packing into a corner of the big triangle, one sees that
a tiling with 
\begin{equation}
l=(n-2)+(n-1)+\cdots 1=T_{n-2}
\label{eq.lmax}
\end{equation}
lozenges (and $n$ isolated unit triangles)
is possible:
\begin{equation}
L_{n,T_{n-2}}\ge 1.
\end{equation}

\section{Example: Side length $n=3$}

\begin{figure}
\includegraphics[scale=0.15]{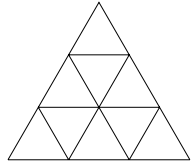}($\times 1$)
\caption{The configuration with 0 lozenges, side length 3  (multiplicity 1).
$L_{3,0}=1$.
There are $n^2=9$ unit triangles.
}
\label{fig.30}
\end{figure}

\begin{figure}
\includegraphics[scale=0.15]{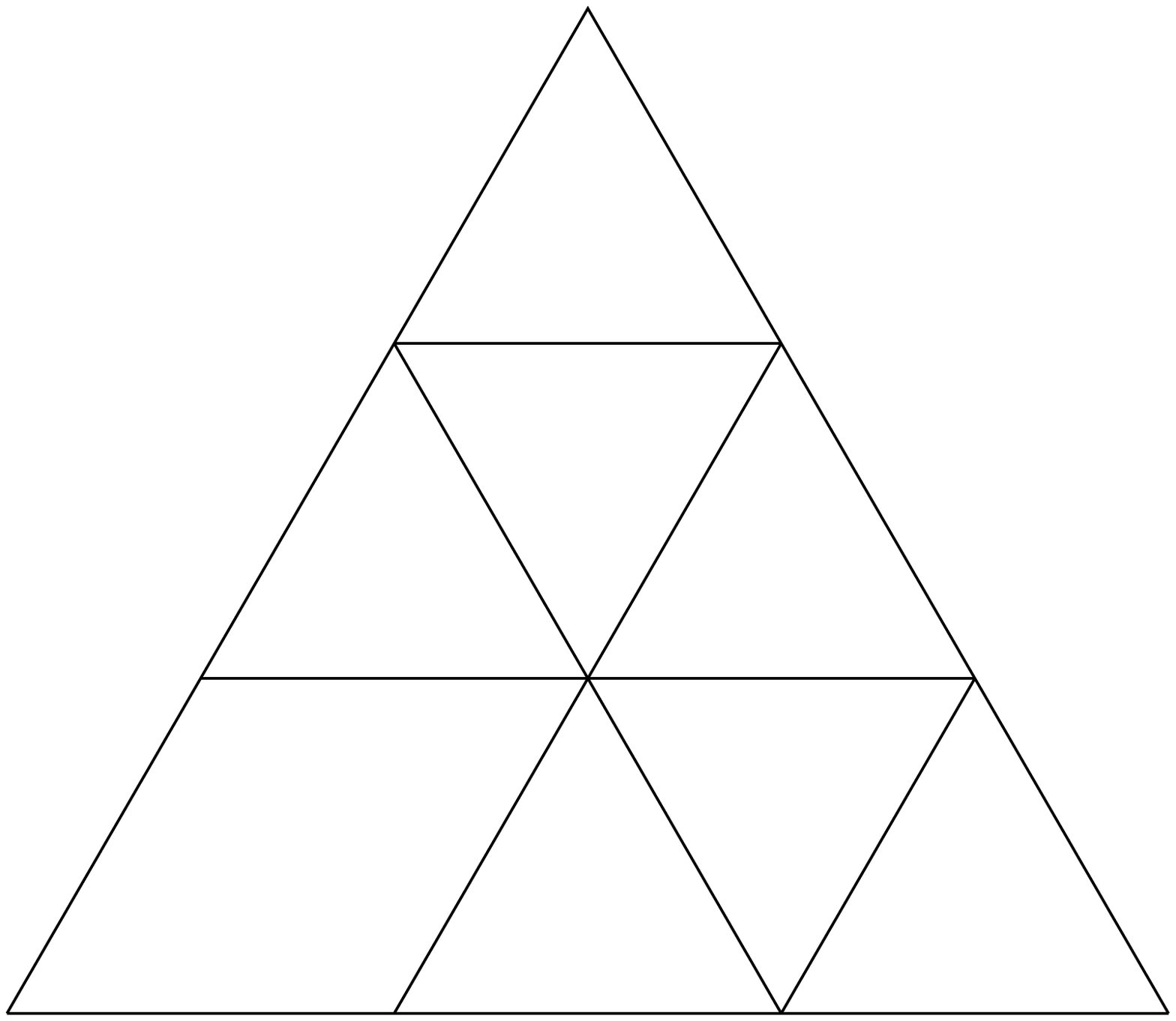}($\times 3$)
\includegraphics[scale=0.15]{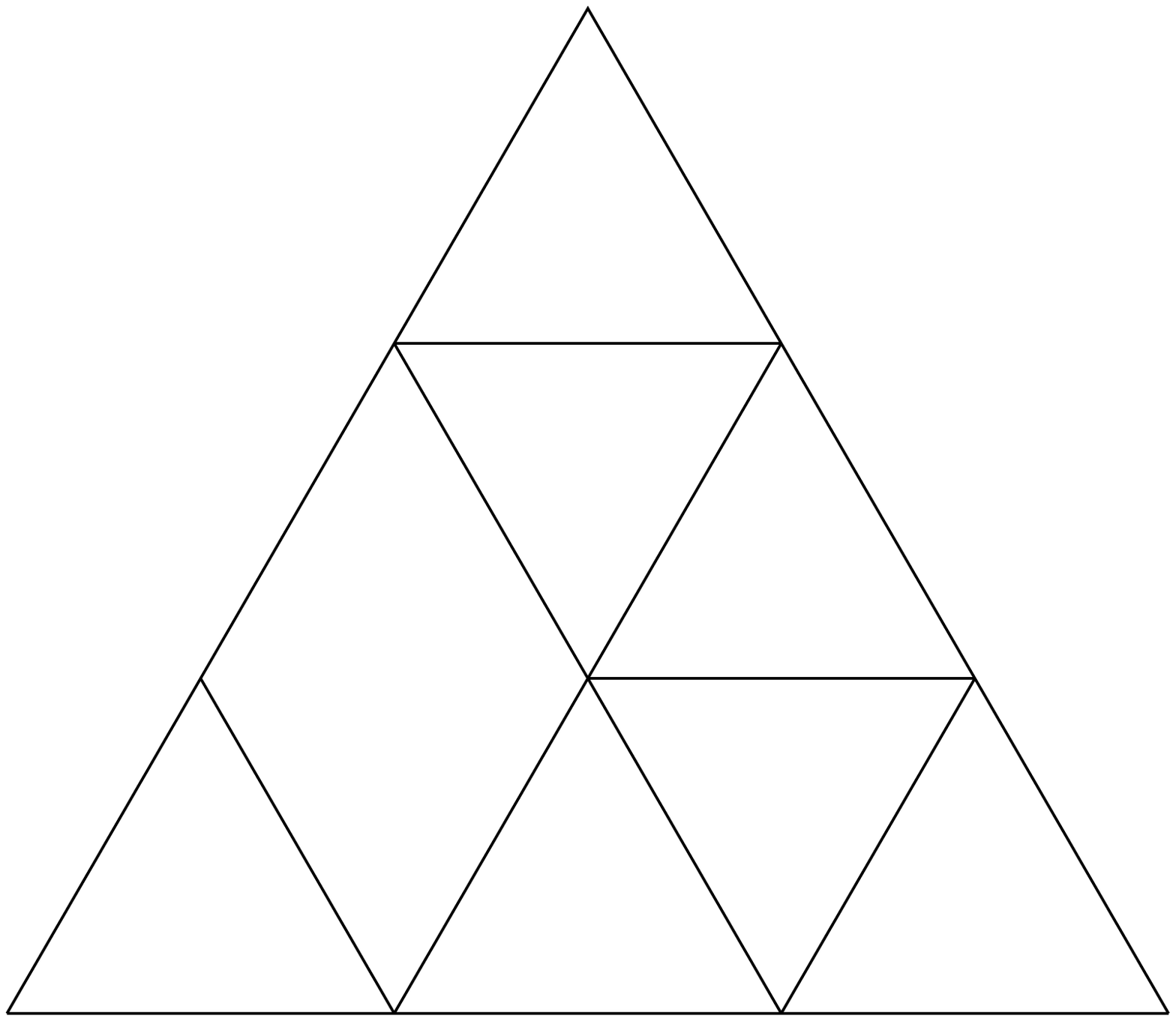}($\times 6$)
\caption{The two configurations with 1 lozenge, side length $n=3$. It either covers a corner
of the big triangle (multiplicity 3) or shares one edge with the middle section
of an edge of the big triangle (multiplicity 6).
$L_{3,1}=3+6=9$}
\label{fig.31}
\end{figure}

\begin{figure}
\includegraphics[scale=0.1]{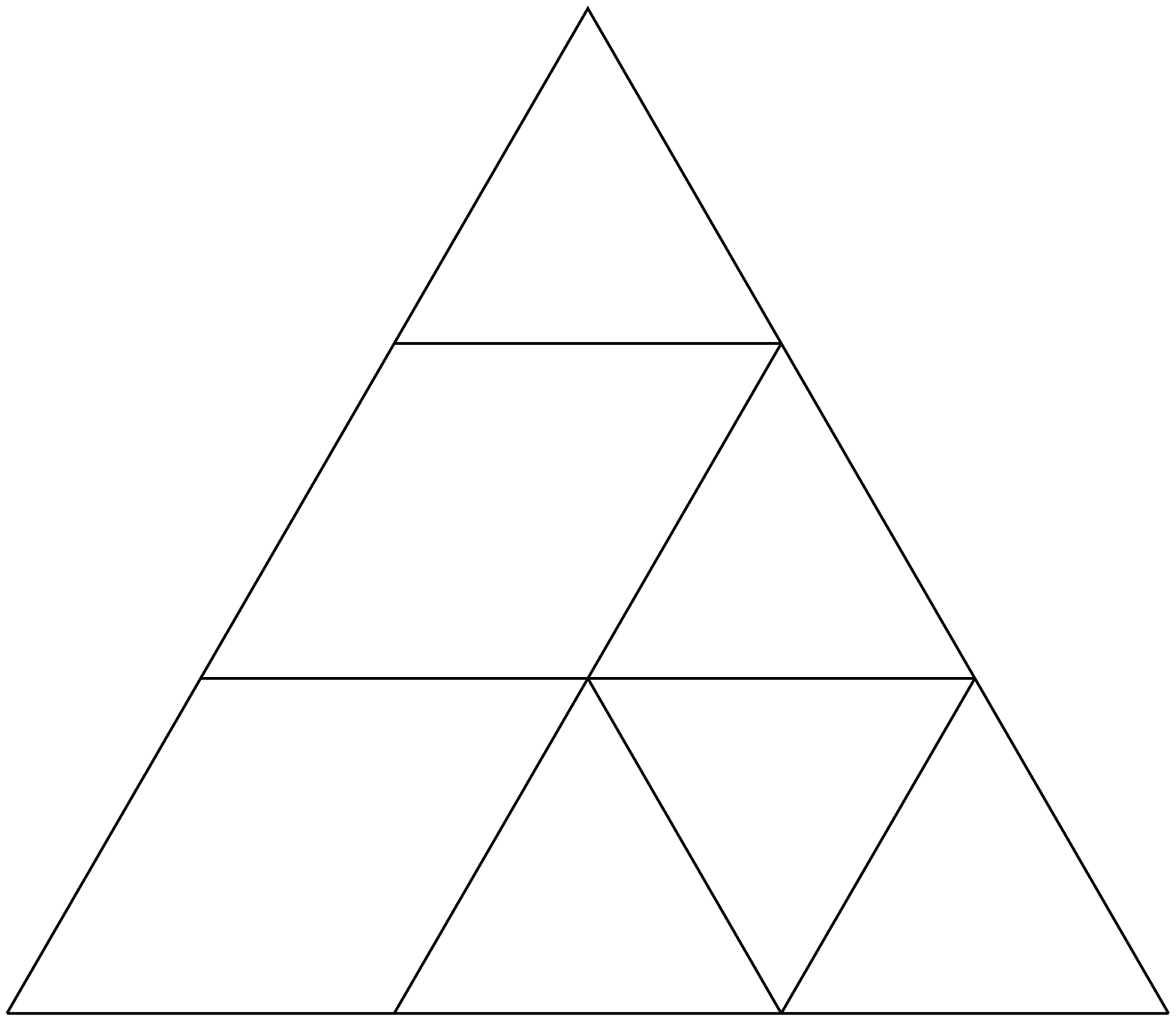}($\times 6$)
\includegraphics[scale=0.1]{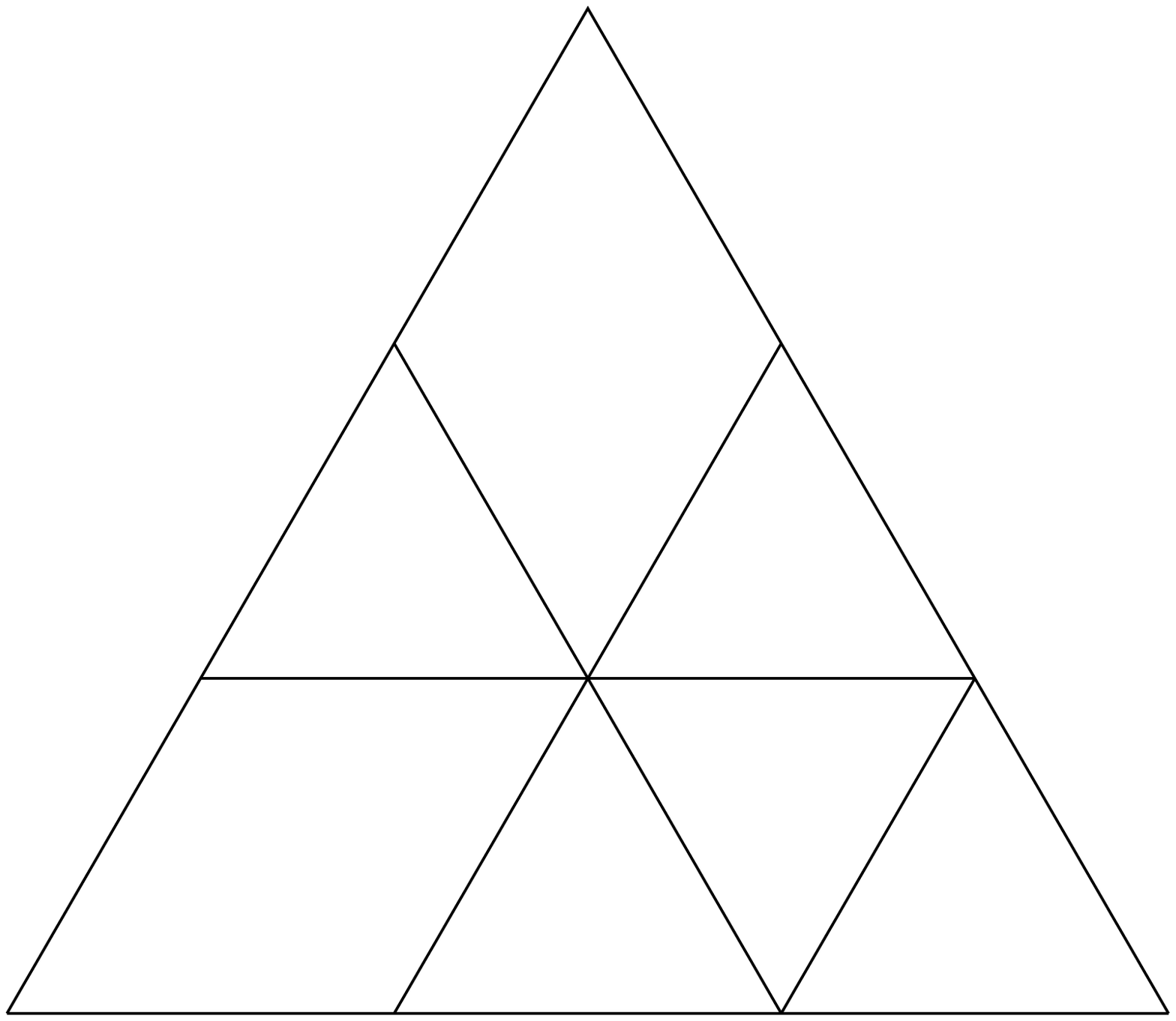}($\times 3$)
\includegraphics[scale=0.1]{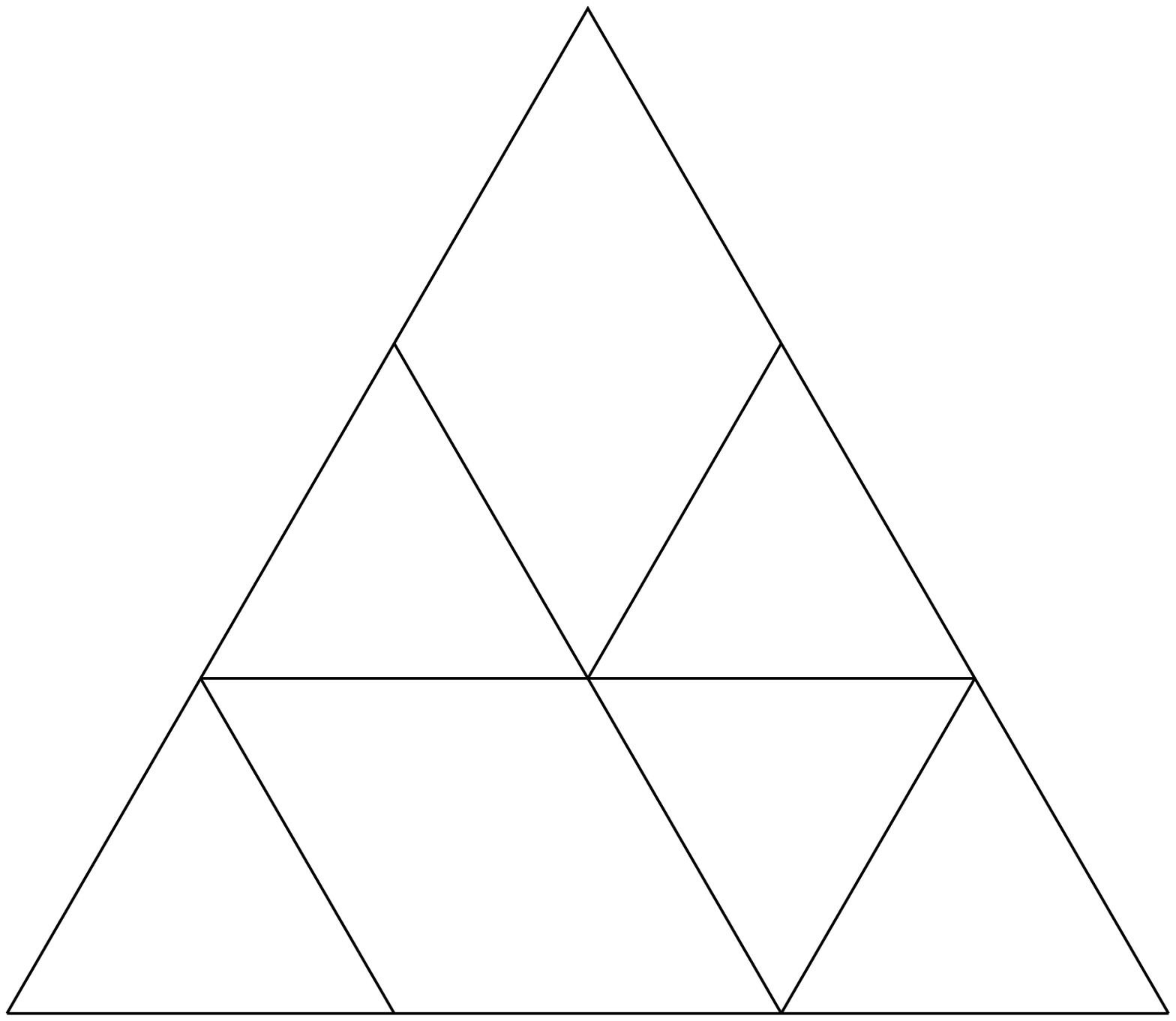}($\times 6$)
\includegraphics[scale=0.1]{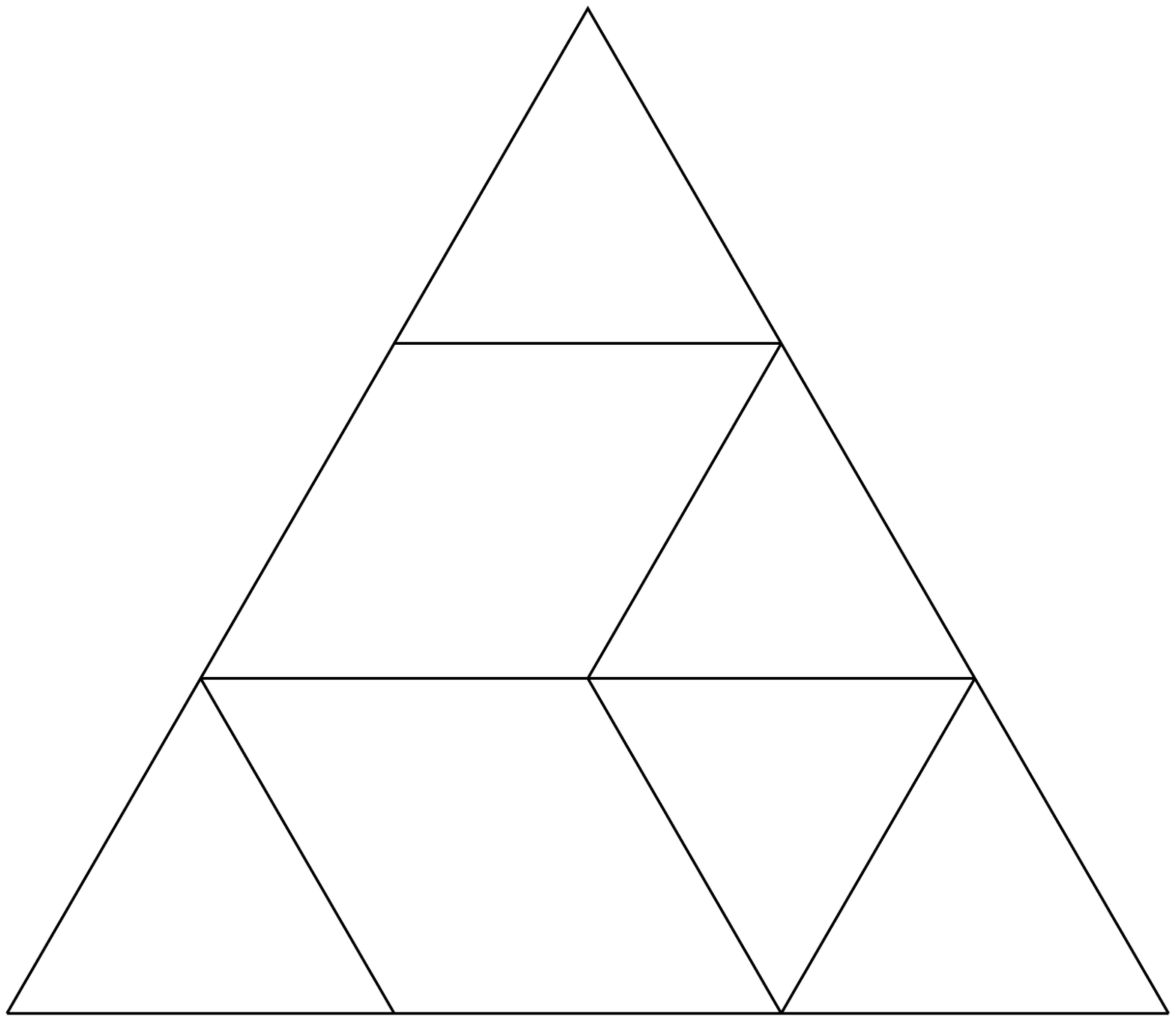}($\times 6$)
\includegraphics[scale=0.1]{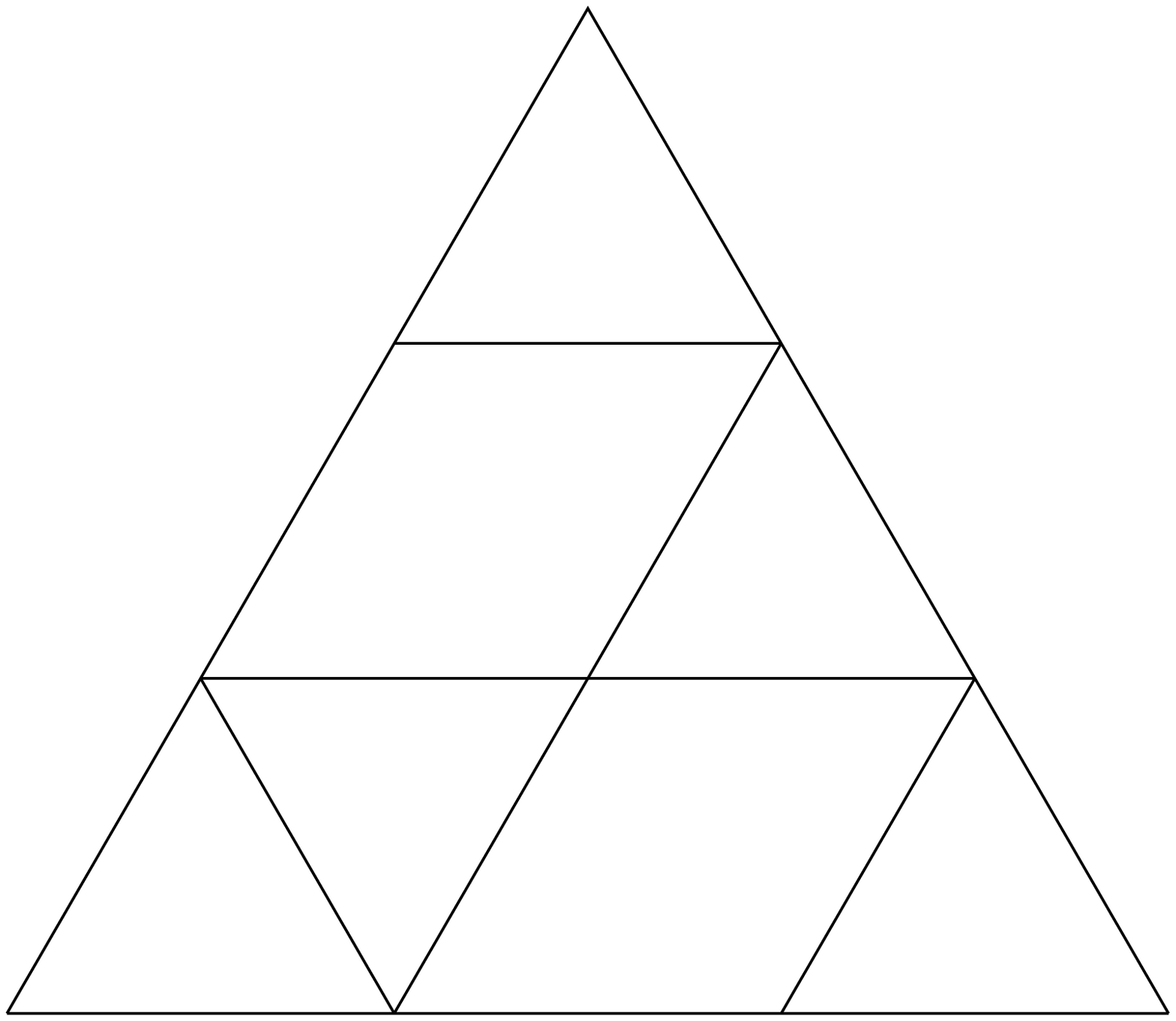}($\times 3$)
\caption{Configurations with 2 lozenges, side length $n=3$. 
They may cover 2/3 of an edge of the big triangle (multiplicity 6).
They may cover two different corners of the big triangle (multiplicity 3).
One may cover a corner of the big triangle and the other the middle section of
the opposite edge (multiplicity 6).
They may share an edge and cover middle sections of two edges of the big triangle
(multiplicity 6).
They may touch in the middle and cover middle sections of two edges of the big triangle
(multiplicity 3).
$L_{3,2}=6+3+6+6+3=24$.
}
\label{fig.32}
\end{figure}

\begin{figure}
\includegraphics[scale=0.1]{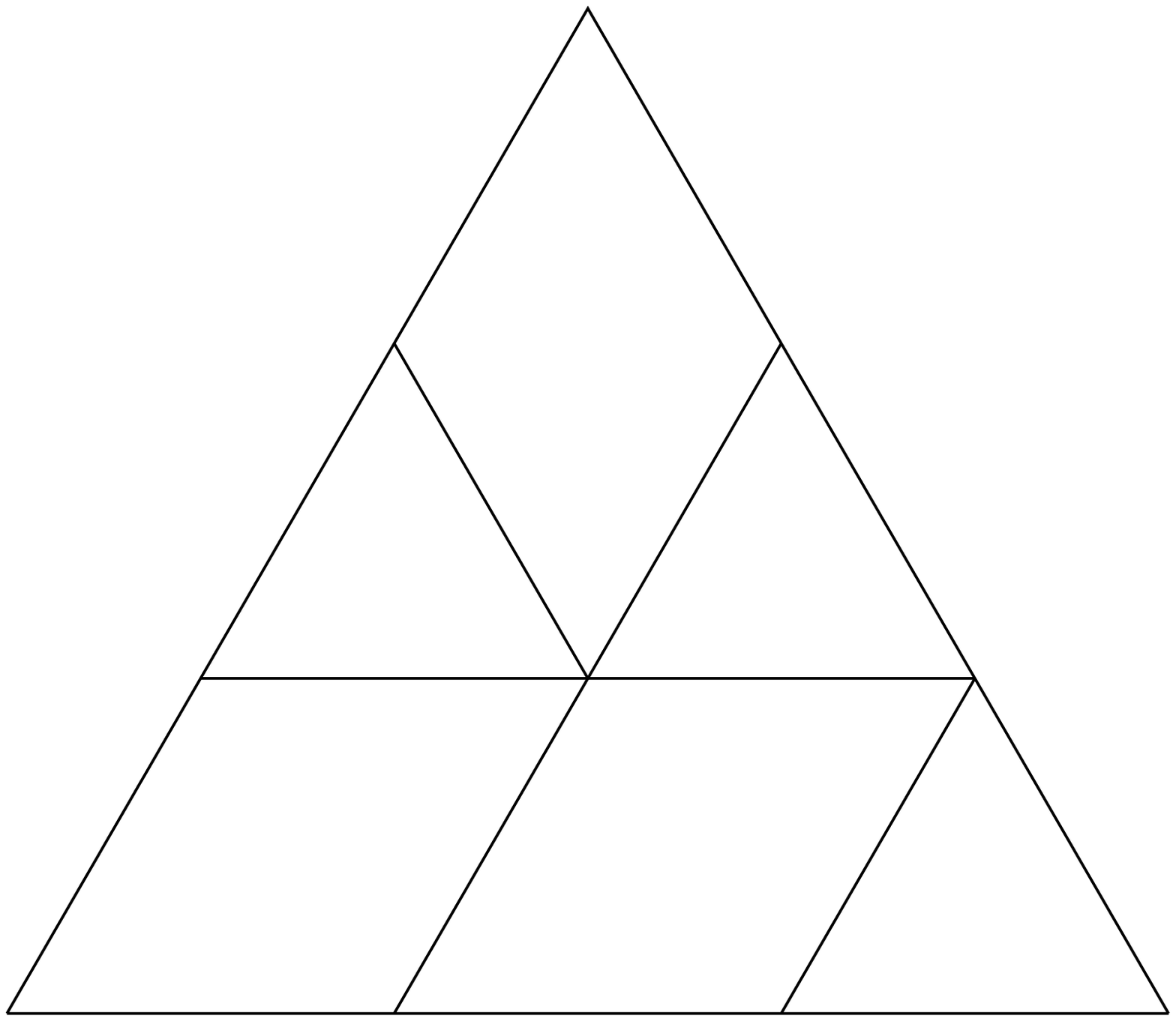}($\times 6$)
\includegraphics[scale=0.1]{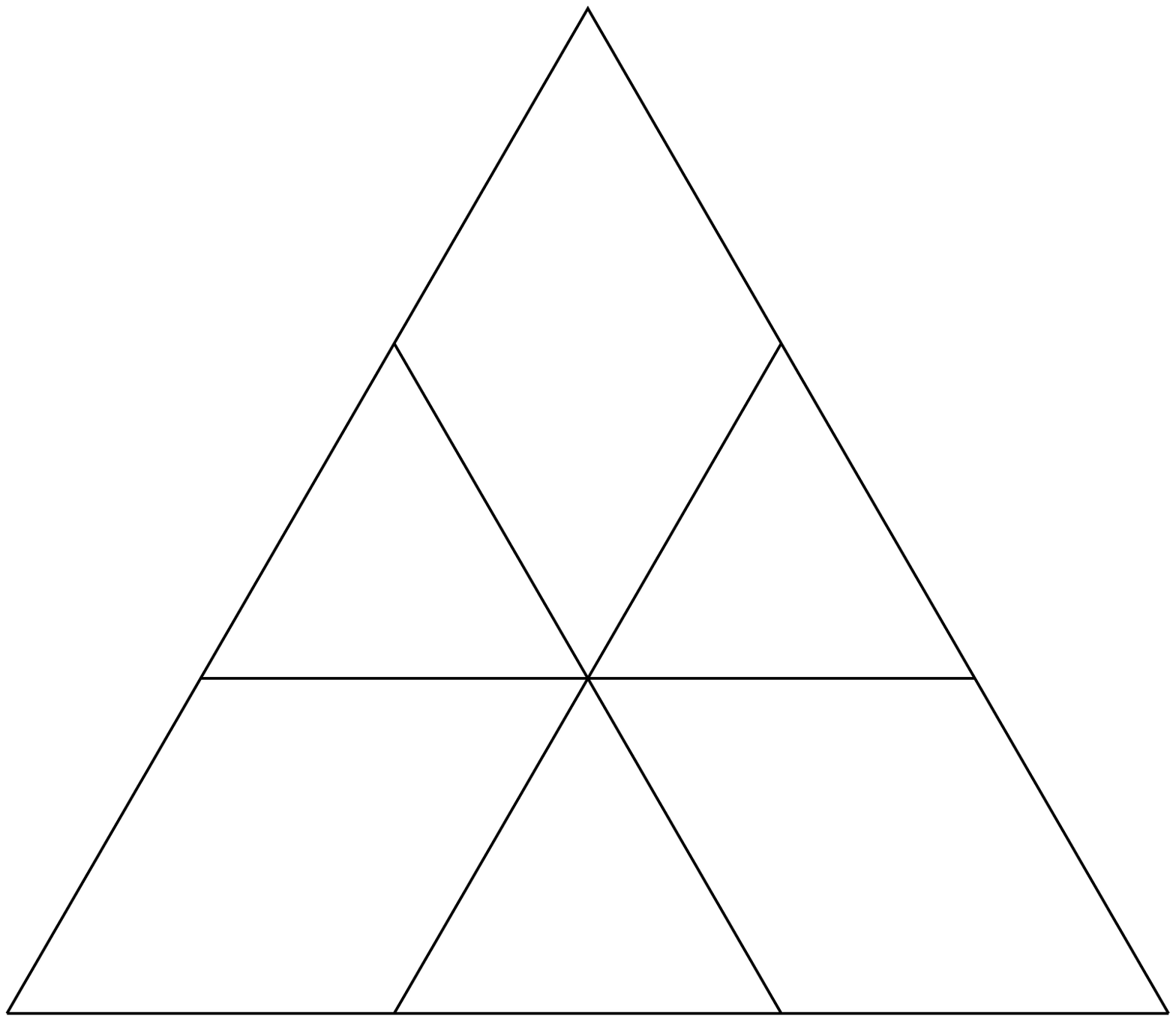}($\times 1$)
\includegraphics[scale=0.1]{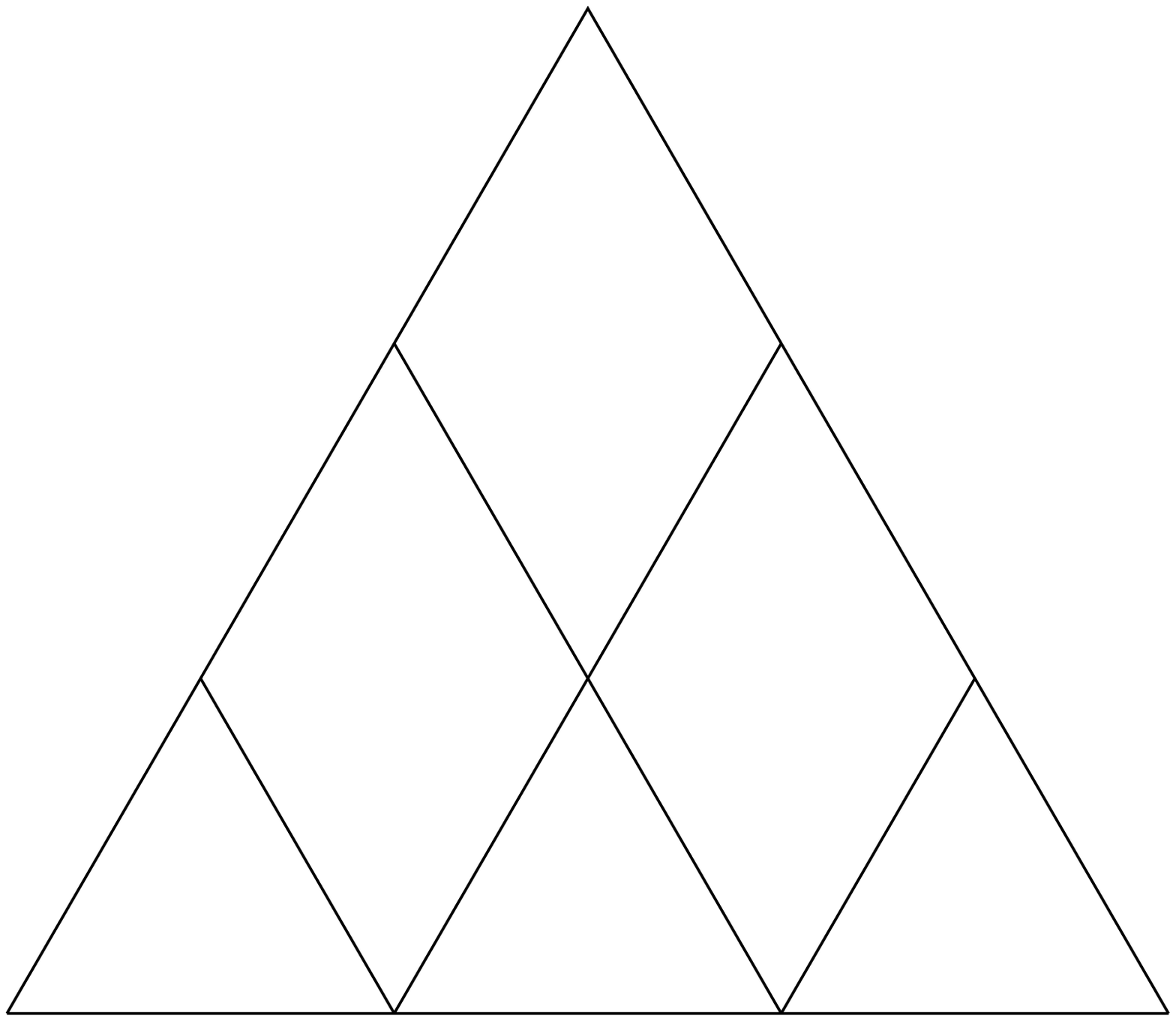}($\times 3$)
\includegraphics[scale=0.1]{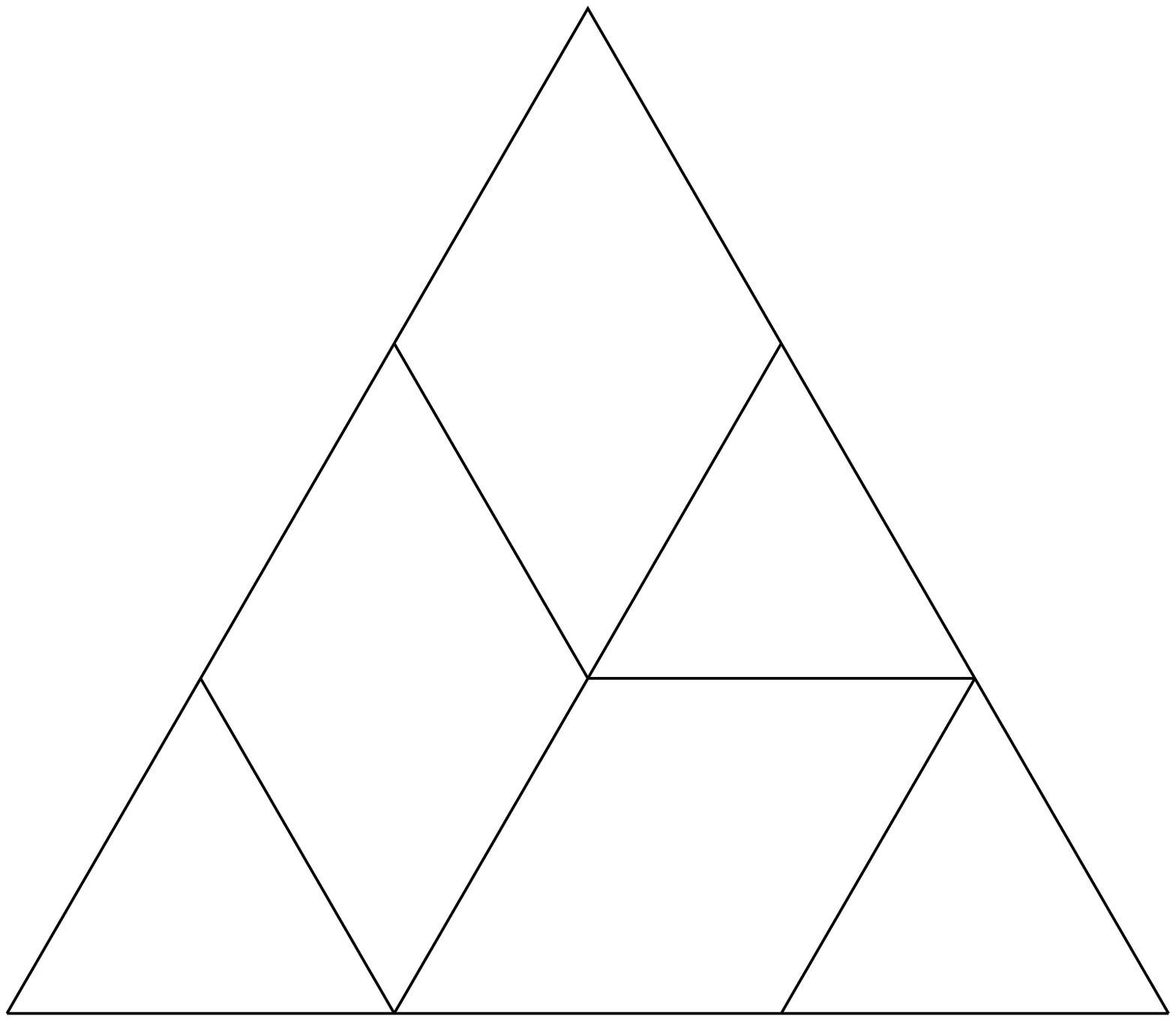}($\times 6$)
\includegraphics[scale=0.1]{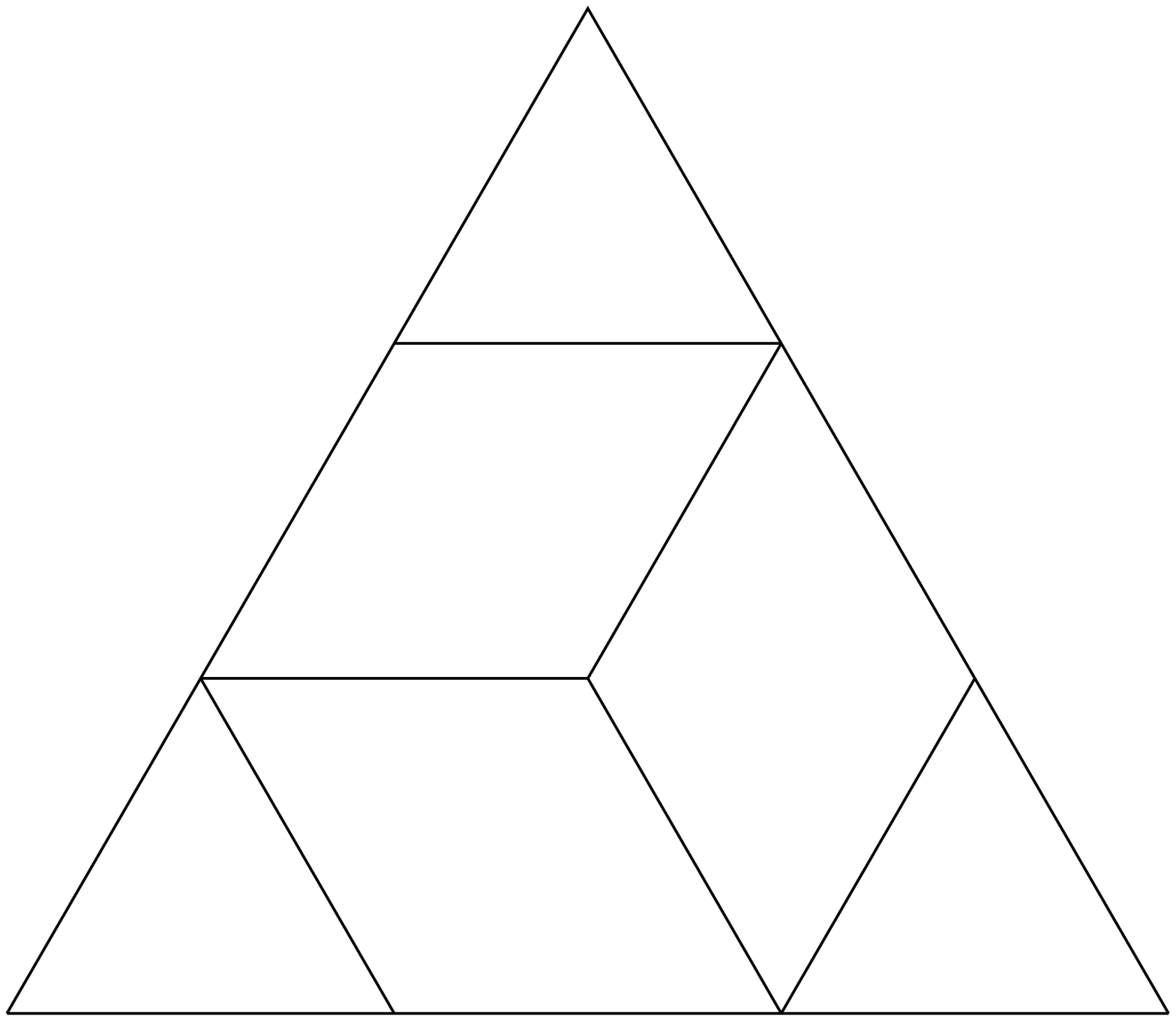}($\times 2$)
\caption{The configurations with 3 lozenges, side length $n=3$ \cite[Fig. 9.57]{SantosACM25}. 
They may cover 2/3 of an edge of the big triangle and one corner (multiplicity 6).
They may cover the three corners of the big triangle (multiplicity 1).
They may cover 1/3 of two edges of the big triangle and have the same orientation (multiplicity 3).
They may cover 2/3 of one edge of the big triangle and the middle of another (multiplicity 6).
They may touch in the middle and cover middle sections of all edges of the big triangle
(multiplicity 2, two circulations, \cite[Fig 3]{KenyonDM152}).
$L_{3,3}=6+1+3+6+2=18$.
}
\label{fig.33}
\end{figure}

The lozenge tilings generated from a big triangle
with side length $n=3$ 
are illustrated in Figures
\ref{fig.30}--\ref{fig.33}, sorted by the number of lozenges
$l=0,\ldots,3$. Some of the diagrams have multiplicities
larger than one if rotations by multiples of 120$^\circ$ or
flips across one of the three lines of symmetry of the big triangle generate
further diagrams of the same shape. (The isosceles triangle has a dihedral symmetry
group of order 12, where the 3 flips along a diagonal have order 2 and the rotations
by 120$^\circ$ or 240$^\circ$ have order 3. The multiplicity is 12 divided
by the order of the symmetry group once the lozenges are inserted.) The configurations generated
by these symmetry operations of the triangle are considered distinct
here; $L_{n,l}$ counts \emph{fixed}, not \emph{free} tilings.

\clearpage

\section{Special Cases}
\subsection{No Lozenge}

The formula
\begin{equation}
L_{n,0}=1
\end{equation}
means that for
each side length $n$ there is one way of not merging any triangles
into lozenges.

\subsection{One Lozenge}

The appearance of the triangular matchstick numbers 
\begin{eqnarray}
L_{n,1}=M_{n-1} \label{eq.Ln1}, \\
\sum_{n\ge 0}L_{n,1}x^n = \frac{3x^2}{(1-x)^3},
\end{eqnarray}
is obvious:
recall that deleting one of the internal edges
creates a lozenge by merging the two triangles that share that edge, 
so \eqref{eq.Ln1} just restates \eqref{eq.inte}. One may refine that total
by classifying the lozenges by the three sets of edges 
that represent the short diagonals, regarding the edges as
arguments to the counts of the structures that can be placed inside the big triangle:
\begin{equation}
V(\includegraphics[scale=0.05]{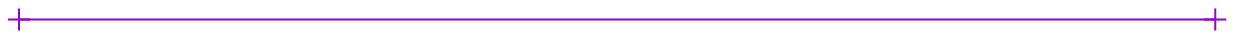})
=
V(\includegraphics[scale=0.02]{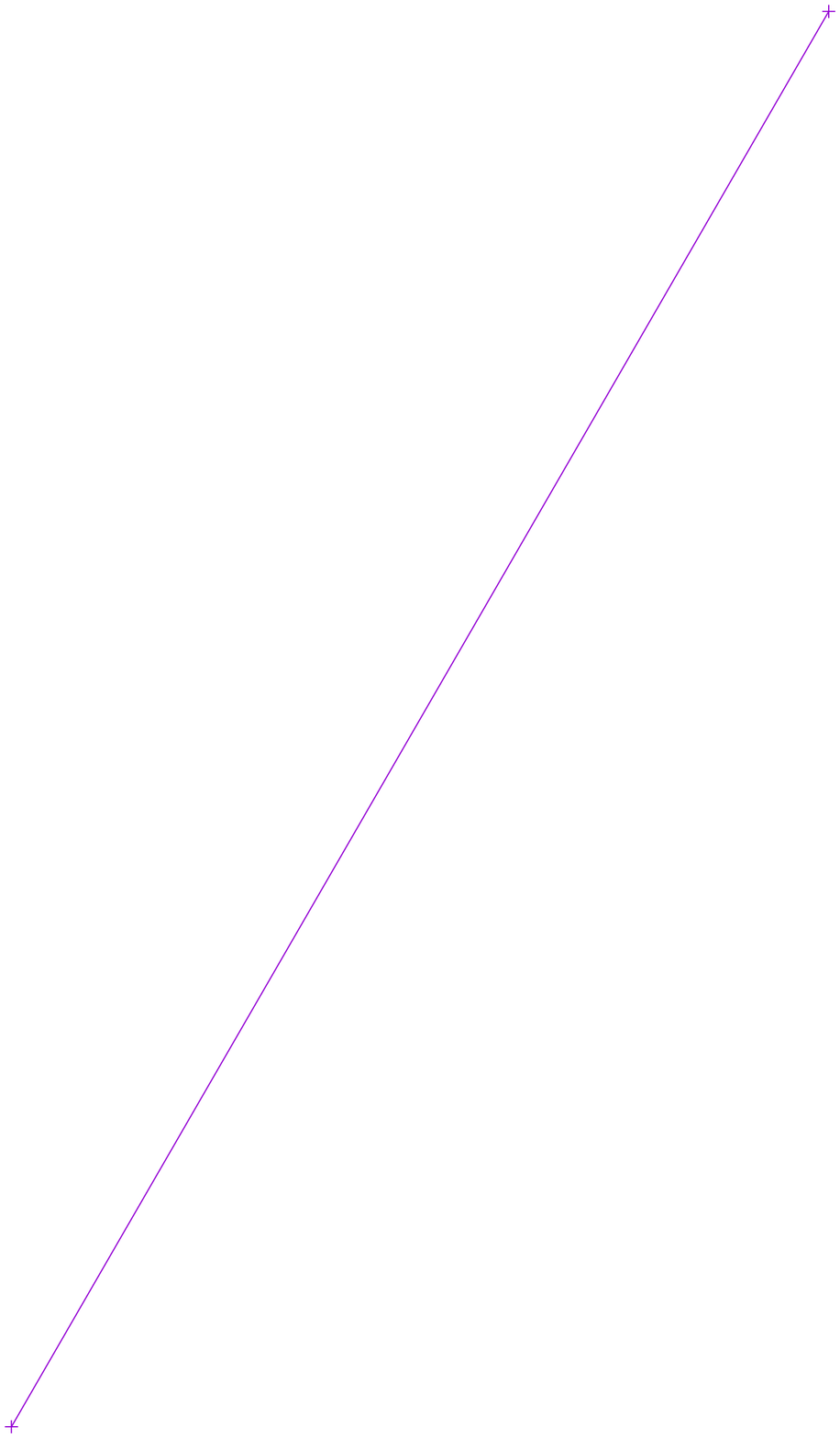})
=
V(\includegraphics[scale=0.02]{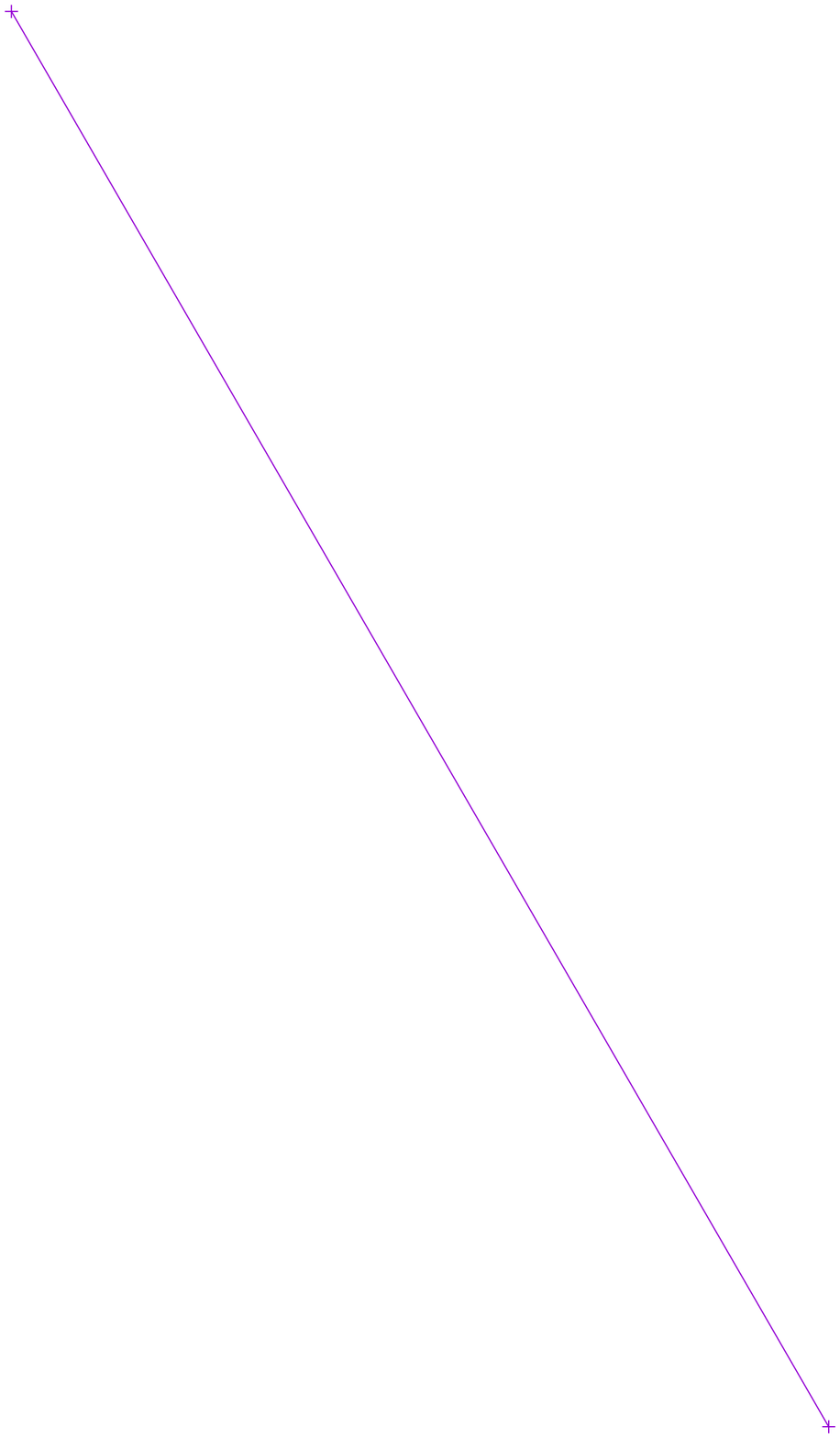})
=T_{n-1}.
\end{equation}

\subsection{Two Lozenges}\label{sec.2L}
Two lozenges are created by deleting two internal edges, which 
can be selected in $\binom{M_{n-1}}{2}$ ways.
Some of these pairs of deleted edges do not represent lozenge tilings because
they are spatially correlated as outlined in Section \ref{sec.corr}.
\begin{defn}
A $V$ subgraph is a pair of internal edges (in the full graph without
lozenges) that share one common vertex,
where the two edge directions differ
by an angle of $60^\circ$.
\end{defn}

There are two distinct sets of $V$'s: the geometries where the $V$ points
to a corner of the big triangle:
\begin{equation}
V(\includegraphics[scale=0.02]{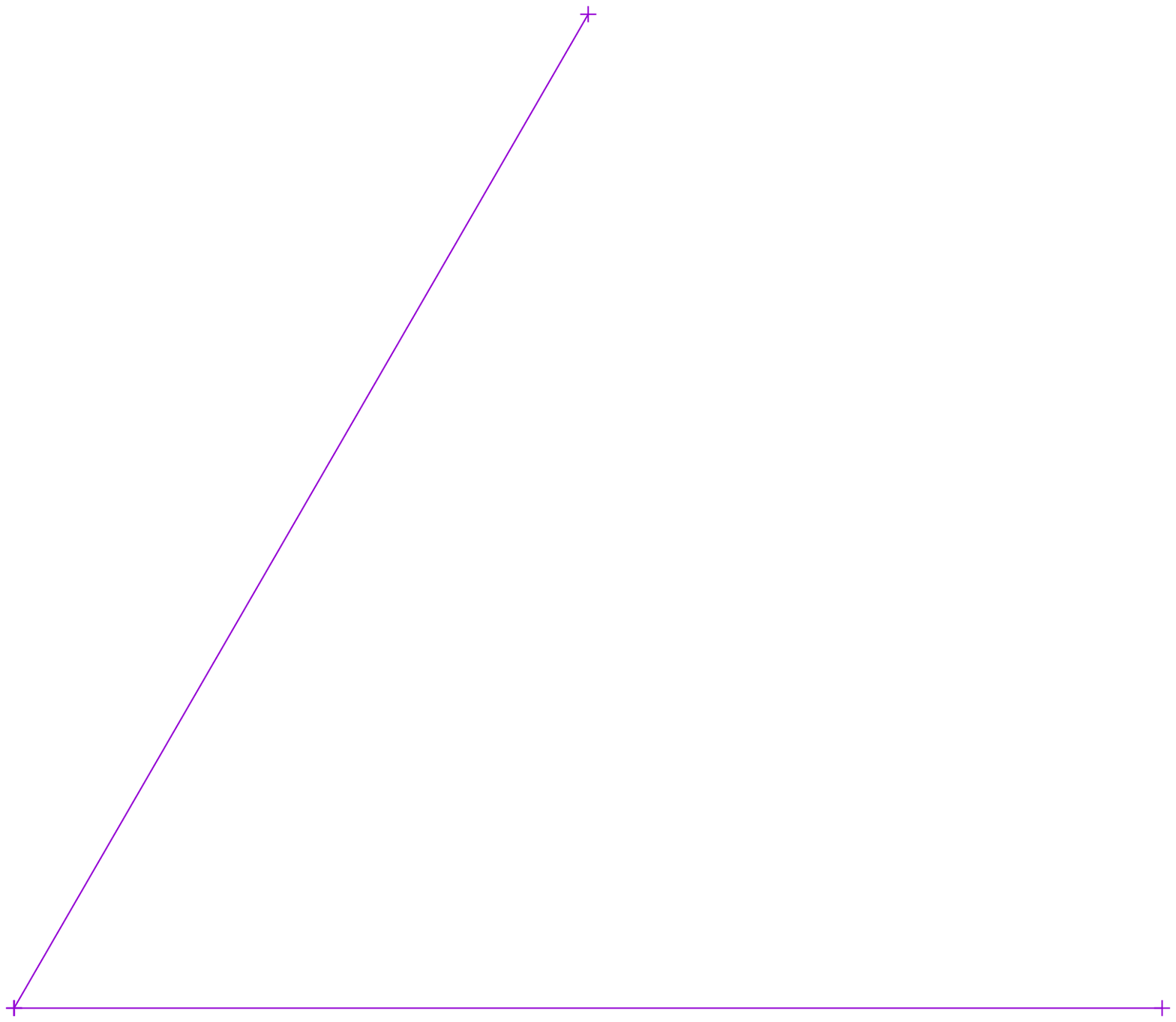})
=
V(\includegraphics[scale=0.02]{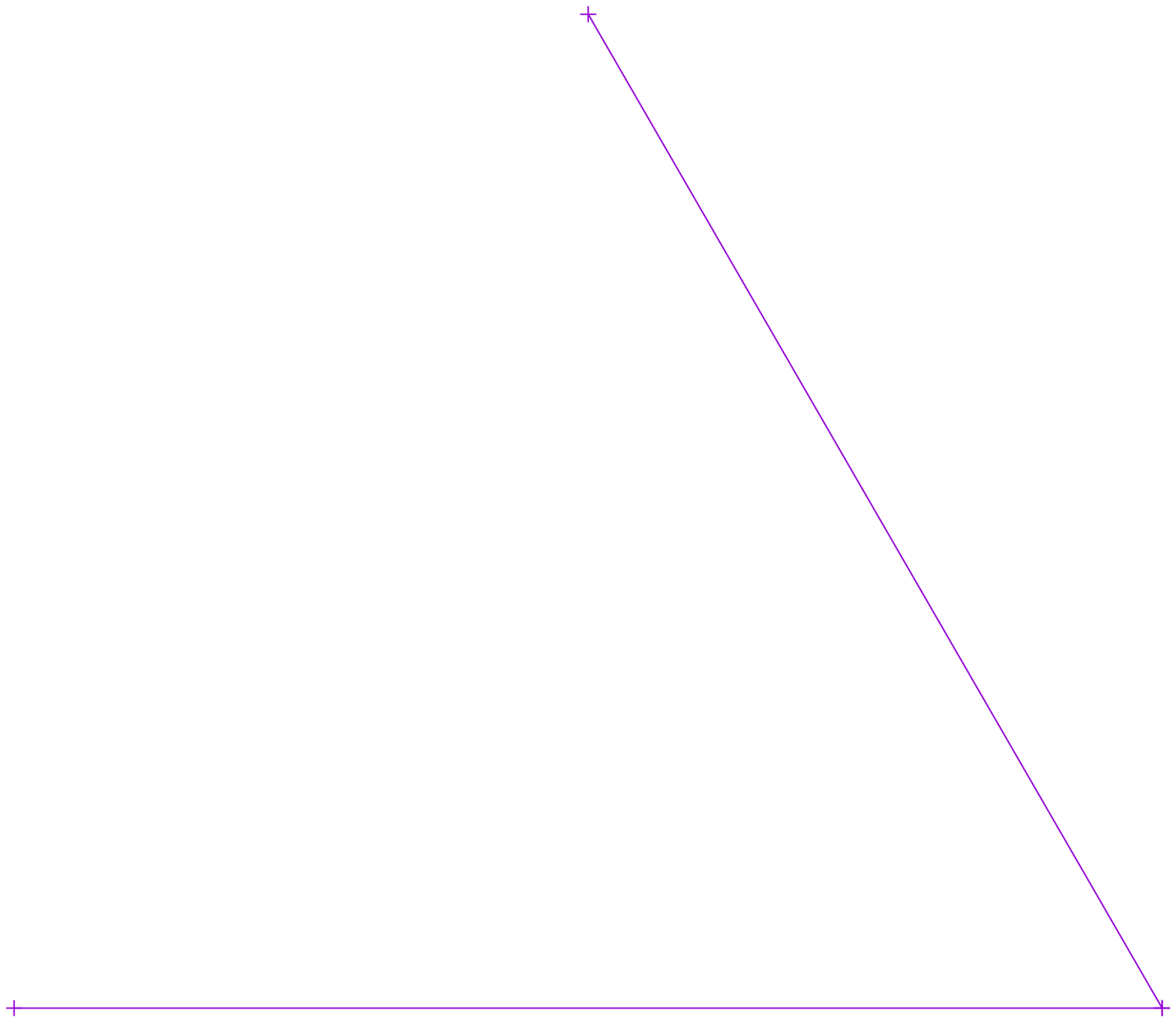})
=
V(\includegraphics[scale=0.02]{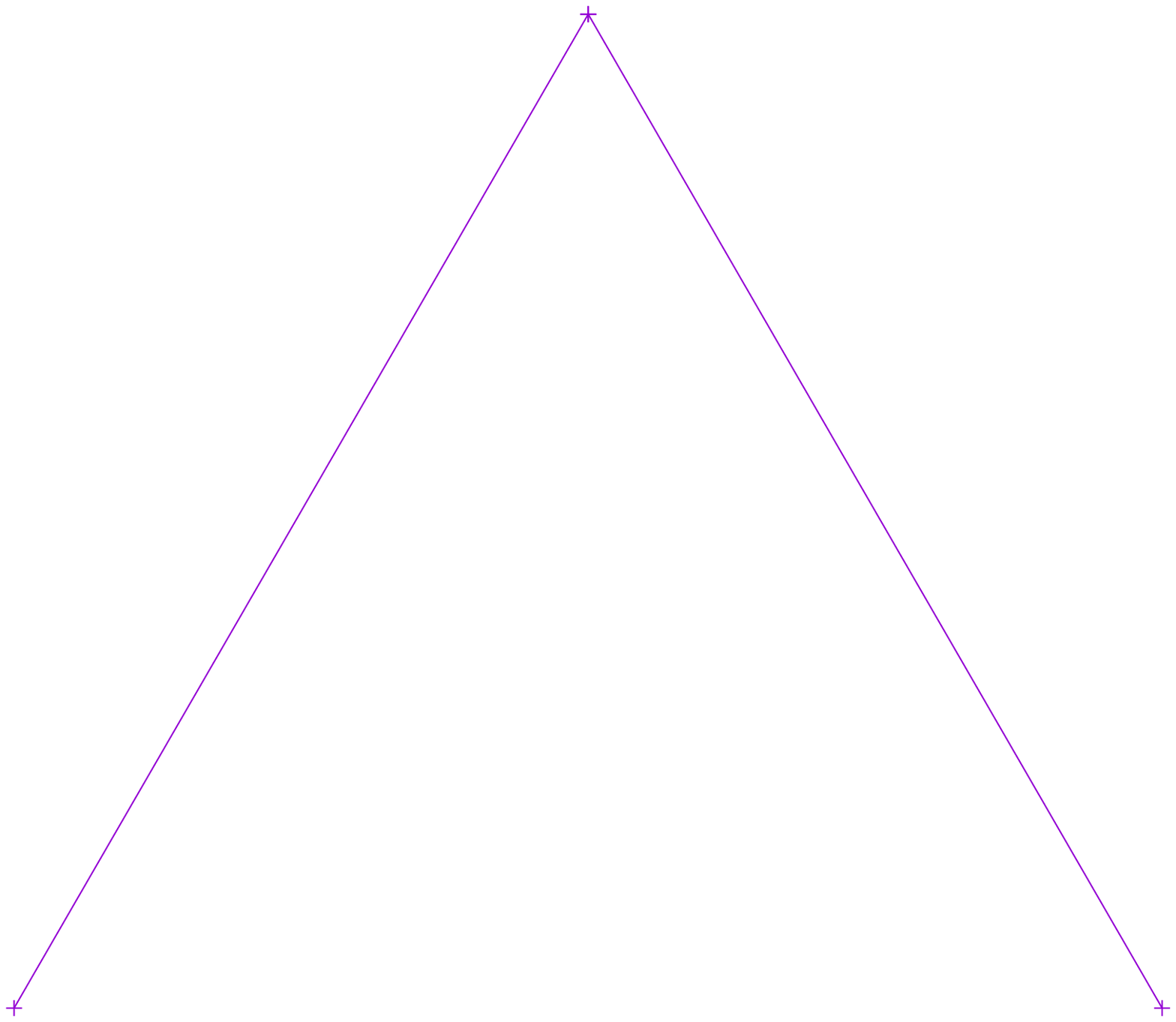})
=T_{n-2},
\label{eq.pl2first}
\end{equation}
and where they point to an edge:
\begin{equation}
V(\includegraphics[scale=0.02]{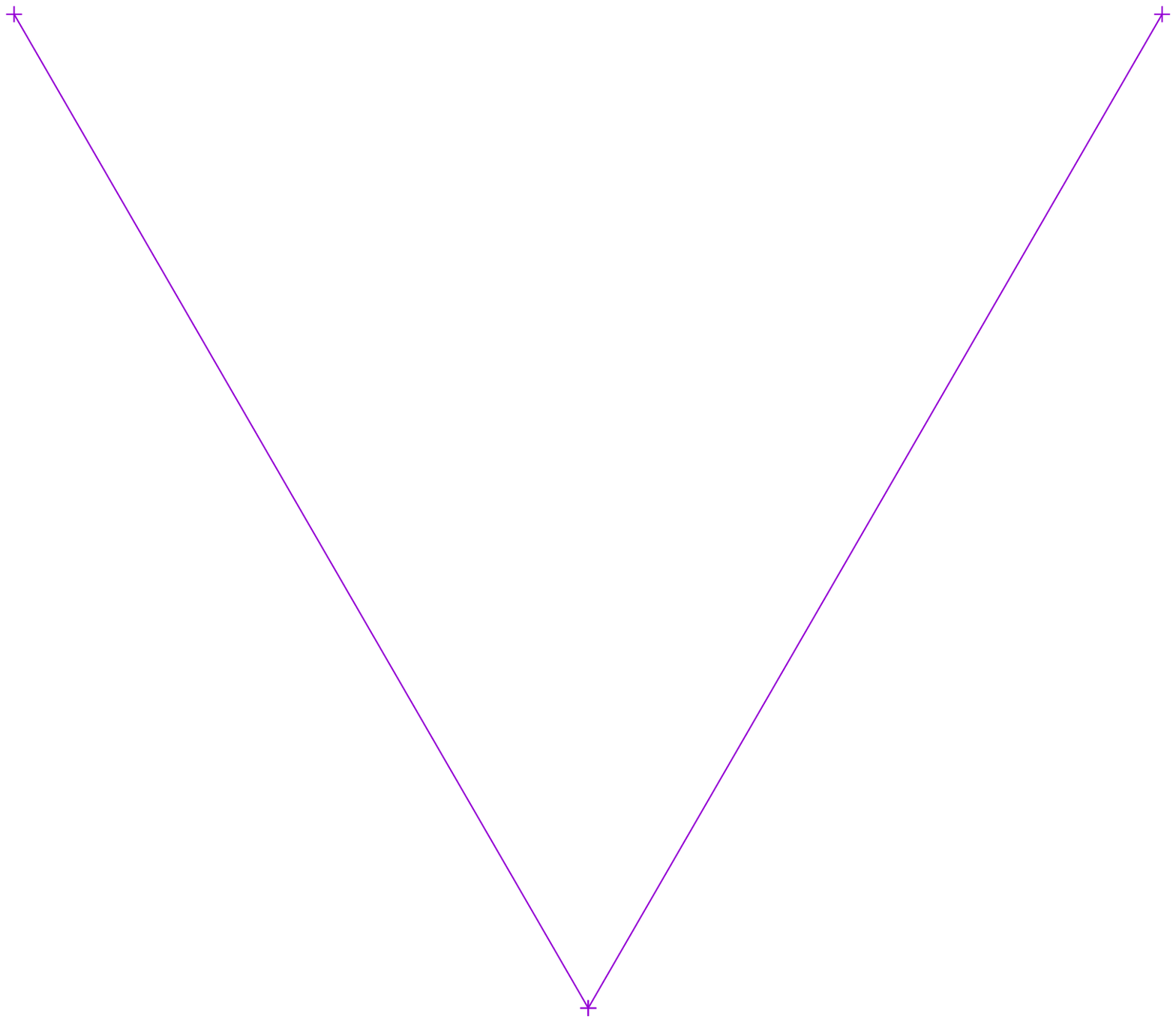})
=
V(\includegraphics[scale=0.02]{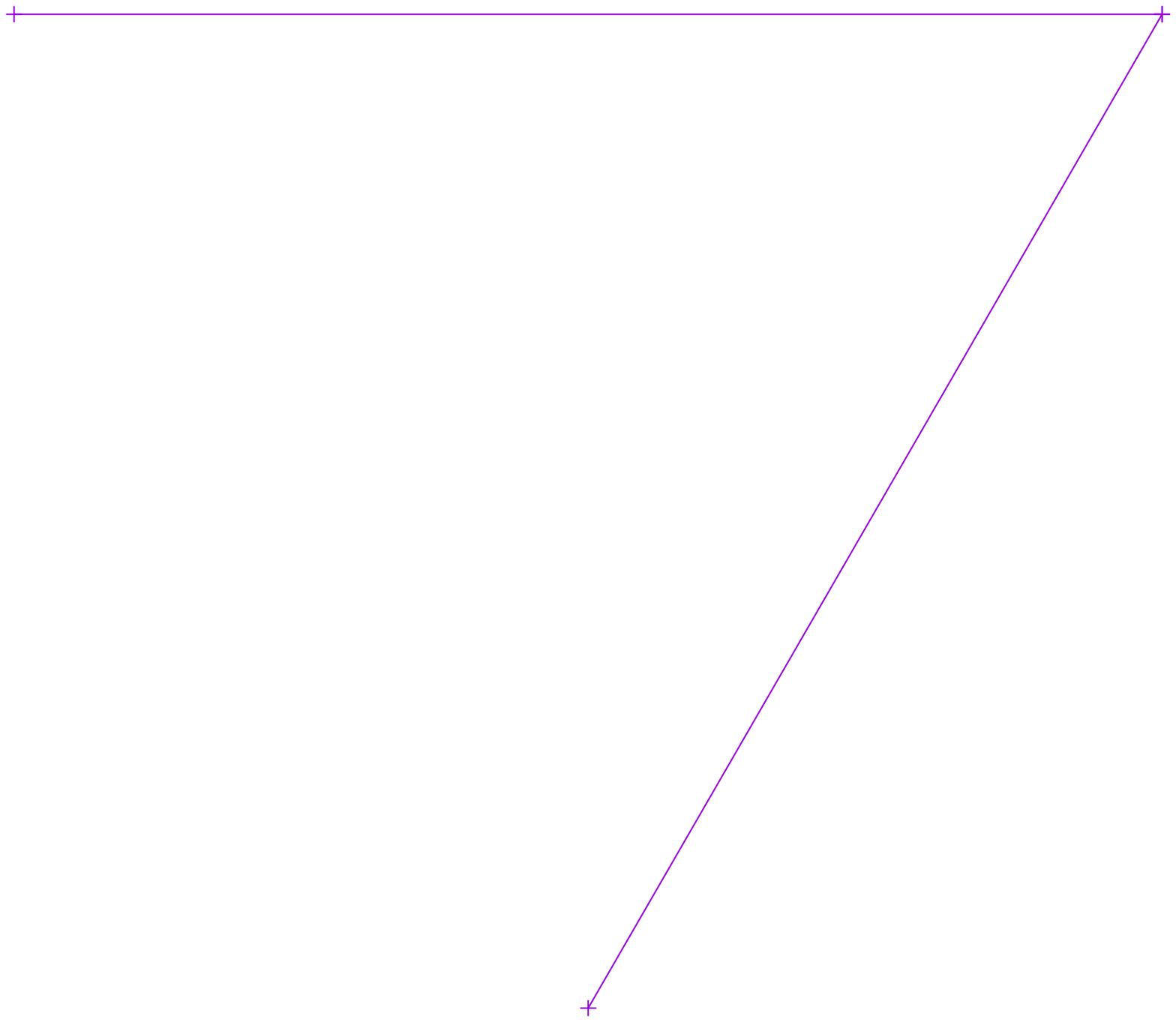})
=
V(\includegraphics[scale=0.02]{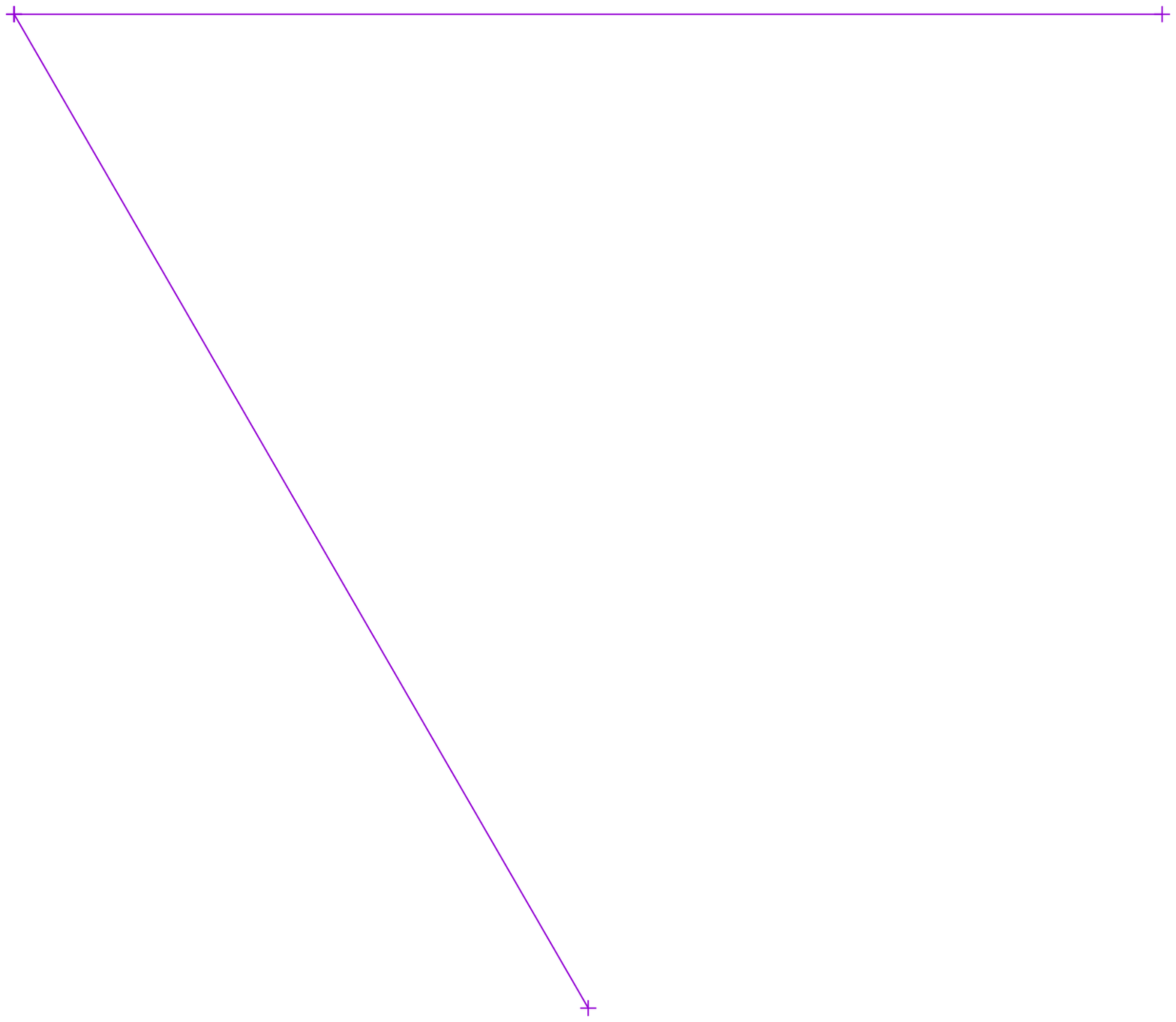})
=T_{n-1}.
\label{eq.pl2last}
\end{equation}
The number places accomodated in the big triangle for some fixed edge set is
always a triangular number, because one can push it as far as possible into
the lower left corner, translate it by units to the left until it touches
the lower right corner (in a number of ways which is $n$ minus a constant),
repeat the procedure one layer higher up (which gives a count that is one less),
until the top corner is reached. The sum of these counts is obviously a triangular
number.
The difference in the counts for the different orientations of the
$V$ between \eqref{eq.pl2first} and \eqref{eq.pl2last} is caused by the fact that
the big triangle is fixed and has only a rotational symmetry axis of order 3, whereas
the (infinite) hexagonal grid has a rotational symmetry axis of order 6---a 
symmetry breaking boundary effect which
will also be observed in the counts later on in this work.

The total number of $V$ subgraphs is 
\begin{equation}
V(\includegraphics[scale=0.02]{pl2_1.eps})
+
V(\includegraphics[scale=0.02]{pl2_3.eps})
+
V(\includegraphics[scale=0.02]{pl2_5.eps})
+
V(\includegraphics[scale=0.02]{pl2_2.eps})
+
V(\includegraphics[scale=0.02]{pl2_4.eps})
+
V(\includegraphics[scale=0.02]{pl2_6.eps})
=3(n-1)^2.
\end{equation}
Each $V$ subgraph reduces the number of lozenge 
tilings by one. 
The number of configurations with 2 lozenges becomes \cite[A326367]{sloane}
\begin{equation}
L_{n,2} 
= \binom{M_{n-1}}{2}-3(n-1)^2
= \frac38(n-1)(n-2)(3n^2+3n-4),\quad n\ge 1 . \\
\end{equation}
The (inverse) binomial transform is
\begin{equation}
L_{n,2} 
= -3+3\binom{n}{1}-3\binom{n}{2}+27\binom{n}{3}+27\binom{n}{4},\quad n\ge 1 . \\
\end{equation}

\subsection{Polyedges}

\begin{figure}
\includegraphics[scale=0.5]{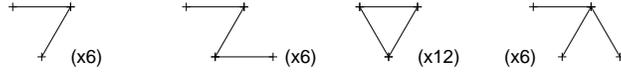}
\caption{The 1 V-shaped forbidden free polyedge in the triangular grid with 2 edges and the 3 forbidden free polyedges 
(zigzag, triangle, fork) with 3 edges.}
\label{fig.poly2}
\end{figure}

\begin{figure}
\includegraphics[scale=0.5]{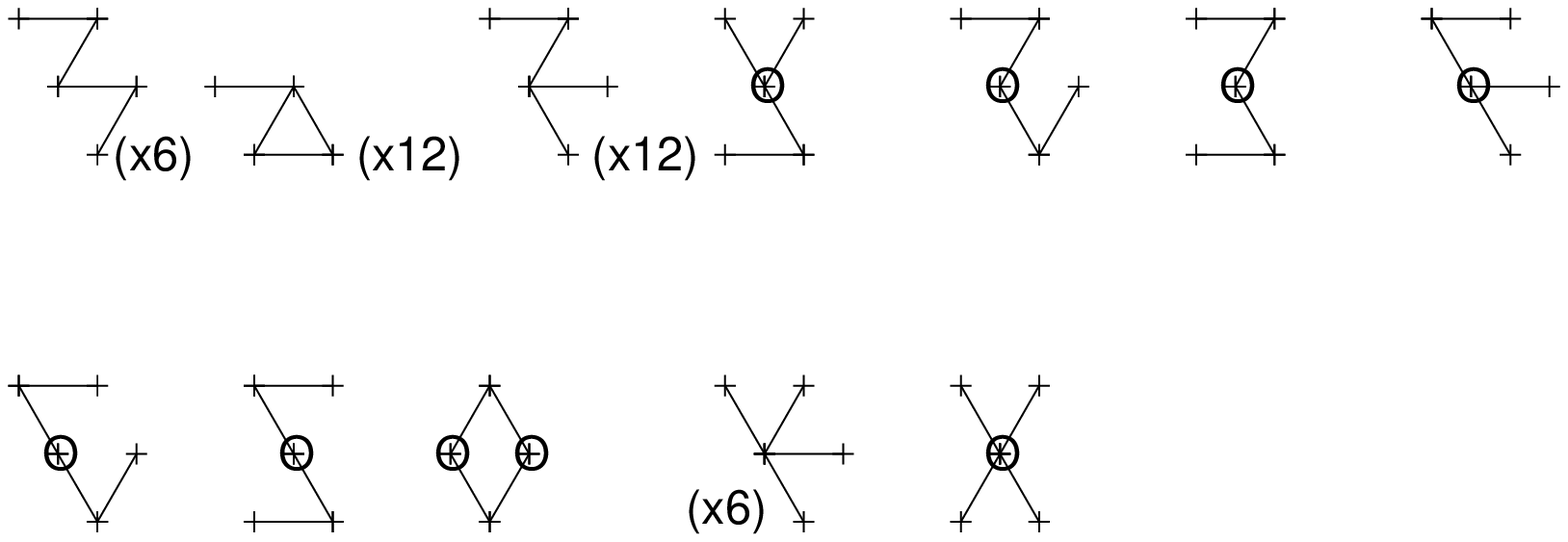}
\caption{The 12 forbidden free polyedges in the triangular grid with 4 edges. 
}
\label{fig.poly4}
\end{figure}

\begin{figure}
\includegraphics[scale=0.5]{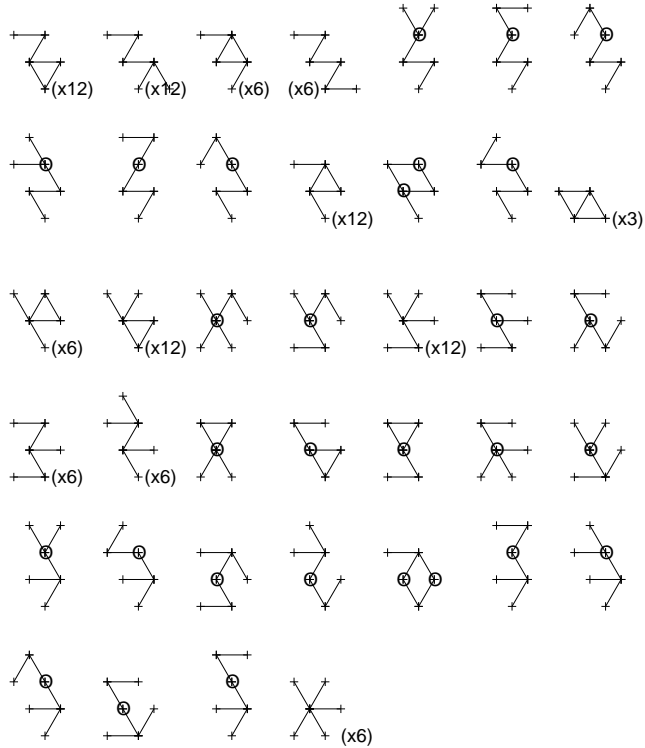}
\caption{The 39 forbidden free polyedges in the triangular grid with 5 edges.}
\label{fig.poly5}
\end{figure}

There exist 1, 3, 12, 60, 375,\ldots connected free polyedges
in the triangular grid with 1, 2, 3,\ldots edges \cite[A159867]{sloane}.
This work focuses on polyedges each edge of which has at least one adjacent edge
that meets at an angle of 60$^\circ$; that relevant subset of ``forbidden''
free connected polyedges contains 
1, 3, 12, 39, 209, 1014,\ldots edges with 2, 3,\ldots edges, illustrated
in Figures \ref{fig.poly2}--\ref{fig.poly5} for up to 5 edges. 
Figure \ref{fig.poly5} will not be employed here; it is  a preview on a possible
extension of our technique to proof the polynomial formula for $L_{n,5}$.

Circles around vertices indicate
that the graph can be cut  at these vertices into two smaller polyedges without missing
forbidden edge sets. From the point of view of the algebra further down, these
polyedges are coincidences of two smaller polyedges that happen to have one or two vertices in common.

\subsection{Three Lozenges}\label{sec.3l}
Three lozenges are created by deleting three internal edges, which 
can be selected in 
$\binom{M_{n-1}}{3}$ 
ways according to \eqref{eq.bino}.
Some of these triples of deleted edges do not represent lozenge tilings because
they are spatially correlated as defined in Section \ref{sec.corr}.

From the 3 free polyedges with 3 edges in Figure \ref{fig.poly2} 
we derive 14 fixed types by rotations
and flips, 6 zigzags, 2 triangles and 6 forks:
\begin{equation}
V(\includegraphics[scale=0.02]{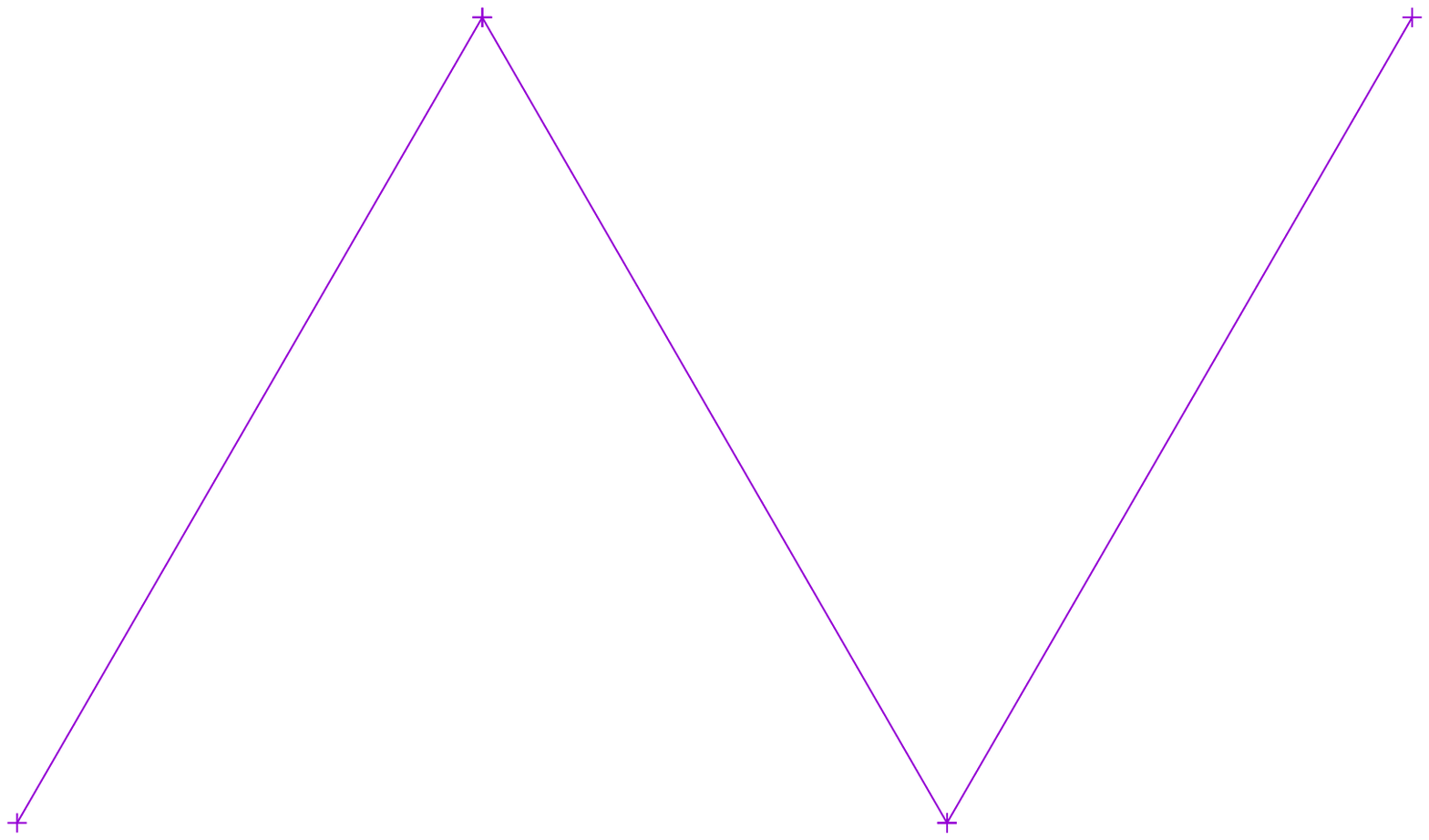})
=
V(\includegraphics[scale=0.018]{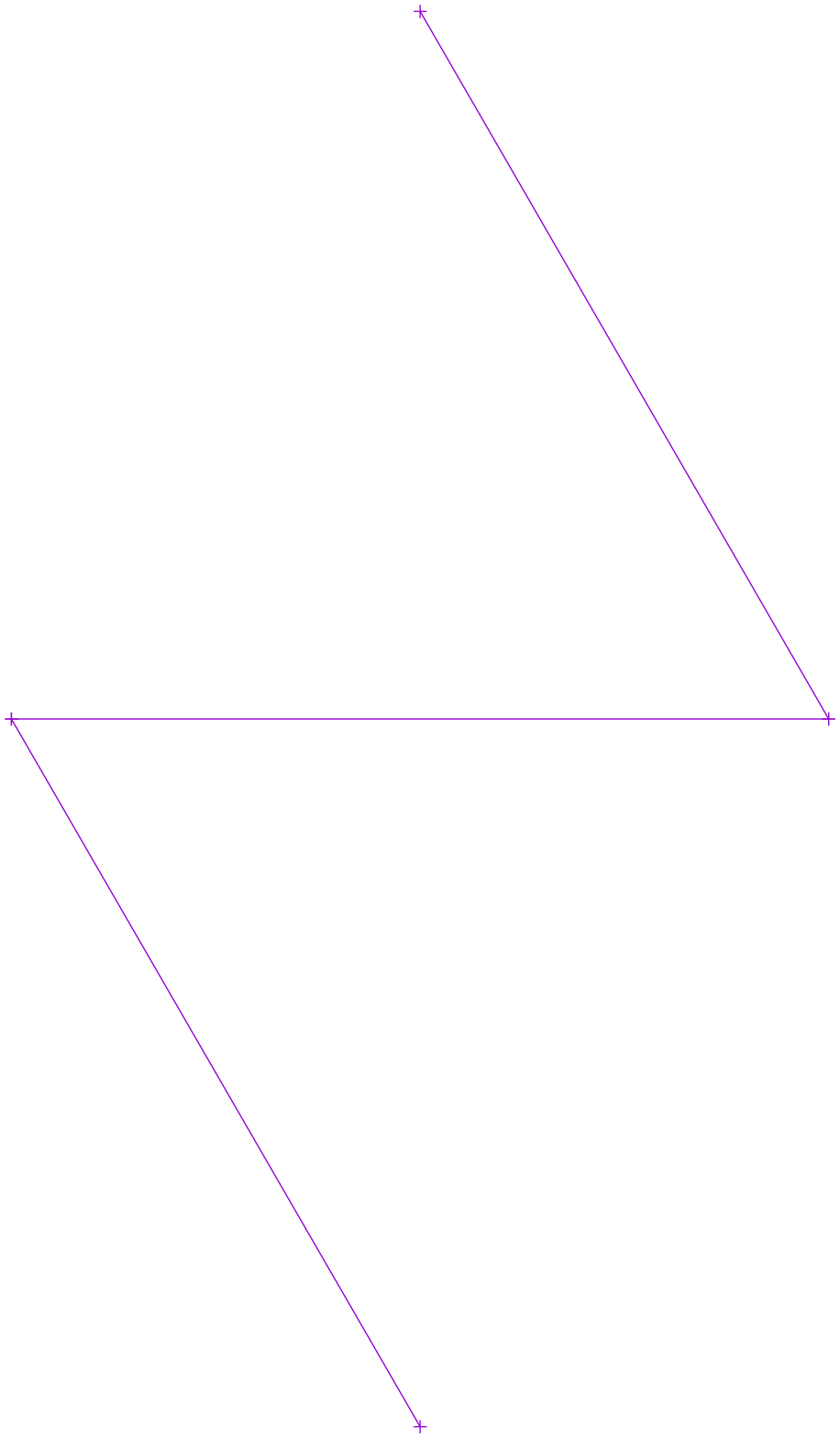})
=
V(\includegraphics[scale=0.02]{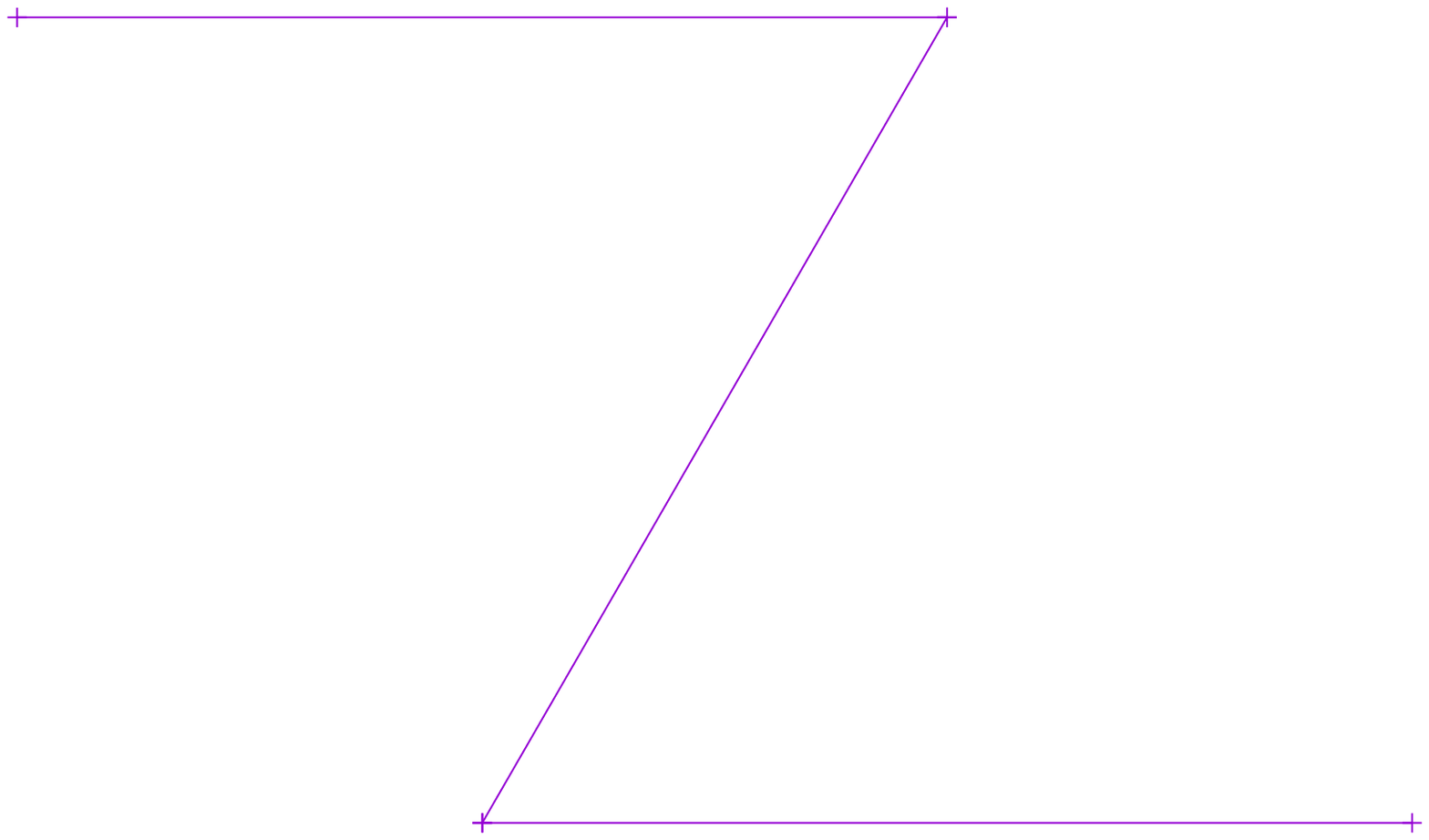})
=
V(\includegraphics[scale=0.02]{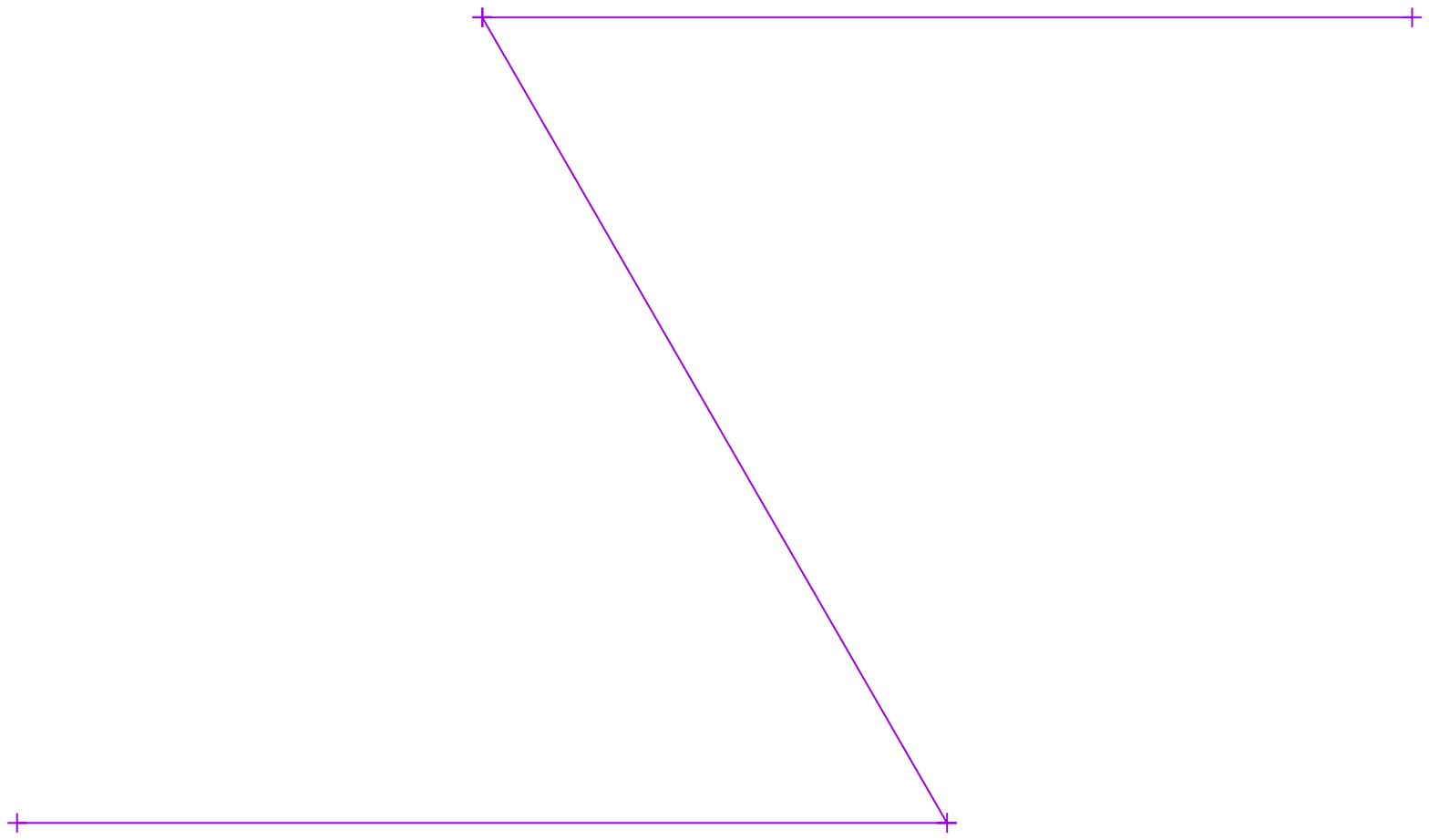})
=
V(\includegraphics[scale=0.018]{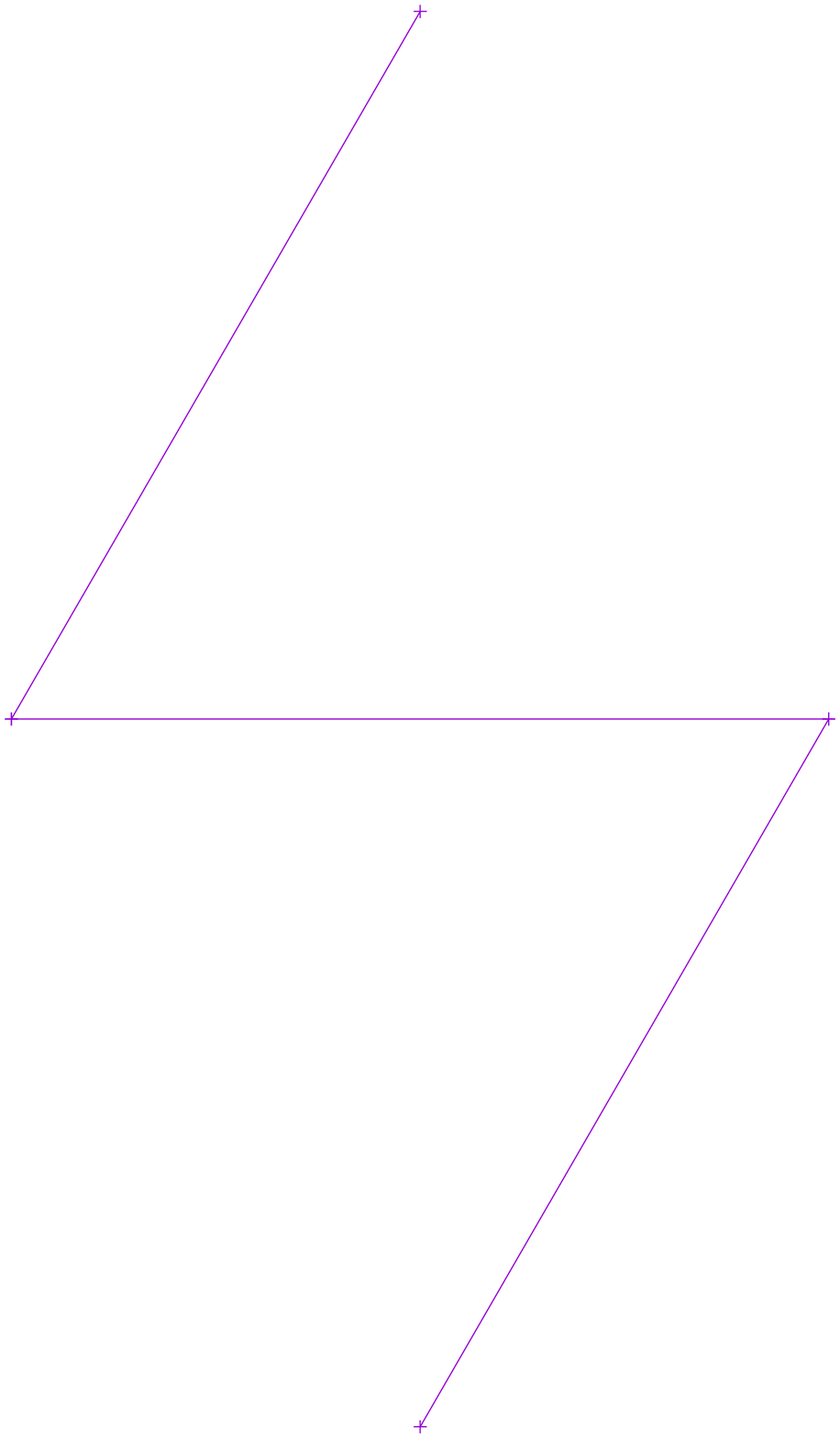})
=
V(\includegraphics[scale=0.02]{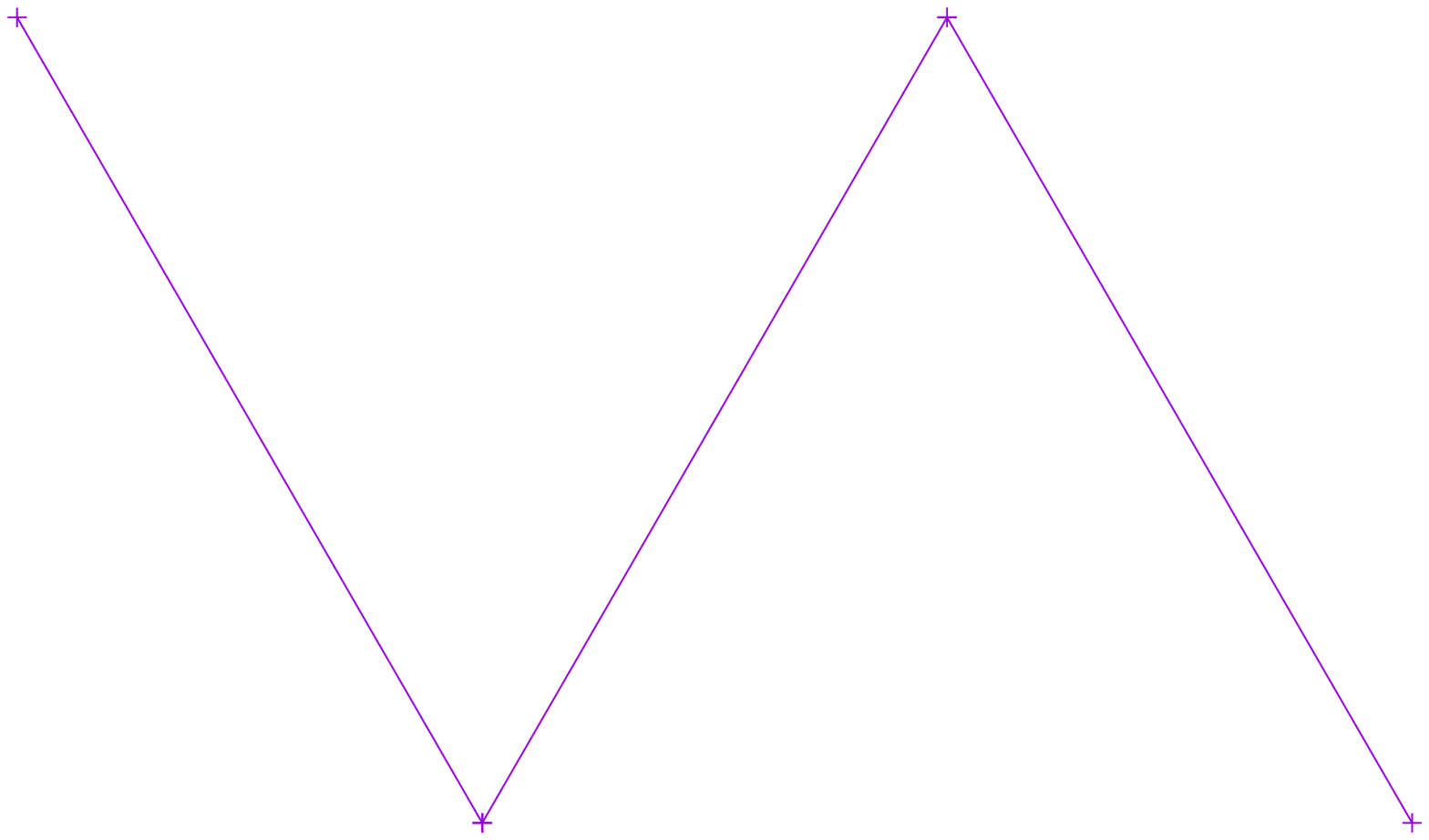})
=T_{n-2};
\label{eq.pl3first}
\end{equation}
\begin{equation}
V(\includegraphics[scale=0.02]{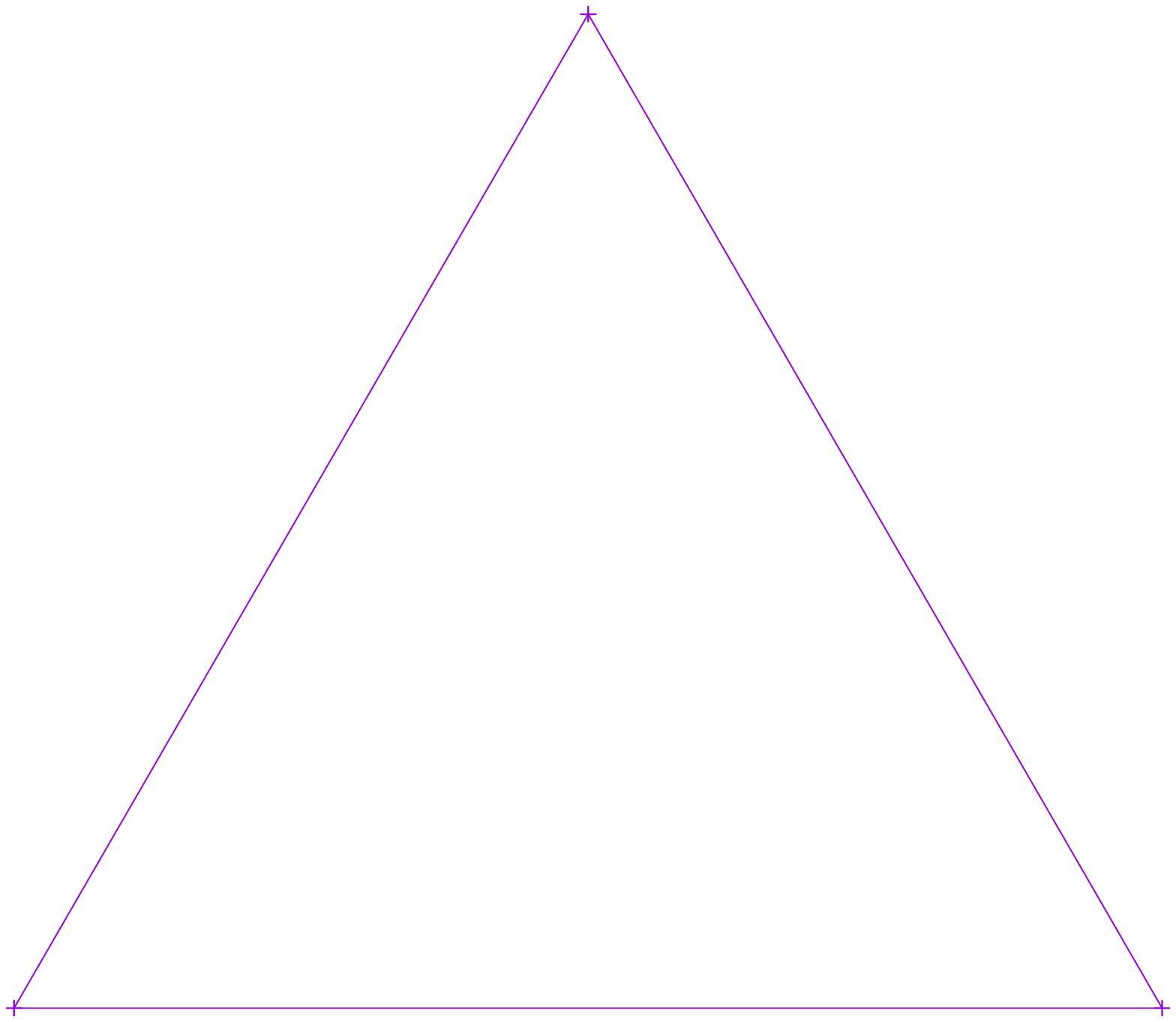})
=T_{n-3}; \quad
V(\includegraphics[scale=0.02]{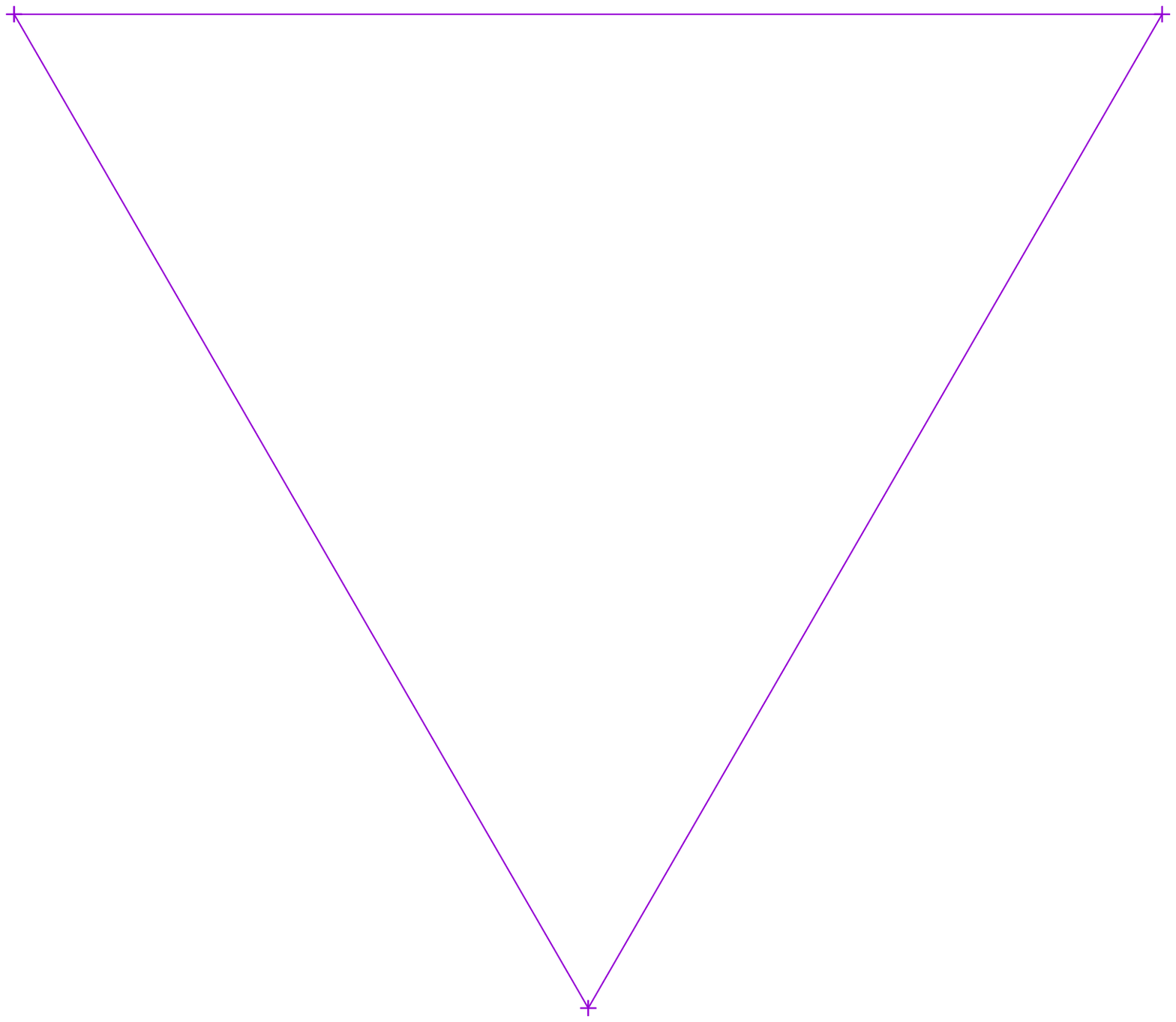})
=T_{n-1};
\end{equation}
\begin{equation}
V(\includegraphics[scale=0.02]{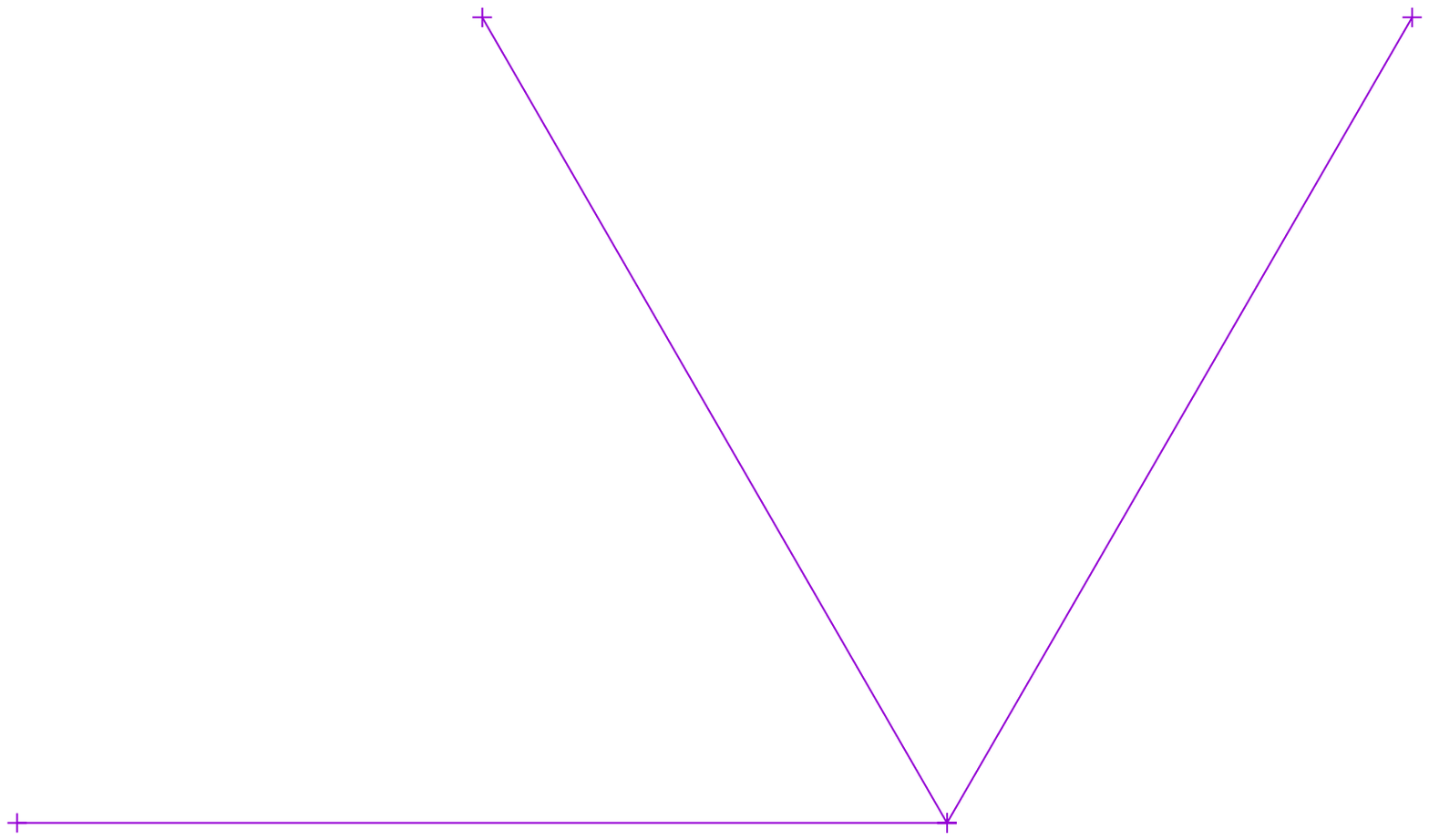})
=
V(\includegraphics[scale=0.018]{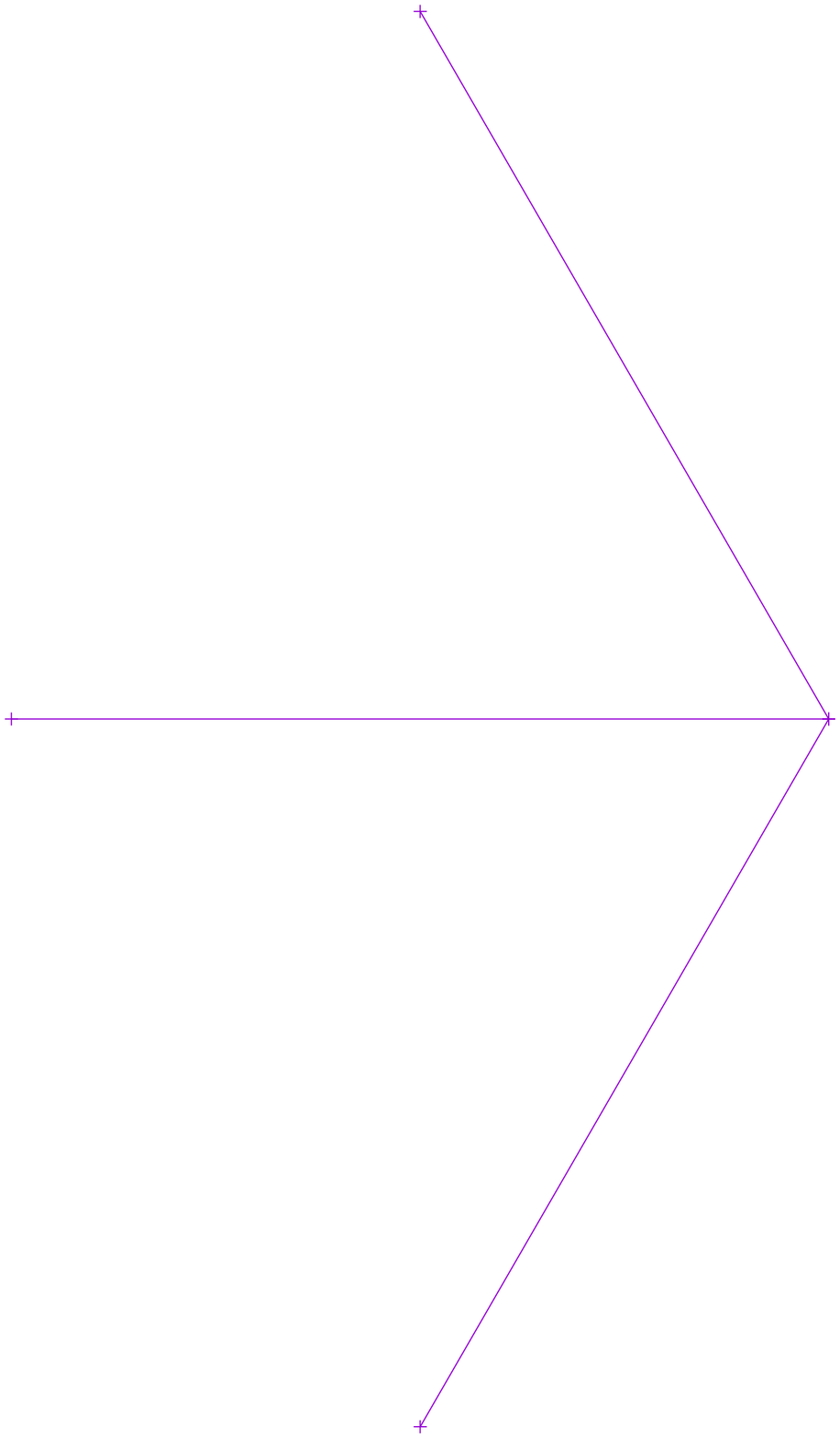})
=
V(\includegraphics[scale=0.02]{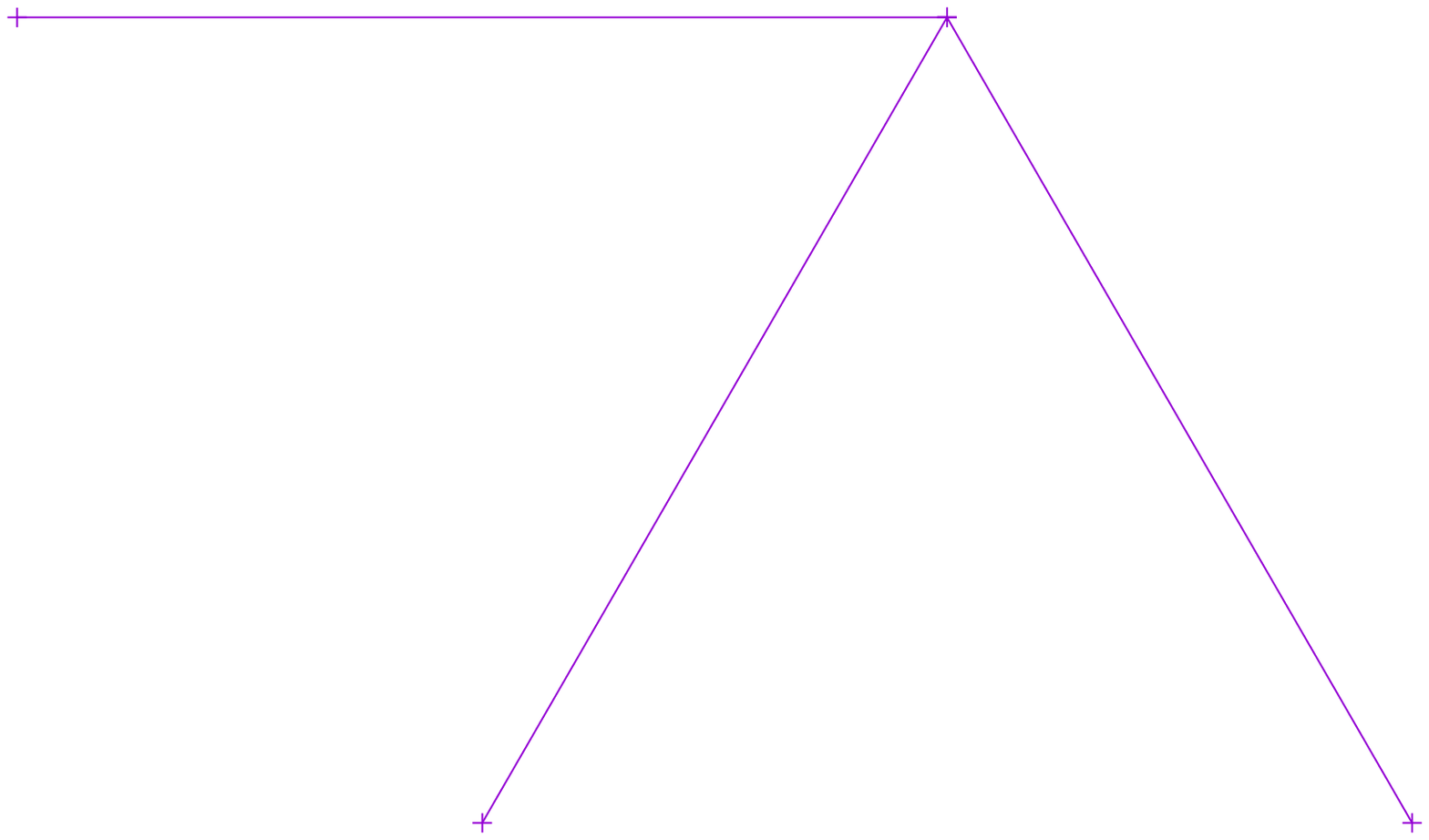})
=
V(\includegraphics[scale=0.02]{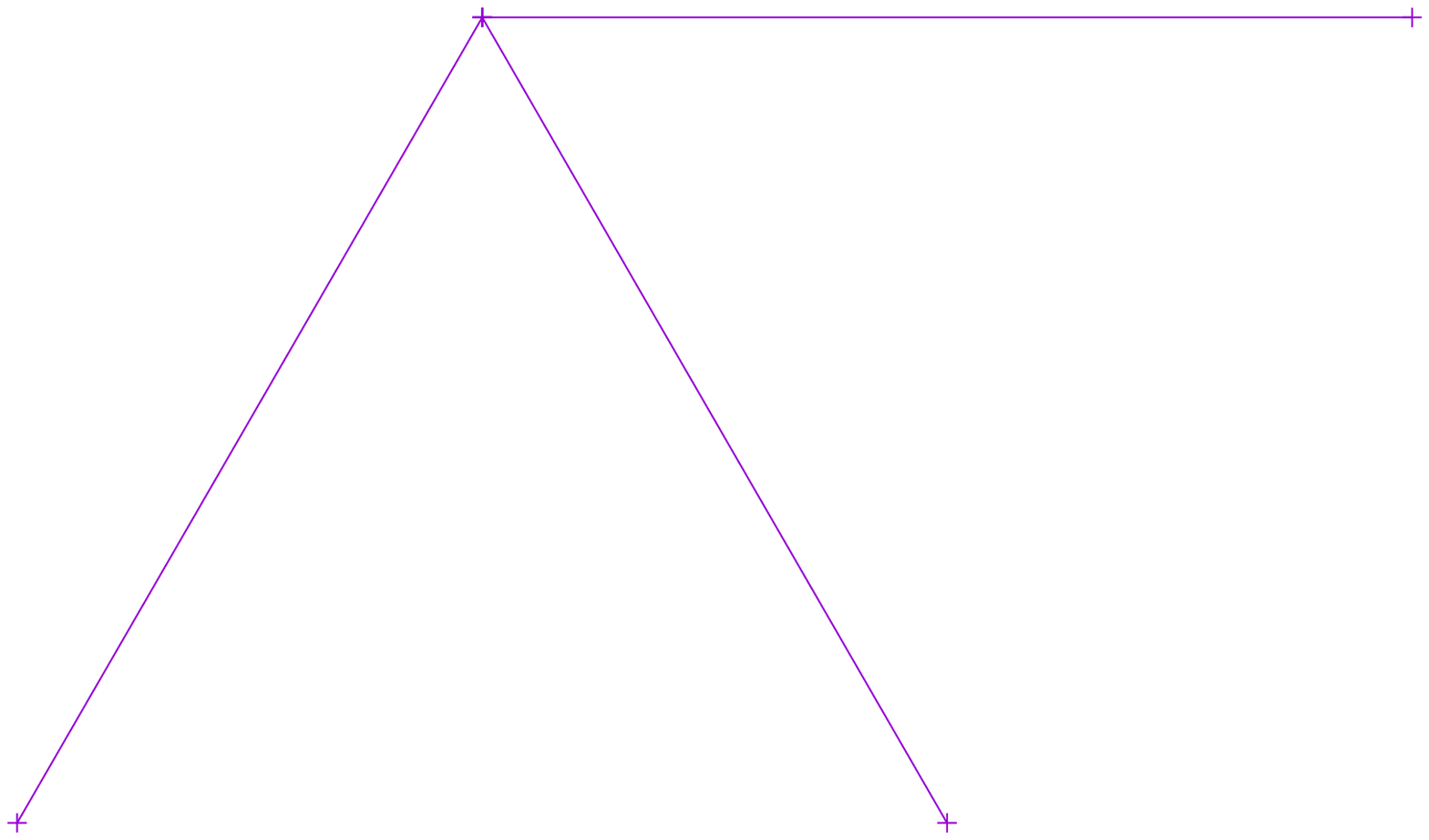})
=
V(\includegraphics[scale=0.018]{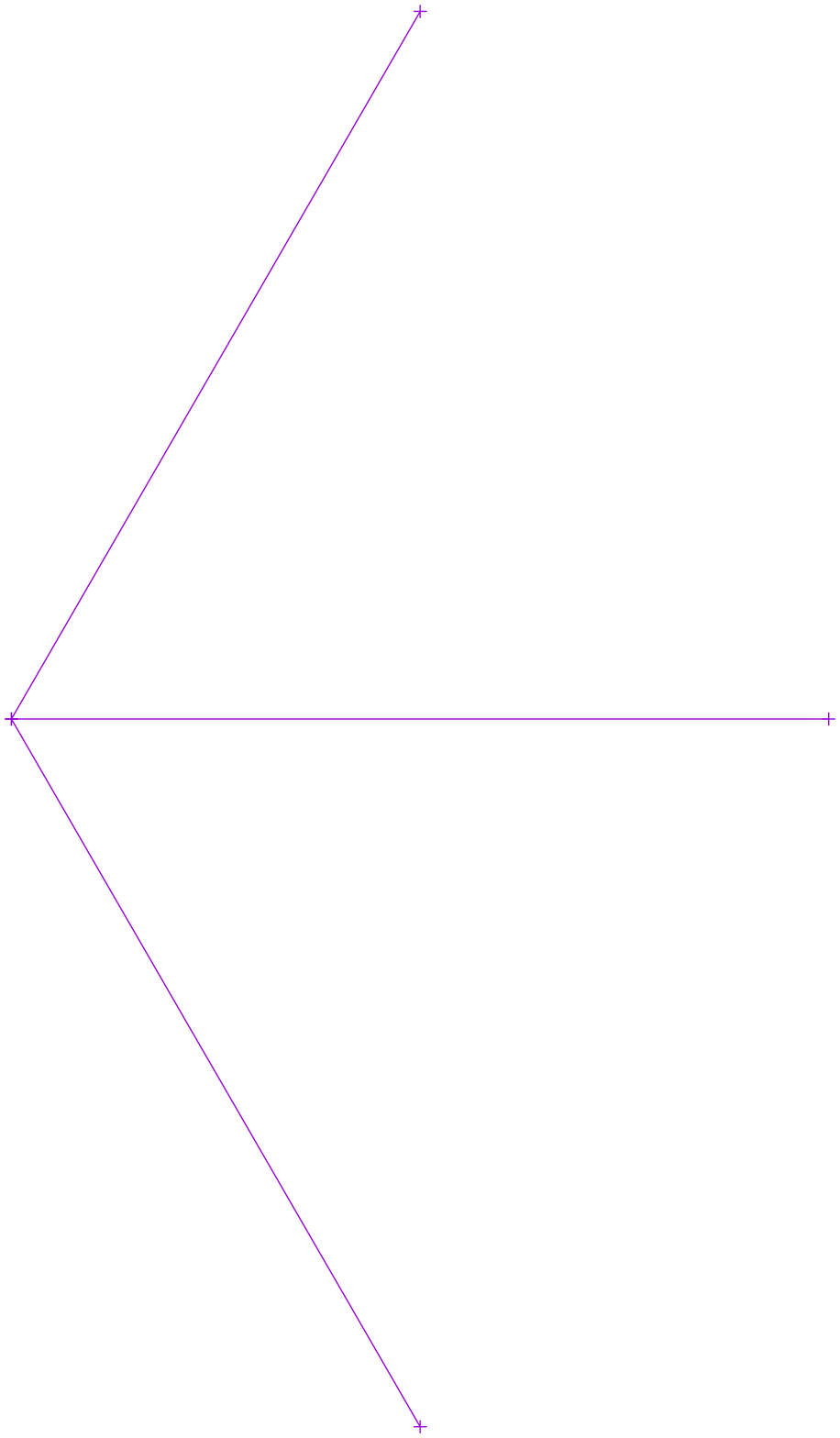})
=
V(\includegraphics[scale=0.02]{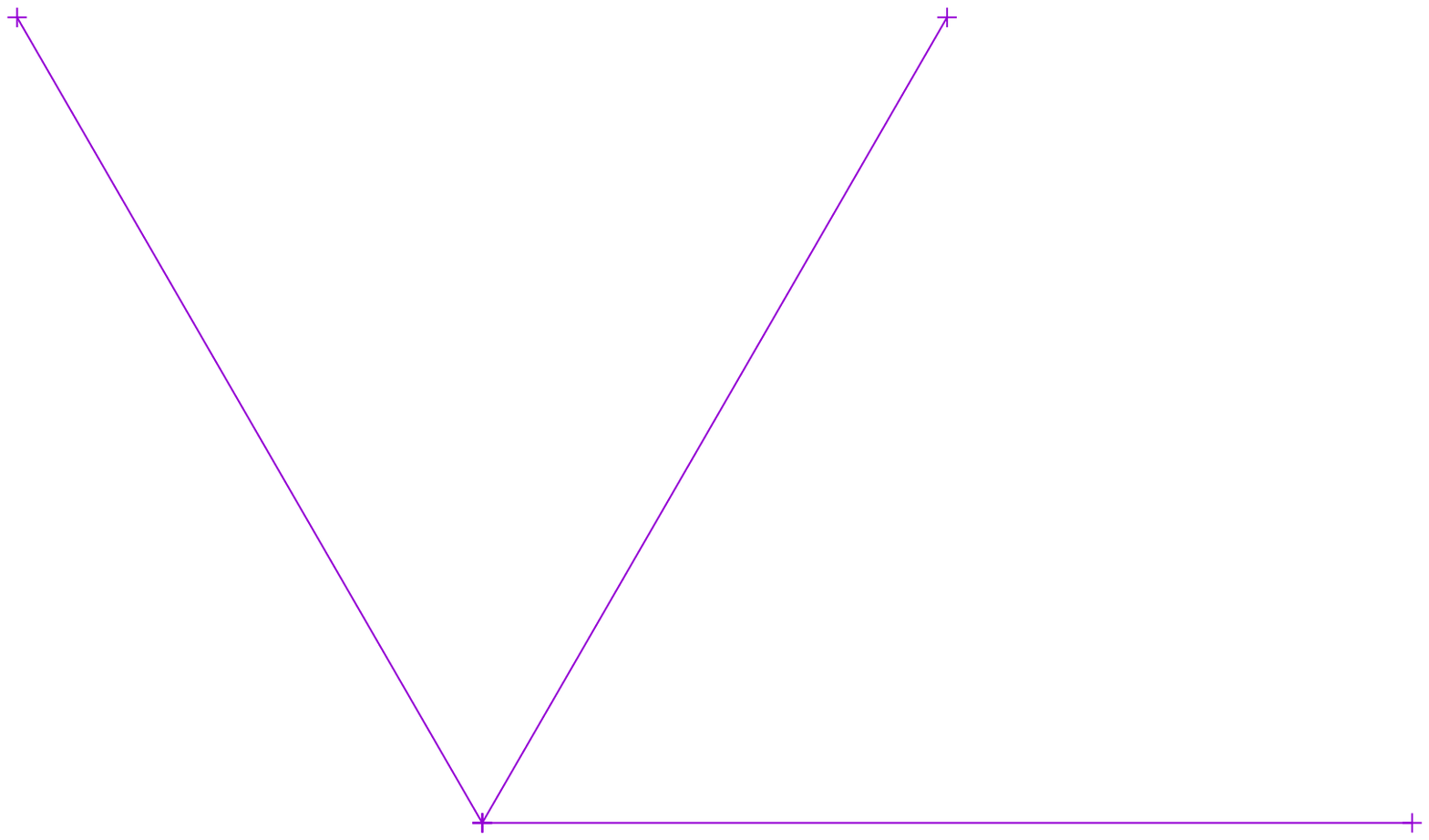})
=T_{n-2}.
\label{eq.pl3last}
\end{equation}
The task is to compute how many copies of the enumerations
in equations \eqref{eq.pl2first}--\eqref{eq.pl2last} and
\eqref{eq.pl3first}--\eqref{eq.pl3last} are to be considered ``forbidden'' 
polyedges in the big triangle.

\begin{figure}
\includegraphics[scale=0.8]{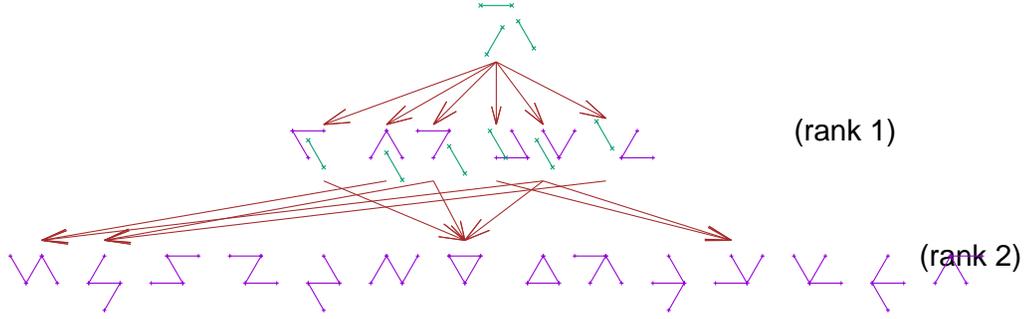}
\caption{Hasse diagram of the deleted edges for $l=3$ deleted unit edges in general, and specializations
where two of
them are correlated in a $V$ graph, and where three of them are correlated in a triangle,
a zigzag or a fork subgraph.}
\label{fig.pos3}
\end{figure}

Polyedges with 3 edges can be classified as (i) polyedges with 3 uncorrelated edges,
(ii) sets of one connected polyedge with 2 edges and another uncorrelated edge, and
(iii) connected polyedges with 3 edges---the standard multiset argument.
This yields Figure \ref{fig.pos3} with three levels, containing 1 set (rank 0), 6 sets (rank 1)
and 14 sets (rank 2)
of polyedges. Green edges symbolize randomly chosen internal edges which do not meet
edges of the 
other (blue) polyedge-substructure
at angles of 60$^\circ$.
Specializing one of the randomly chosen internal edges so the joints at
60$^\circ$ increase in number means increase the rank by 1,
in the language of posets the diagrams at lower rank up are $\ge$ than the 
diagrams at higher rank, and diagrams at the same rank are not comparable.
The 6 brown arrows from the rank 0  to the rank 1
sets indicate how $V$-graphs are generated by correlating
two edges. The 9 brown arrows from the rank 1 to the rank 2 sets are an
(incomplete) illustration how any of the 6 zigzag  or 6 fork graphs is 
a specialized version of two $V$-graphs, and of how the 2 triangular graphs
are specialized versions of three $V$-graphs.
The strategy is to enumerate the fixed polyedges at ranks $\ge 1$ 
of the Hasse diagram.

Actually the enumeration of the diagrams of rank 2 is already
completed with equations \eqref{eq.pl3first}--\eqref{eq.pl3last}.
The enumeration of the polyedges at rank 1 is completed here
by Mobius inversion of the enumerations with unions of the sets of rank 1 and 2
\cite{RotaZWV2}.
The functions $\bar V$ count sets of complementary polyedges which have a specific
$V$-graph and a third edge selected from any other of the remaining $M_{n-1}-2$
inner edges of the triangle. The distinction with $V$ is that $\bar V$ does not
care whether the third edge creates a further 60$^\circ$ angle with any of the
two edges of the $V$. Regard the overbar as a closure or superset counting
function. The enumeration is the product of the enumerations of 
placing the associate $V$ somewhere and the third edge elsewhere:
\begin{eqnarray}
\bar V(\includegraphics[scale=0.018]{pl2_1.eps}\,\includegraphics[scale=0.012]{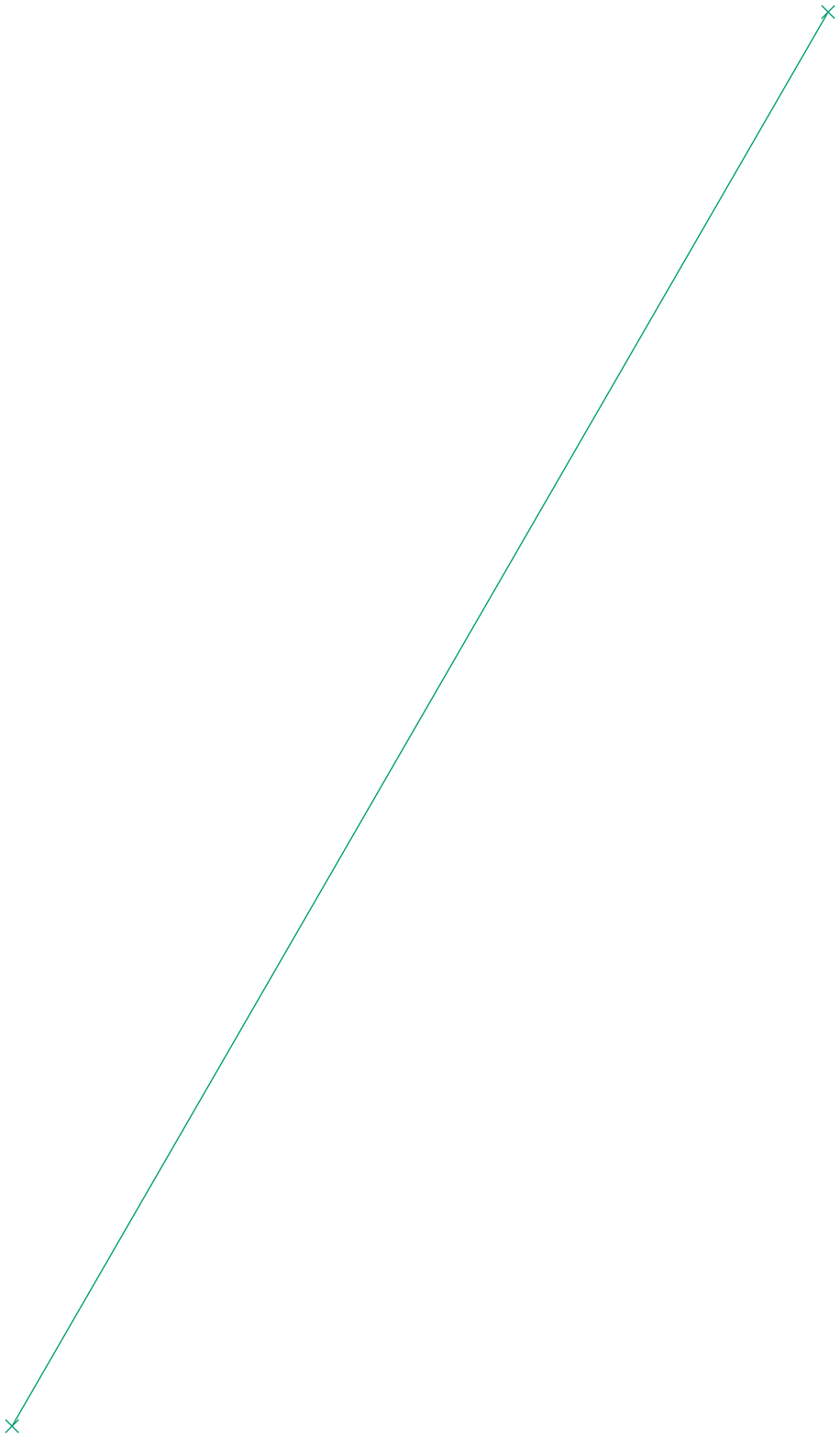}) &=& 
  V(\includegraphics[scale=0.018]{pl2_1.eps}) (M_{n-1}-2) ;\\
\bar V(\includegraphics[scale=0.018]{pl2_2.eps}\,\includegraphics[scale=0.012]{pl1_2C.eps}) &=& 
  V(\includegraphics[scale=0.018]{pl2_2.eps}) (M_{n-1}-2) ;\\
\bar V(\includegraphics[scale=0.018]{pl2_3.eps}\,\includegraphics[scale=0.012]{pl1_2C.eps}) &=& 
  V(\includegraphics[scale=0.018]{pl2_3.eps}) (M_{n-1}-2) ;\\
\bar V(\includegraphics[scale=0.018]{pl2_4.eps}\,\includegraphics[scale=0.012]{pl1_2C.eps}) &=& 
  V(\includegraphics[scale=0.018]{pl2_4.eps}) (M_{n-1}-2) ;\\
\bar V(\includegraphics[scale=0.018]{pl2_5.eps}\,\includegraphics[scale=0.012]{pl1_2C.eps}) &=& 
  V(\includegraphics[scale=0.018]{pl2_5.eps}) (M_{n-1}-2) ;\\
\bar V(\includegraphics[scale=0.018]{pl2_6.eps}\,\includegraphics[scale=0.012]{pl1_2C.eps}) &=& 
  V(\includegraphics[scale=0.018]{pl2_6.eps}) (M_{n-1}-2).
\end{eqnarray}
The green edge in the arguments of these functions indicate an edge with any of the 3 orientations.
The $\bar V$ are the sums of the enumeration of the $V$ at some point and
all $V$ explicitly less than this in the Hasse diagram. Because there is no layer
between the intermediate and lowest level in Figure \ref{fig.pos3}, all Mobius functions
between pairs of these are $-1$. By selecting all 30 brown arrows in the diagram
(9 were only shown)
from the 6 sets of the intermediate to the 14 sets at the lowest level we obtain
\begin{eqnarray}
 V(\includegraphics[scale=0.018]{pl2_1.eps}\;\includegraphics[scale=0.015]{pl1_2C.eps}) &=& 
\bar V(\includegraphics[scale=0.018]{pl2_1.eps}\;\includegraphics[scale=0.015]{pl1_2C.eps}) 
 - V(\includegraphics[scale=0.018]{pl3_3.eps})
 - V(\includegraphics[scale=0.018]{pl3_5.eps})
 - V(\includegraphics[scale=0.018]{pl3_7.eps})
 - V(\includegraphics[scale=0.018]{pl3_13.eps})
 - V(\includegraphics[scale=0.018]{pl3_14.eps})
  ;\\
 V(\includegraphics[scale=0.018]{pl2_2.eps}\;\includegraphics[scale=0.015]{pl1_2C.eps}) &=& 
\bar V(\includegraphics[scale=0.018]{pl2_2.eps}\;\includegraphics[scale=0.015]{pl1_2C.eps}) 
 - V(\includegraphics[scale=0.018]{pl3_1.eps})
 - V(\includegraphics[scale=0.018]{pl3_6.eps})
 - V(\includegraphics[scale=0.018]{pl3_8.eps})
 - V(\includegraphics[scale=0.018]{pl3_9.eps})
 - V(\includegraphics[scale=0.018]{pl3_14.eps})
  ;\\
 V(\includegraphics[scale=0.018]{pl2_3.eps}\;\includegraphics[scale=0.015]{pl1_2C.eps}) &=& 
\bar V(\includegraphics[scale=0.018]{pl2_3.eps}\;\includegraphics[scale=0.015]{pl1_2C.eps}) 
 - V(\includegraphics[scale=0.018]{pl3_2.eps})
 - V(\includegraphics[scale=0.018]{pl3_4.eps})
 - V(\includegraphics[scale=0.018]{pl3_7.eps})
 - V(\includegraphics[scale=0.018]{pl3_9.eps})
 - V(\includegraphics[scale=0.018]{pl3_10.eps})
   ;\\
 V(\includegraphics[scale=0.018]{pl2_4.eps}\;\includegraphics[scale=0.015]{pl1_2C.eps}) &=& 
\bar V(\includegraphics[scale=0.018]{pl2_4.eps}\;\includegraphics[scale=0.015]{pl1_2C.eps}) 
 - V(\includegraphics[scale=0.018]{pl3_3.eps})
 - V(\includegraphics[scale=0.018]{pl3_5.eps})
 - V(\includegraphics[scale=0.018]{pl3_8.eps})
 - V(\includegraphics[scale=0.018]{pl3_10.eps})
 - V(\includegraphics[scale=0.018]{pl3_11.eps})
  ;\\
 V(\includegraphics[scale=0.018]{pl2_5.eps}\;\includegraphics[scale=0.015]{pl1_2C.eps}) &=& 
\bar V(\includegraphics[scale=0.018]{pl2_5.eps}\;\includegraphics[scale=0.015]{pl1_2C.eps}) 
 - V(\includegraphics[scale=0.018]{pl3_1.eps})
 - V(\includegraphics[scale=0.018]{pl3_6.eps})
 - V(\includegraphics[scale=0.018]{pl3_7.eps})
 - V(\includegraphics[scale=0.018]{pl3_11.eps})
 - V(\includegraphics[scale=0.018]{pl3_12.eps})
  ;\\
 V(\includegraphics[scale=0.018]{pl2_6.eps}\includegraphics[scale=0.015]{pl1_2C.eps}) &=& 
\bar V(\includegraphics[scale=0.018]{pl2_6.eps}\includegraphics[scale=0.015]{pl1_2C.eps}) 
 - V(\includegraphics[scale=0.018]{pl3_2.eps})
 - V(\includegraphics[scale=0.018]{pl3_4.eps})
 - V(\includegraphics[scale=0.018]{pl3_8.eps})
 - V(\includegraphics[scale=0.018]{pl3_13.eps})
 - V(\includegraphics[scale=0.018]{pl3_12.eps})
  .
\end{eqnarray}
Gathering the 6 terms of rank 1 and the 14 terms of
of rank 2 as a correction to \eqref{eq.bino} yields
\cite{DresdenLozenge}\cite[A326368]{sloane}
\begin{multline}
L_{n,3}=\binom{M_{n-1}}{3}
-\big[
V(\includegraphics[scale=0.02]{pl3_1.eps})
+
V(\includegraphics[scale=0.018]{pl3_2.eps})
+
V(\includegraphics[scale=0.02]{pl3_3.eps})
+
V(\includegraphics[scale=0.02]{pl3_4.eps})
+
V(\includegraphics[scale=0.018]{pl3_5.eps})
+
V(\includegraphics[scale=0.02]{pl3_6.eps})
+
V(\includegraphics[scale=0.02]{pl3_7.eps})
+
V(\includegraphics[scale=0.02]{pl3_8.eps})
\\
+
V(\includegraphics[scale=0.02]{pl3_9.eps})
+
V(\includegraphics[scale=0.018]{pl3_10.eps})
+
V(\includegraphics[scale=0.02]{pl3_11.eps})
+
V(\includegraphics[scale=0.02]{pl3_12.eps})
+
V(\includegraphics[scale=0.018]{pl3_13.eps})
+
V(\includegraphics[scale=0.02]{pl3_14.eps})
\\
+
 V(\includegraphics[scale=0.018]{pl2_1.eps}\;\includegraphics[scale=0.015]{pl1_2C.eps}) 
+
 V(\includegraphics[scale=0.018]{pl2_2.eps}\;\includegraphics[scale=0.015]{pl1_2C.eps}) 
+
 V(\includegraphics[scale=0.018]{pl2_3.eps}\;\includegraphics[scale=0.015]{pl1_2C.eps}) 
+
 V(\includegraphics[scale=0.018]{pl2_4.eps}\;\includegraphics[scale=0.015]{pl1_2C.eps}) 
+
 V(\includegraphics[scale=0.018]{pl2_5.eps}\;\includegraphics[scale=0.015]{pl1_2C.eps}) 
+
 V(\includegraphics[scale=0.018]{pl2_6.eps}\;\includegraphics[scale=0.015]{pl1_2C.eps}) 
\big]
=
\\
= \frac{1}{16}(n - 2) (9 n^5  -  9 n^4  -  81 n^3  +  81 n^2  +  160 n  -  192),\quad n\ge 2.
\end{multline}
The (inverse) binomial transform is
\begin{equation}
L_{n,3} = 24-22\binom{n}{1}+20\binom{n}{2}+378\binom{n}{4}+810\binom{n}{5}+405\binom{n}{6},\quad n\ge 2.
\end{equation}

\subsection{Four Lozenges}
\begin{figure}
\includegraphics[scale=0.7]{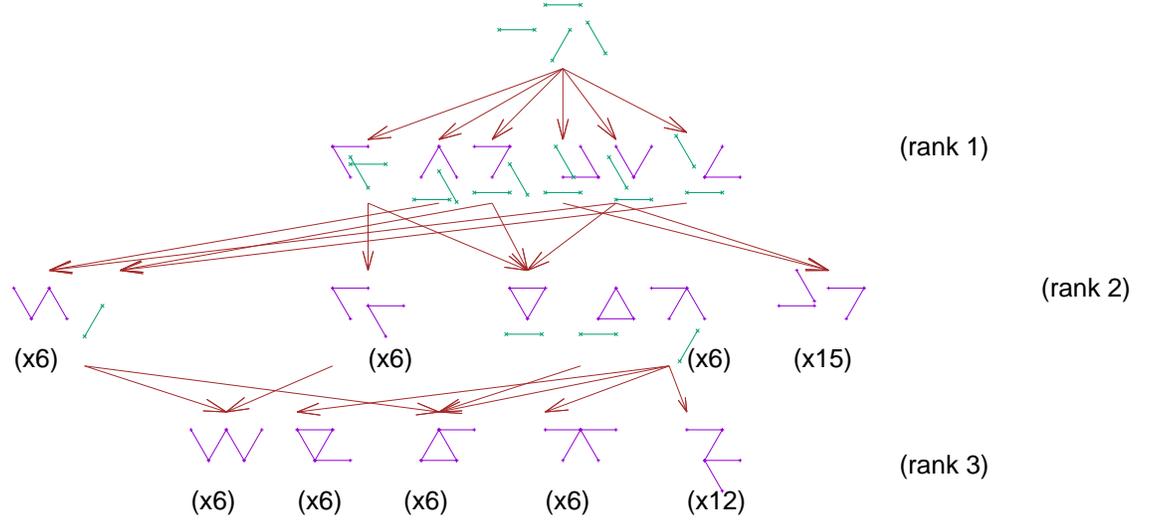}
\caption{Sketch of the Hasse diagram of the deleted edges for $l=4$ deleted unit edges in general, 
and specializations stepping down in the layers of the diagram.}
\label{fig.pos4}
\end{figure}

Figure \ref{fig.poly4} supports 36 fixed polyedges with 4 edges
which fit as follows in the big triangle:
\begin{equation}
V(\includegraphics[scale=0.02]{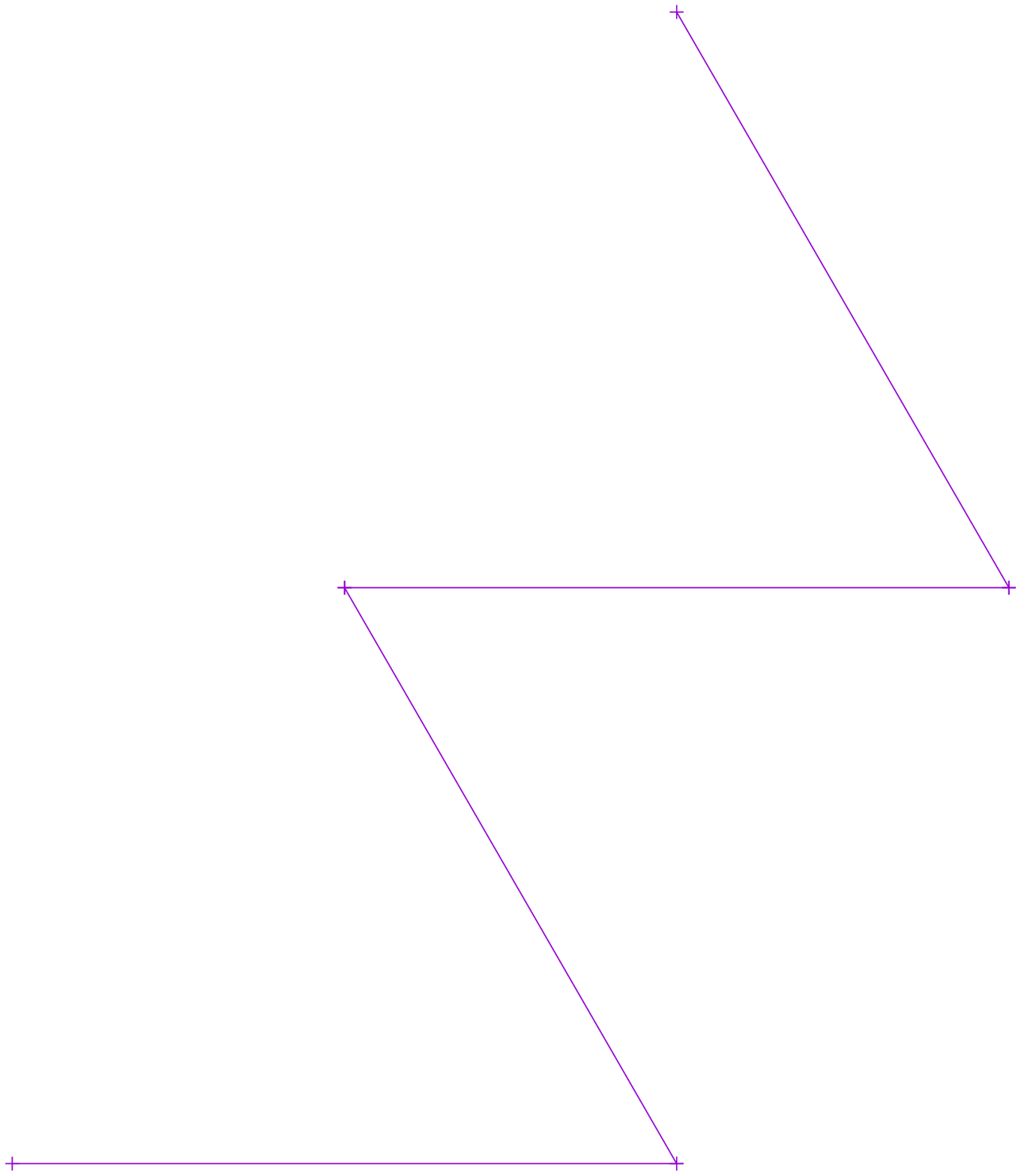})
=
V(\includegraphics[scale=0.018]{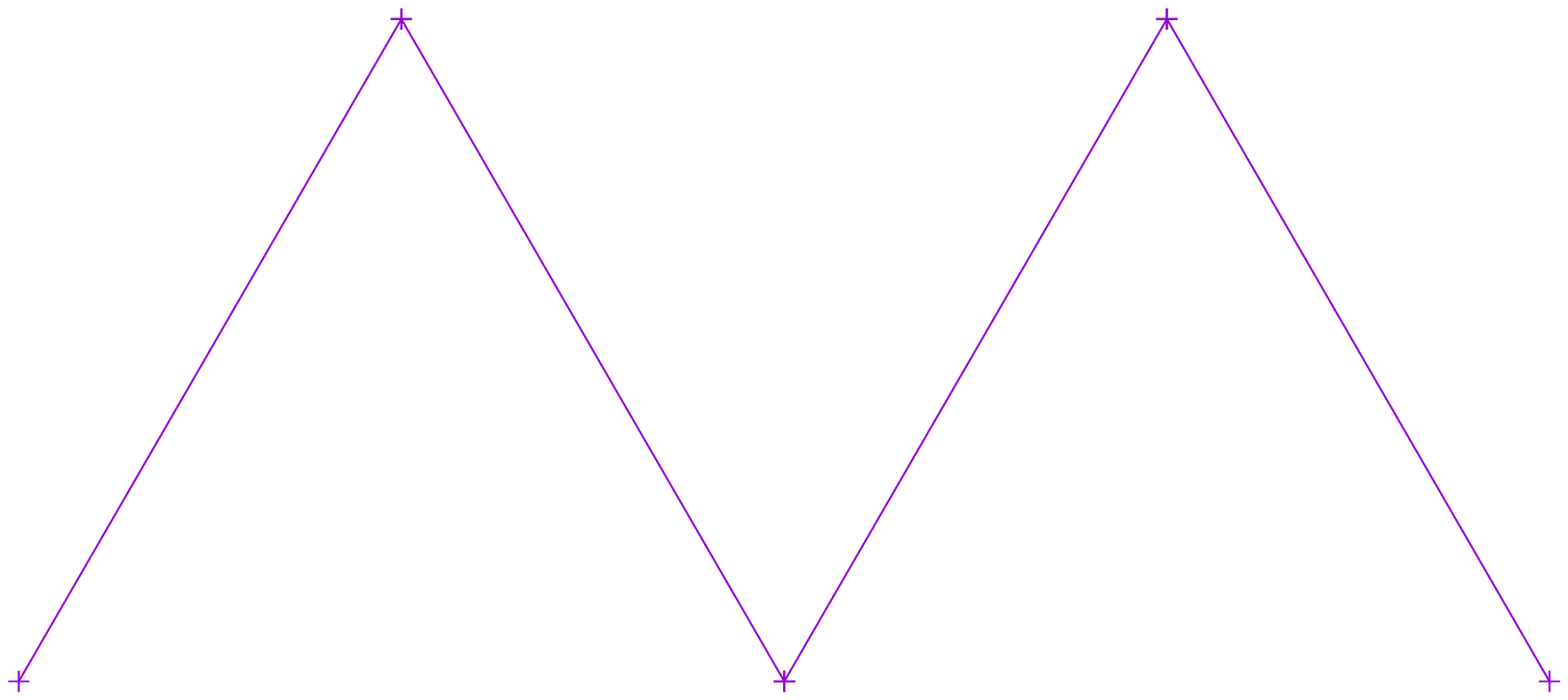})
=
V(\includegraphics[scale=0.02]{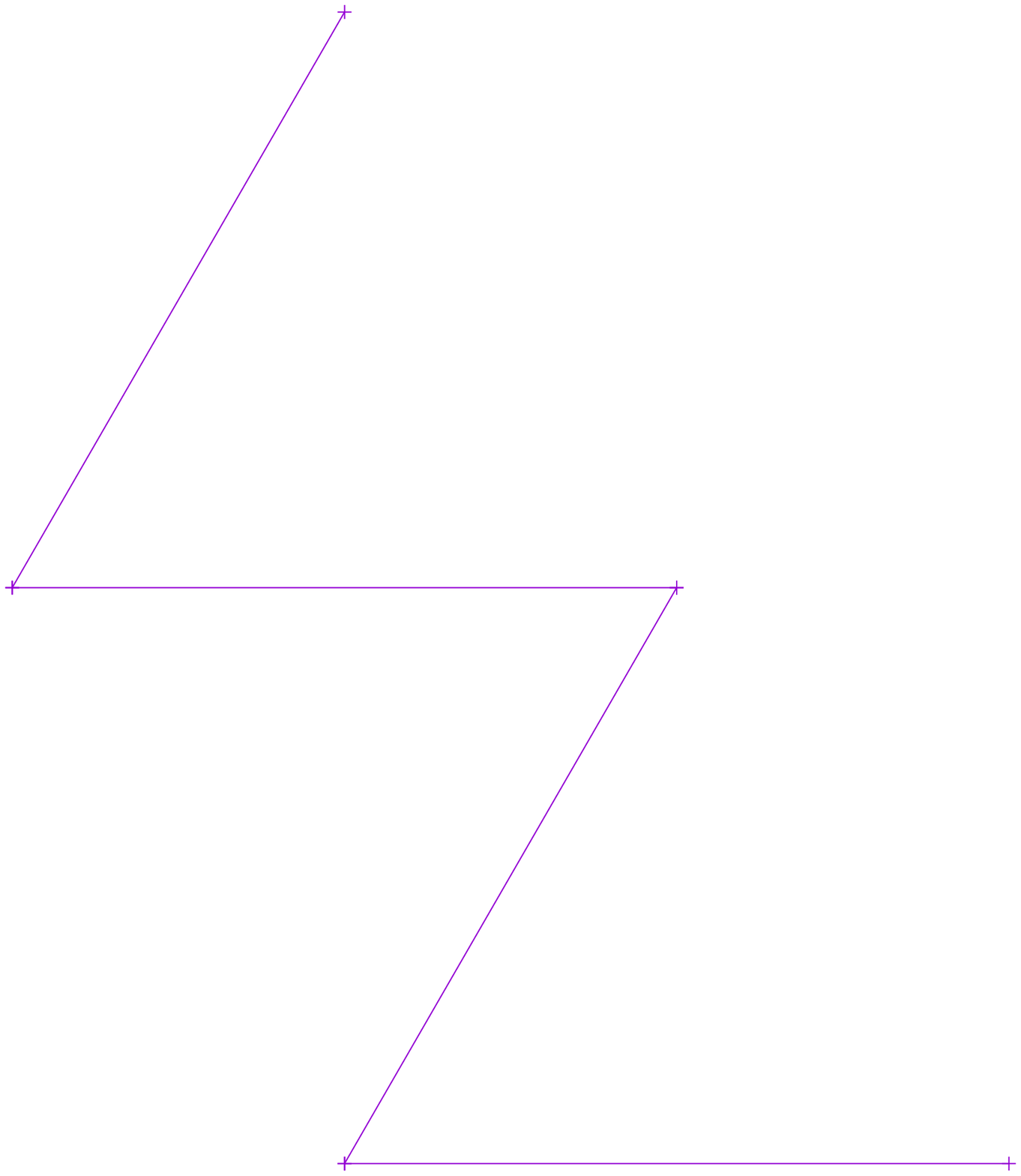})
=T_{n-3};\quad
V(\includegraphics[scale=0.02]{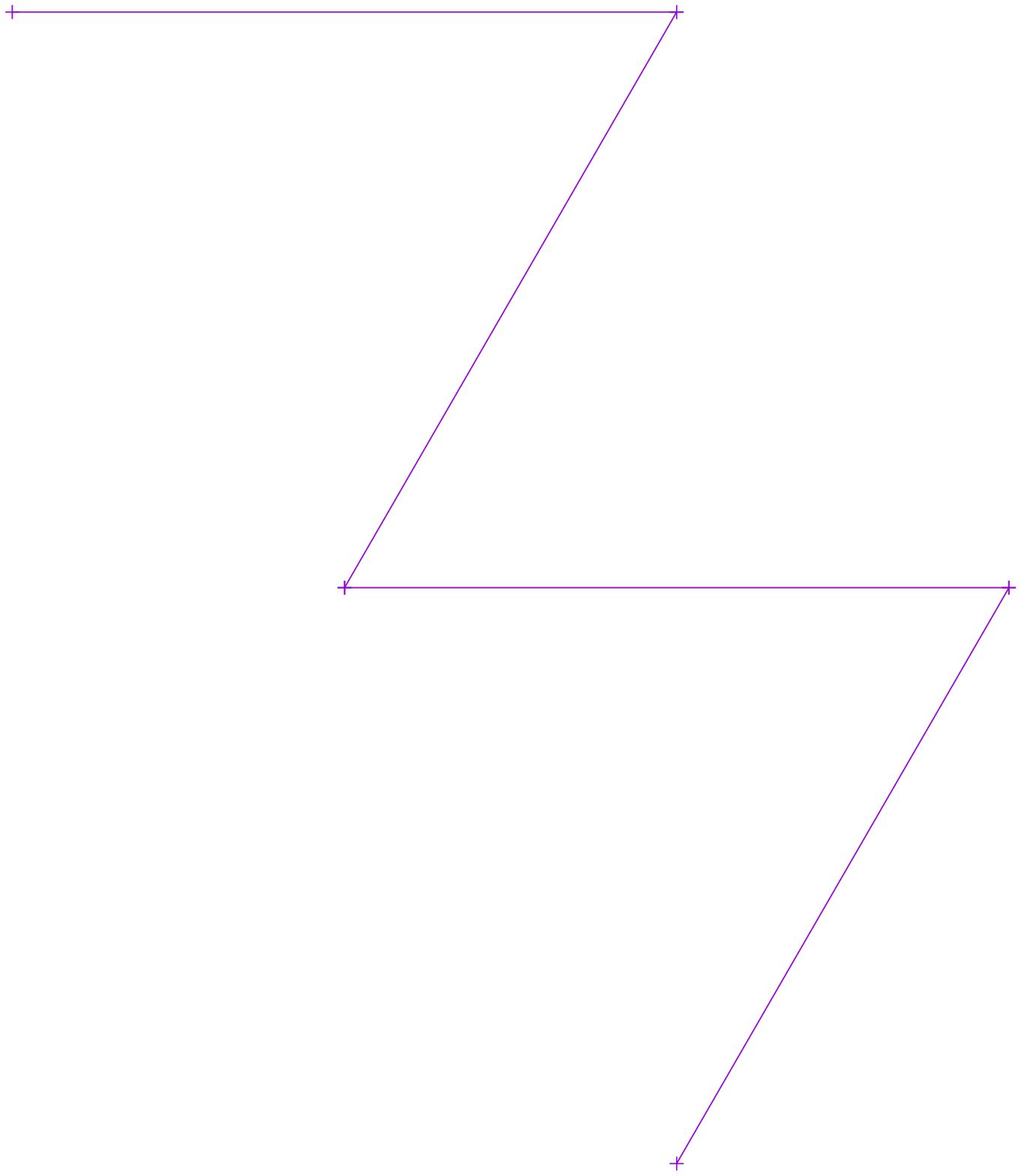})
=
V(\includegraphics[scale=0.02]{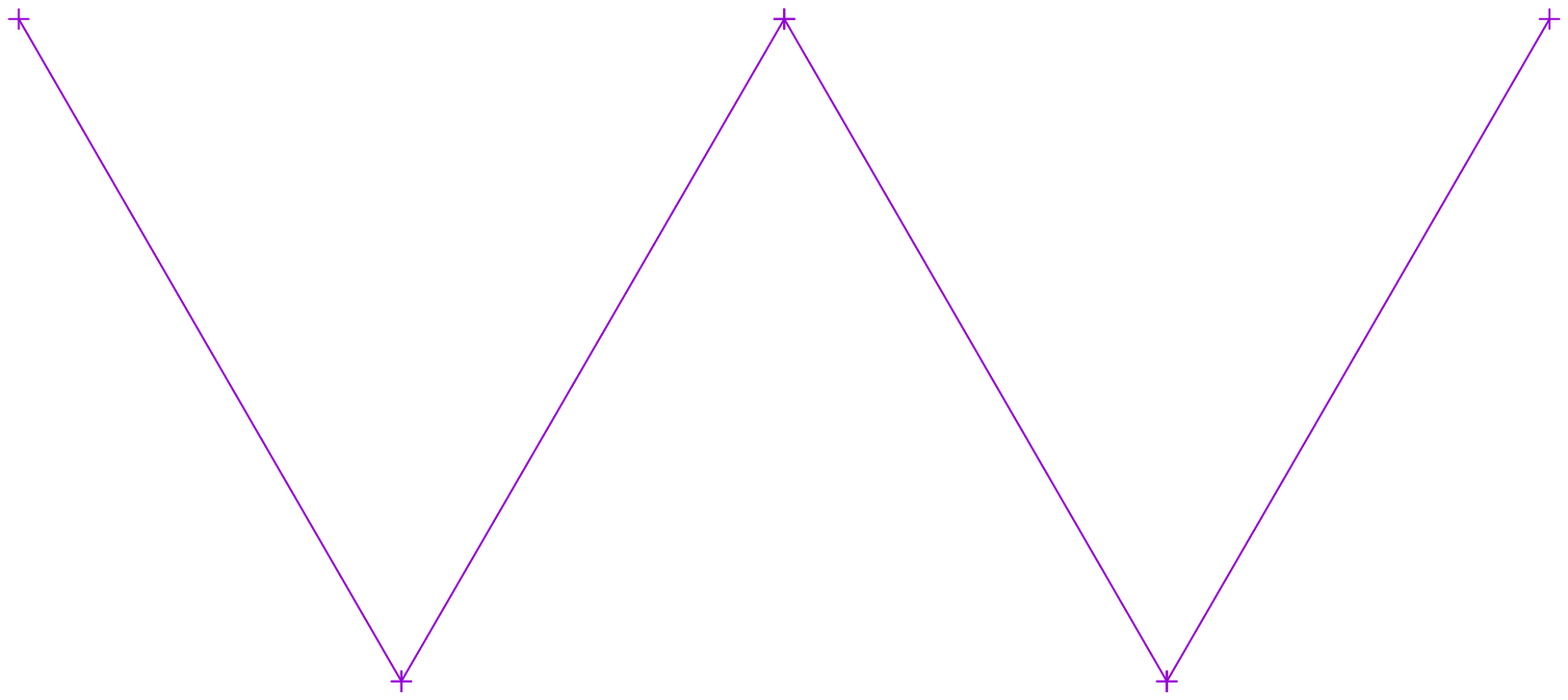})
=
V(\includegraphics[scale=0.018]{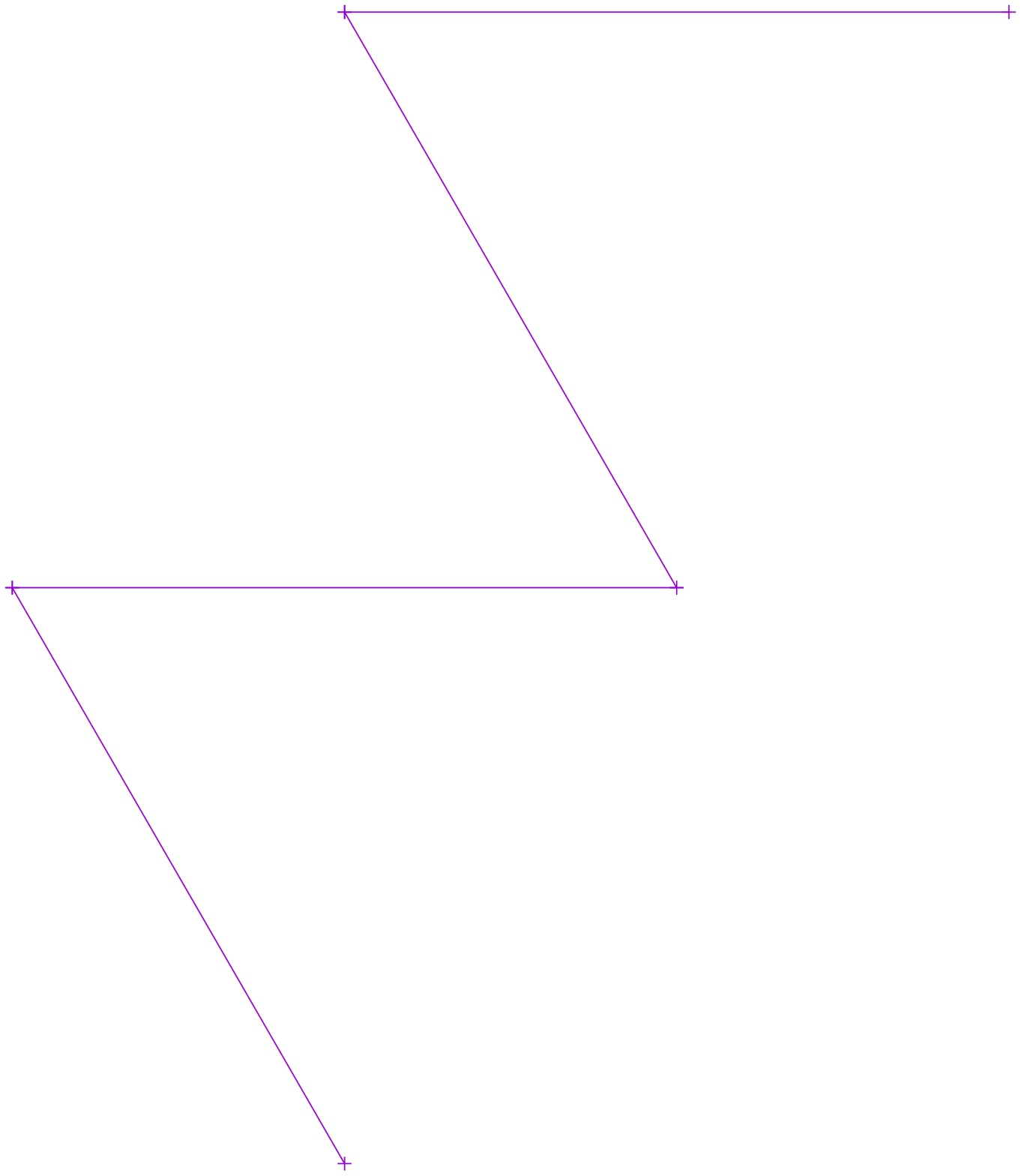})
=T_{n-2};
\label{eq.pl4first}
\end{equation}

\begin{equation}
V(\includegraphics[scale=0.02]{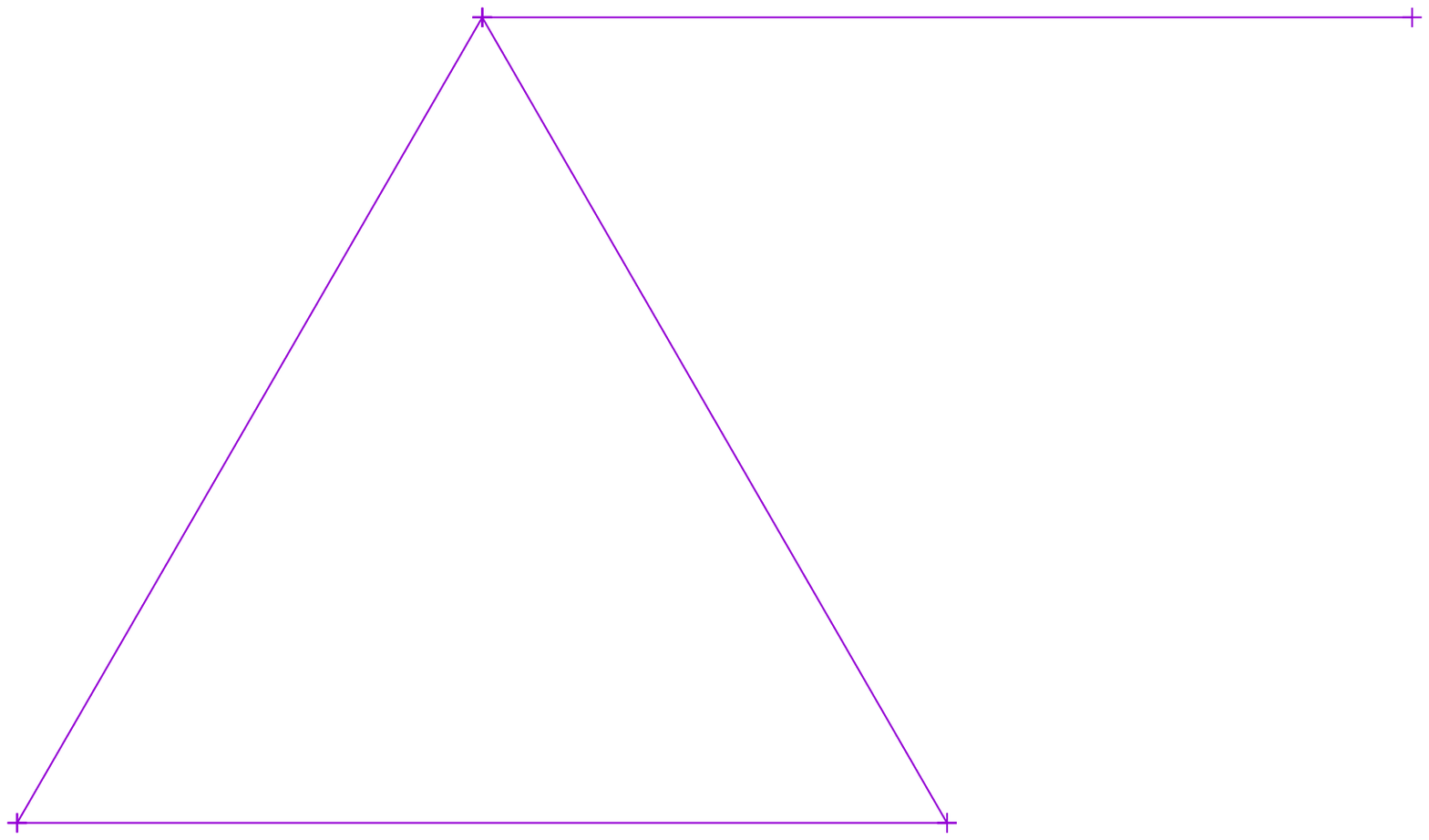})
=
V(\includegraphics[scale=0.02]{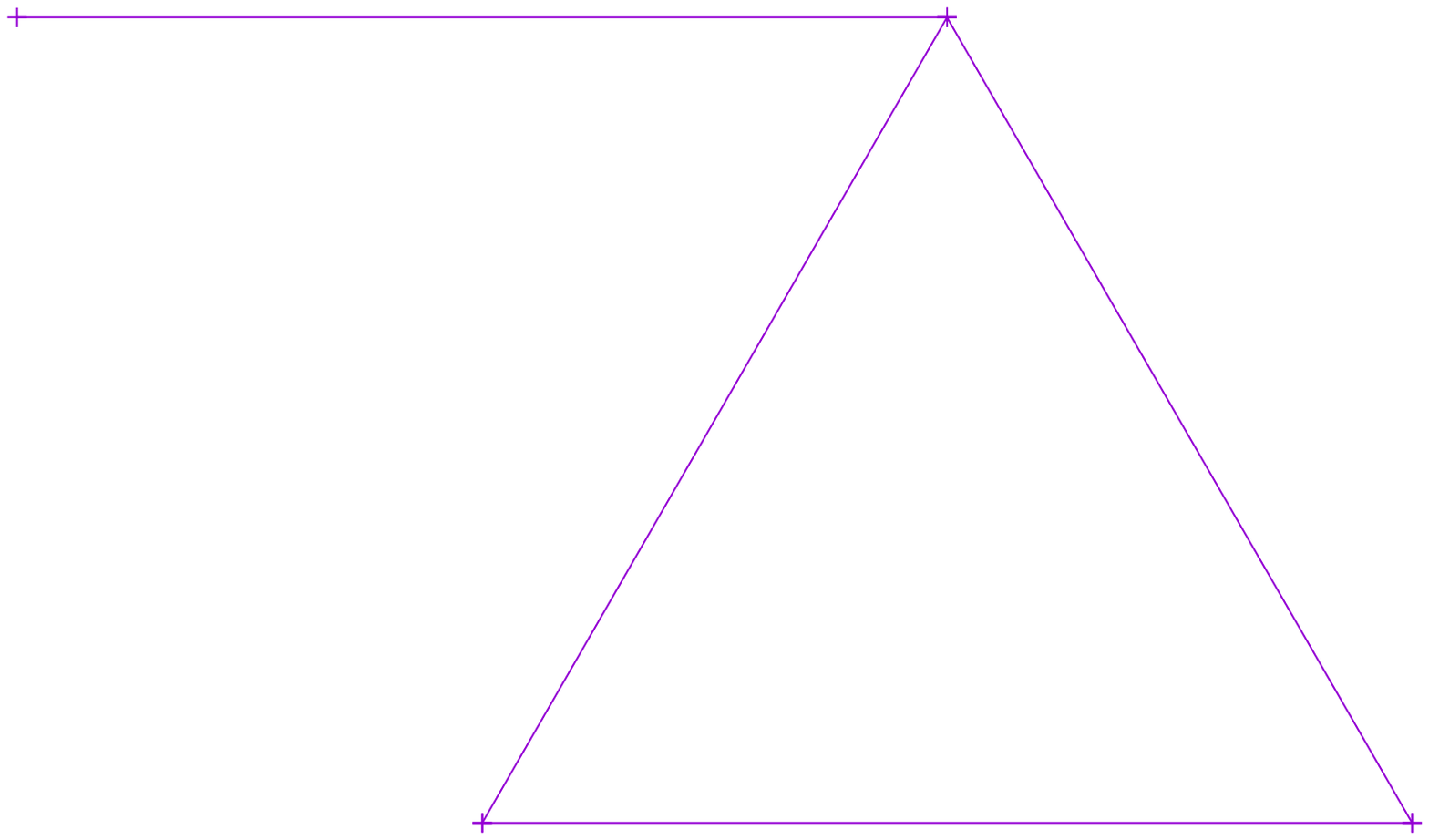})
=
V(\includegraphics[scale=0.02]{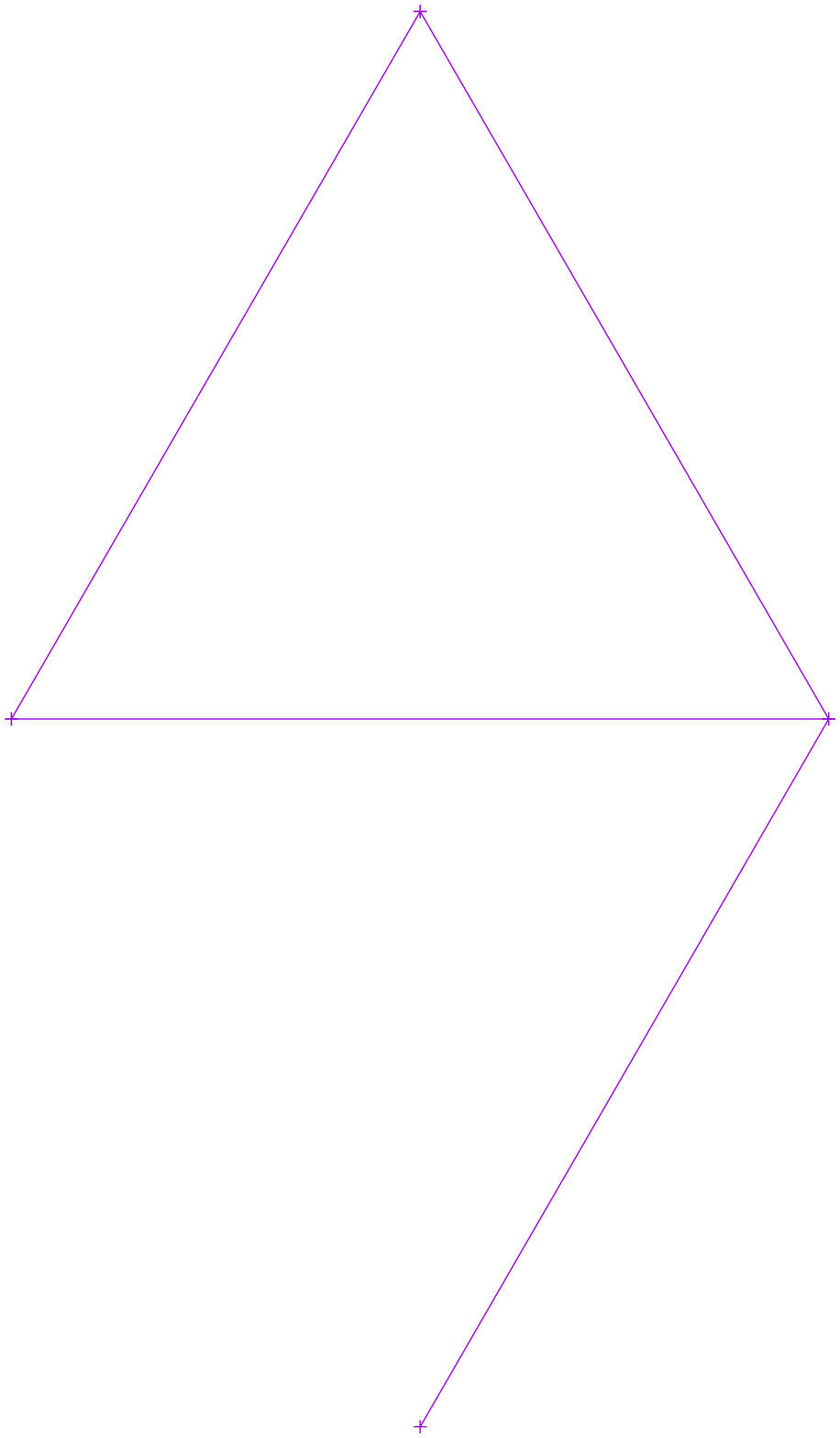})
=
V(\includegraphics[scale=0.02]{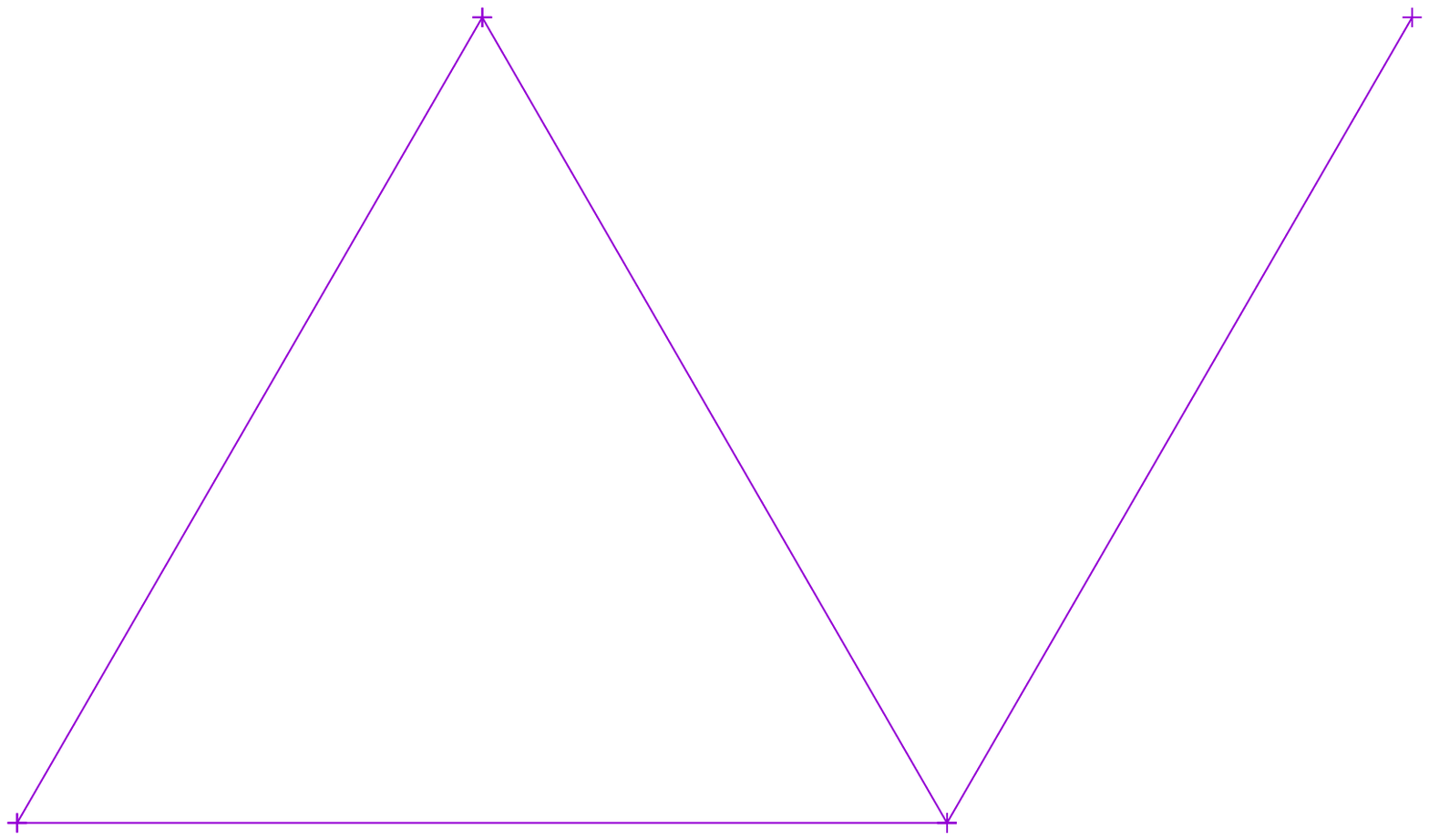})
=
V(\includegraphics[scale=0.02]{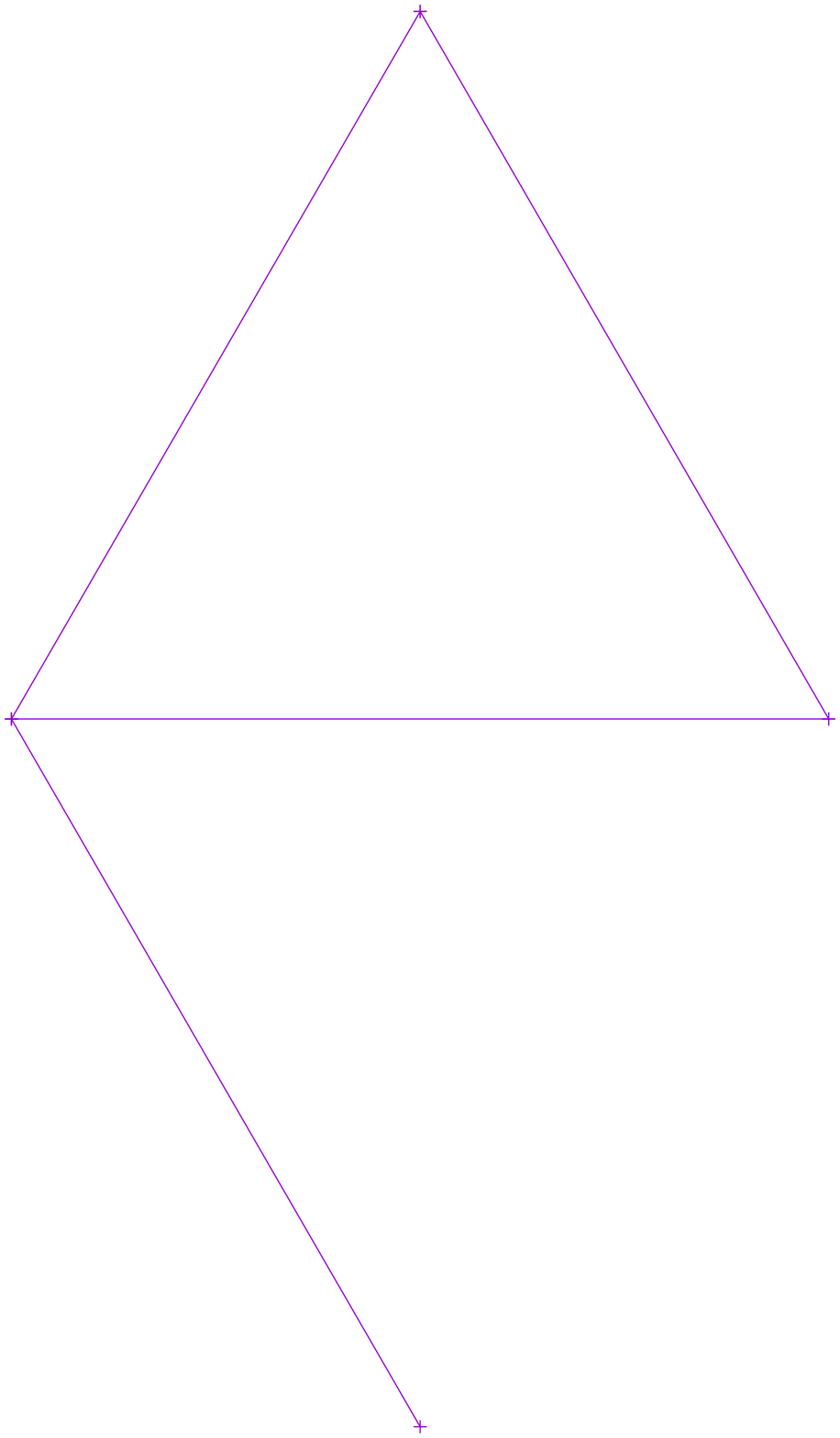})
=
V(\includegraphics[scale=0.02]{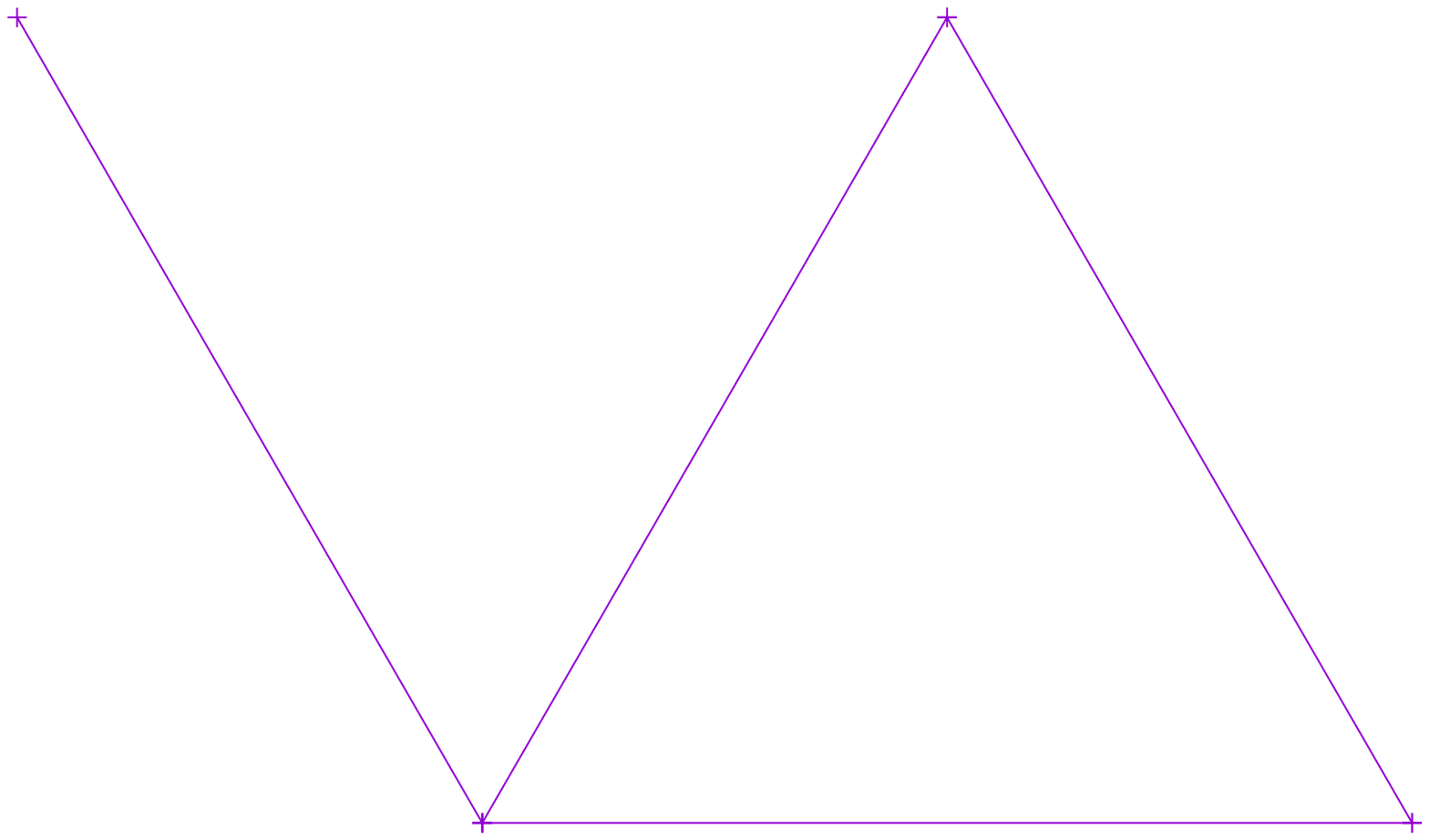})
=T_{n-3};
\end{equation}

\begin{equation}
V(\includegraphics[scale=0.02]{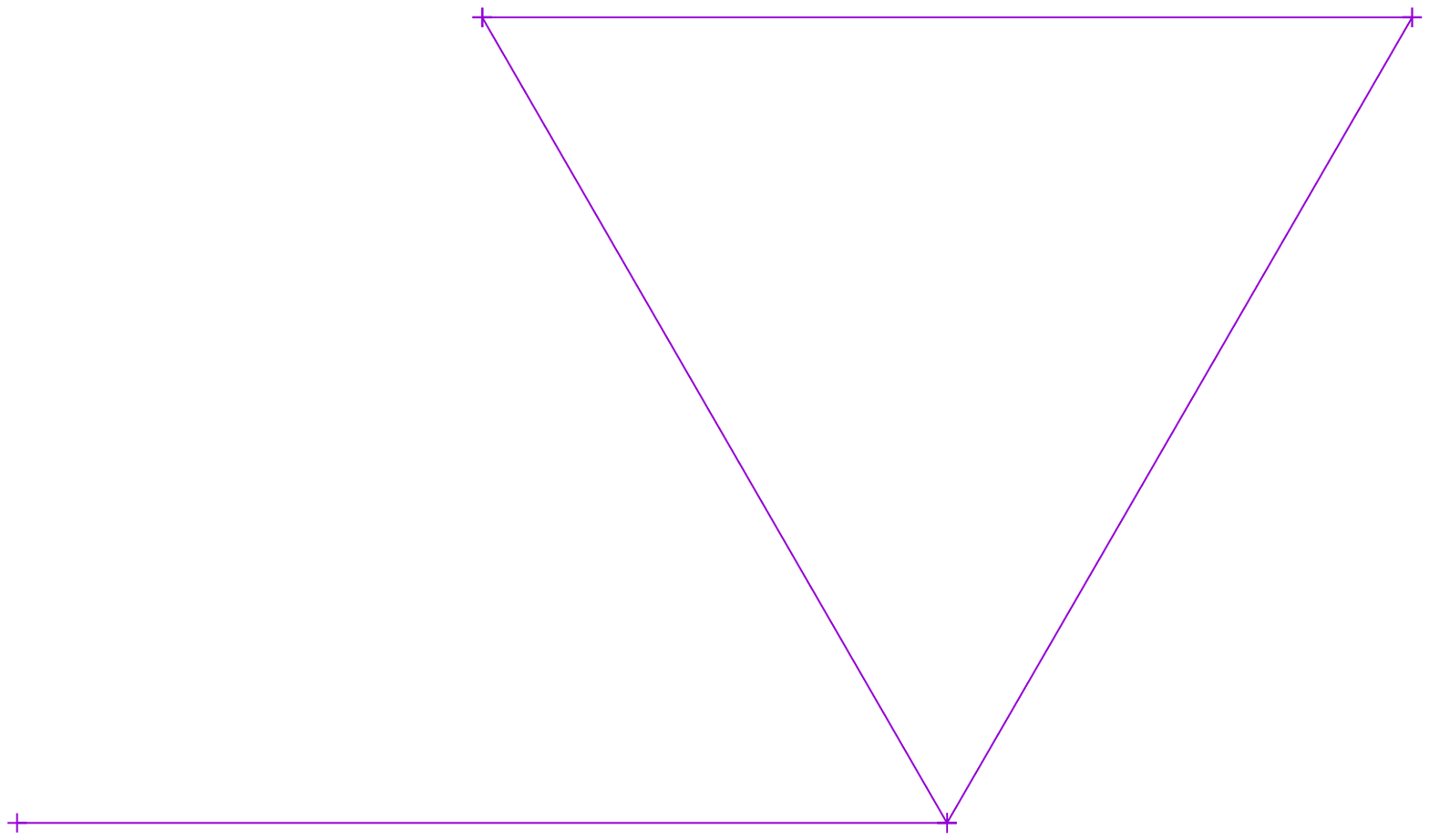})
=
V(\includegraphics[scale=0.02]{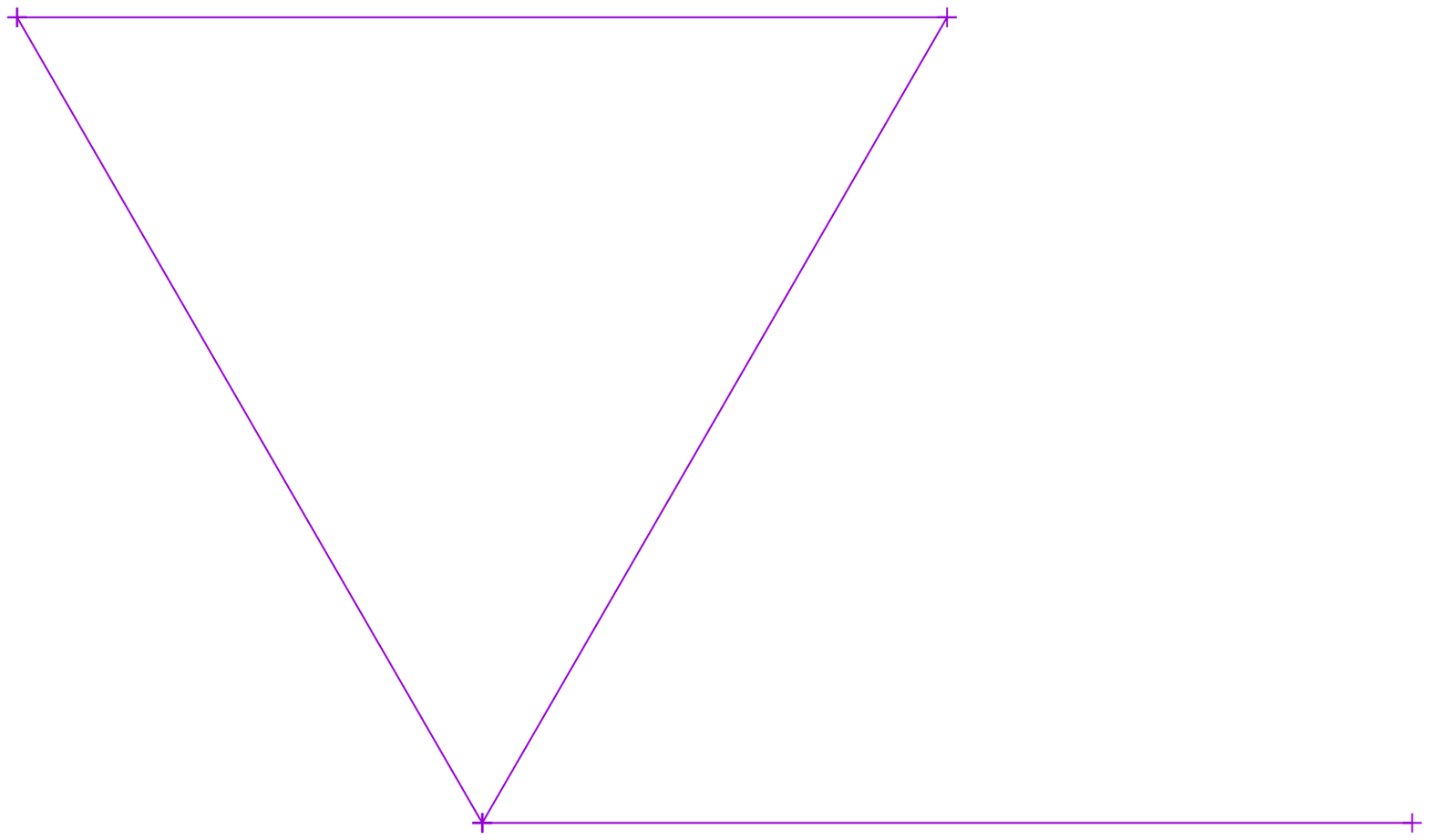})
=
V(\includegraphics[scale=0.02]{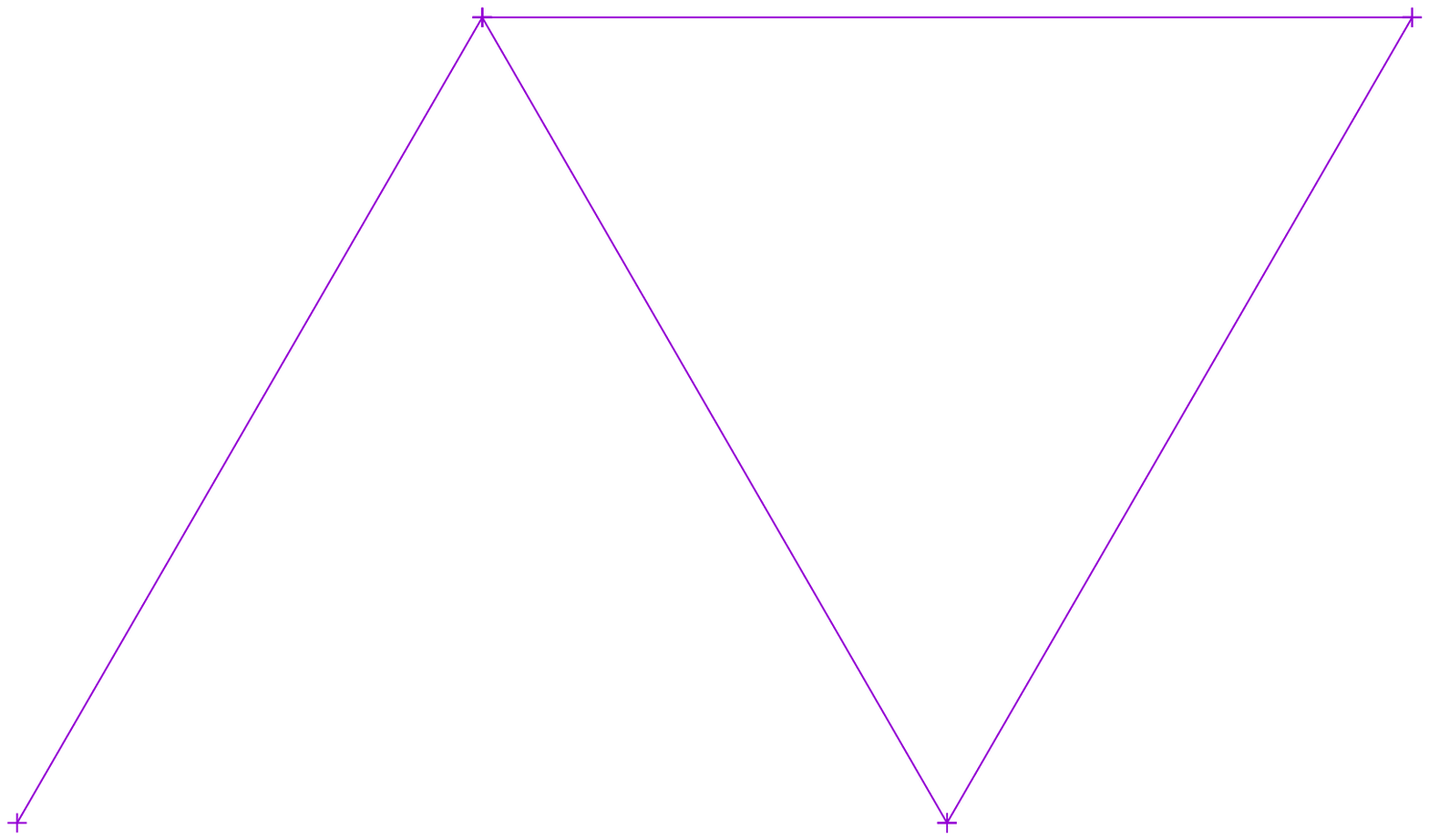})
=
V(\includegraphics[scale=0.02]{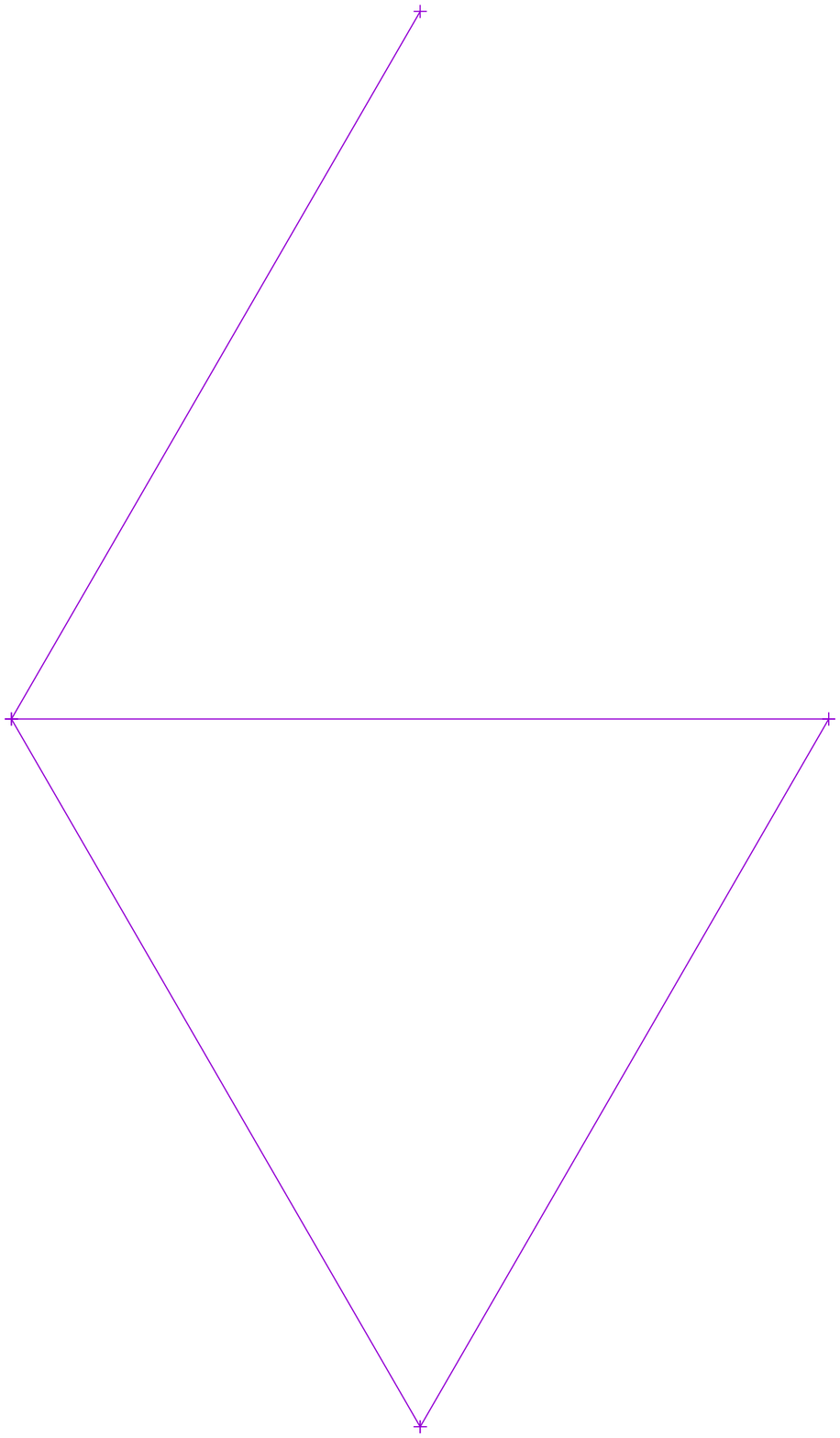})
=
V(\includegraphics[scale=0.02]{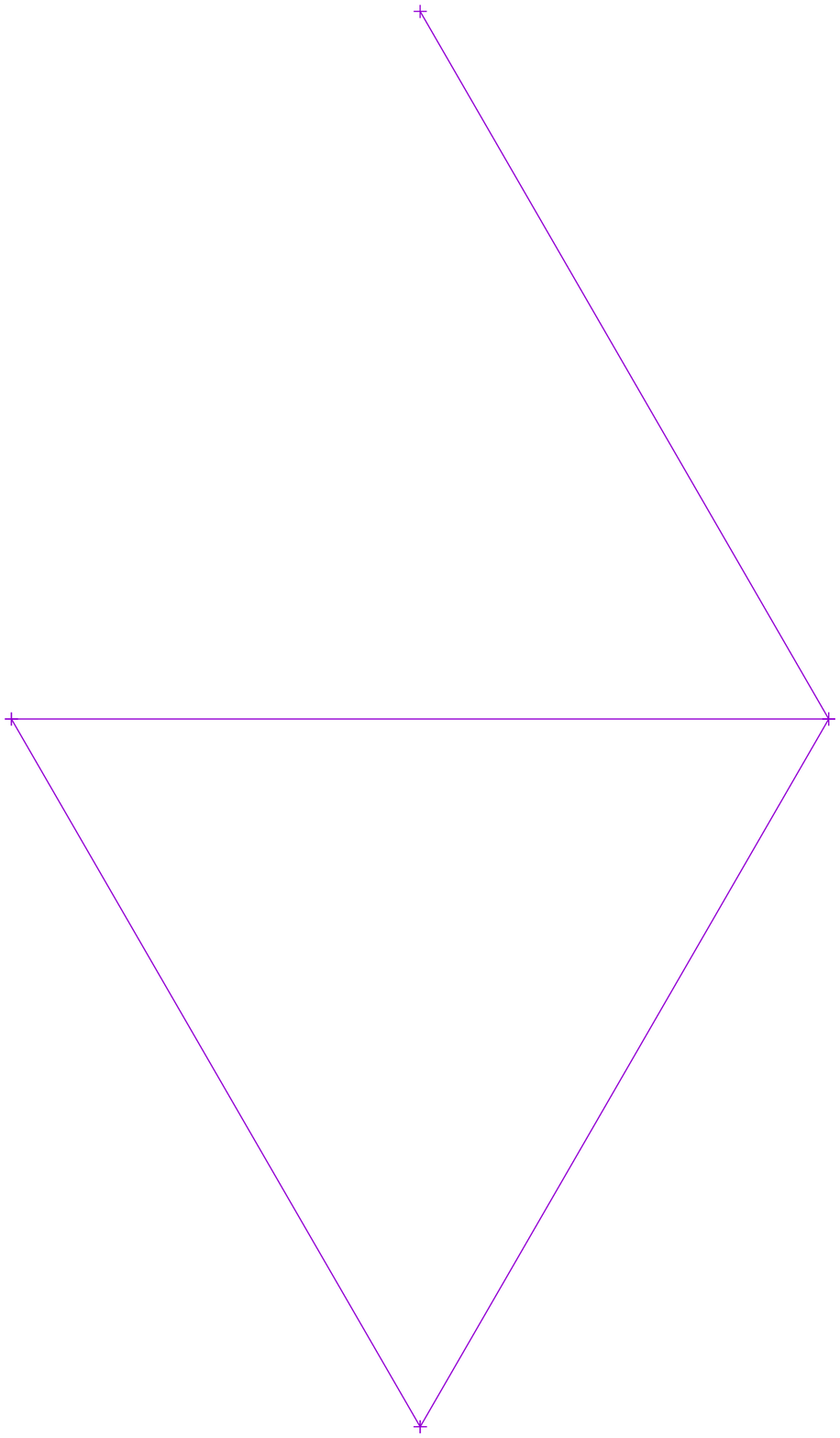})
=
V(\includegraphics[scale=0.02]{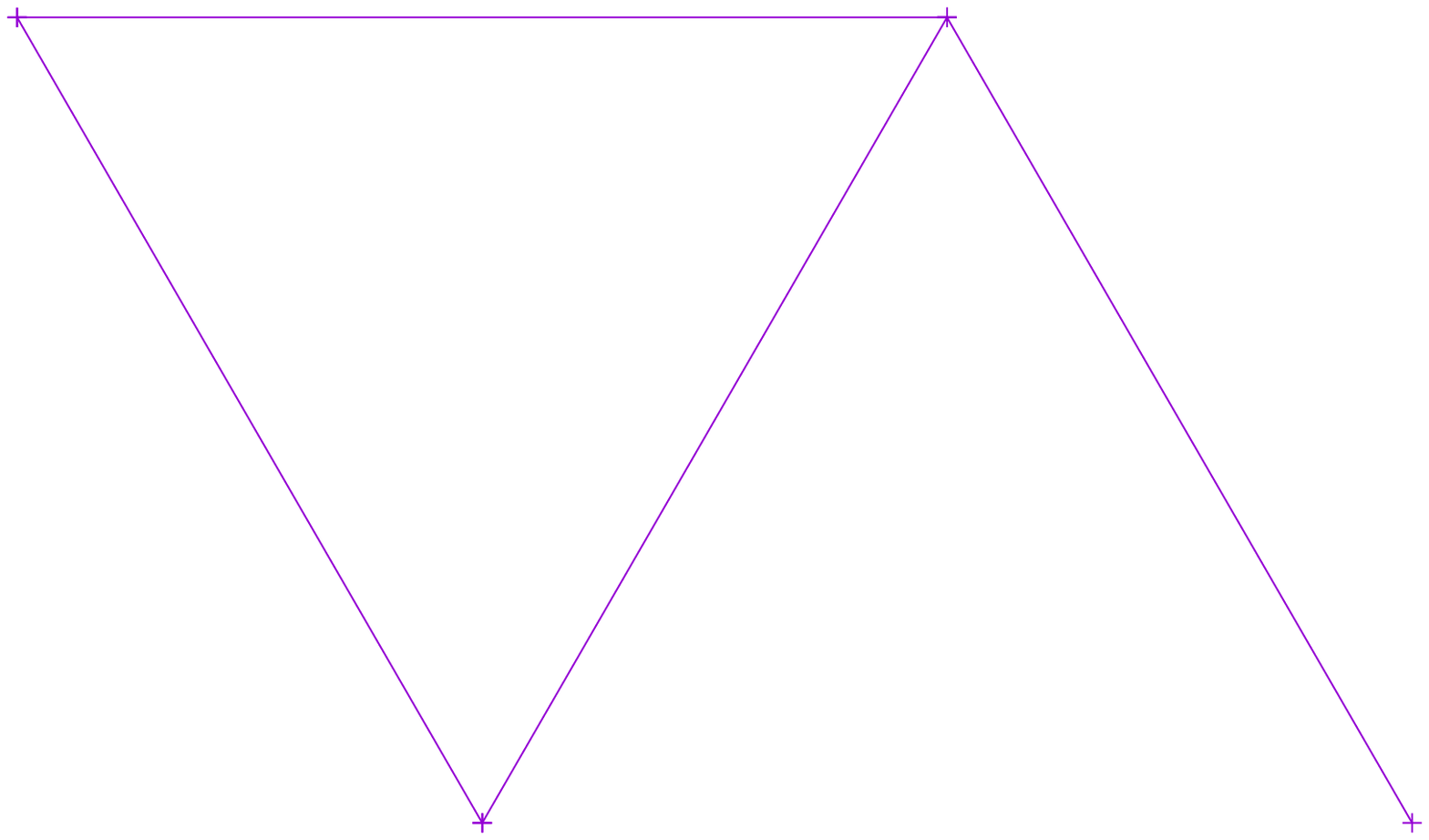})
=T_{n-2};
\end{equation}

\begin{equation}
V(\includegraphics[scale=0.02]{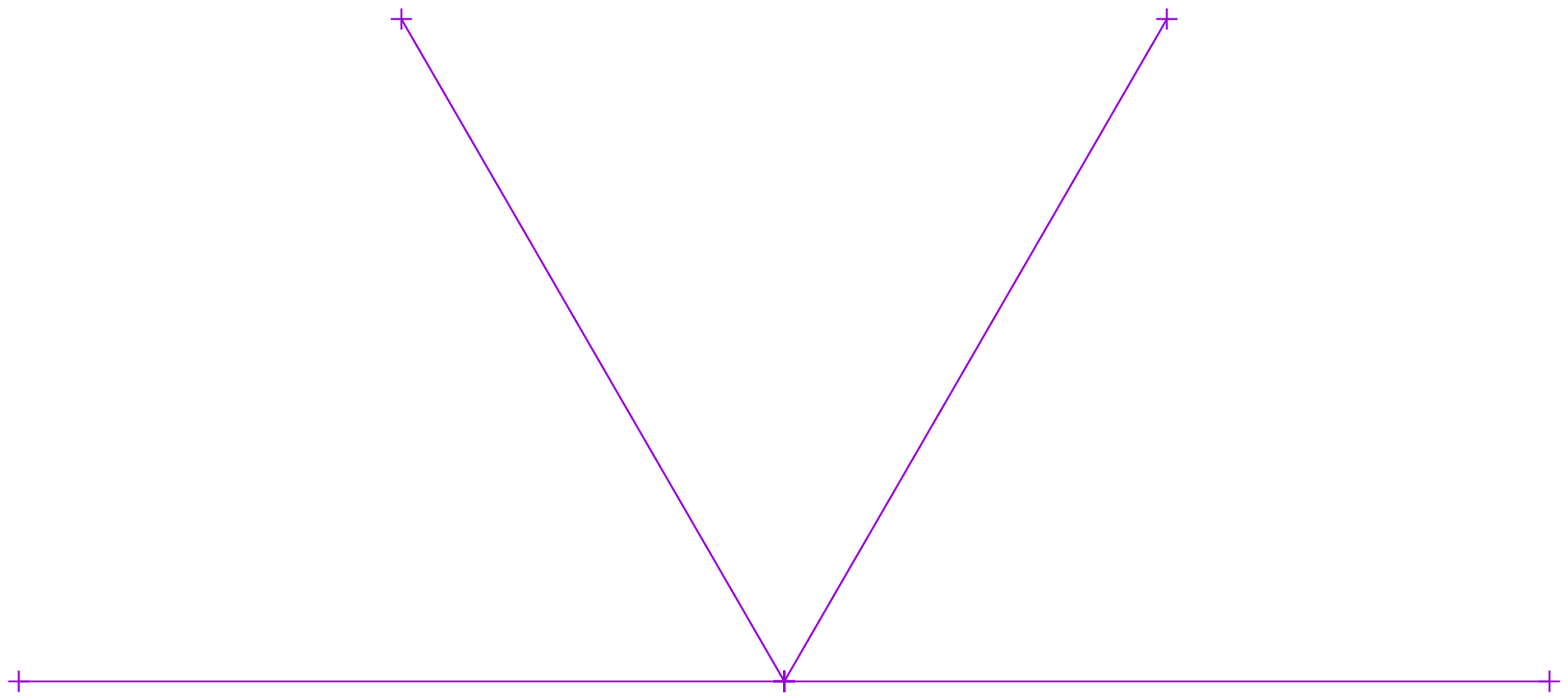})
=
V(\includegraphics[scale=0.02]{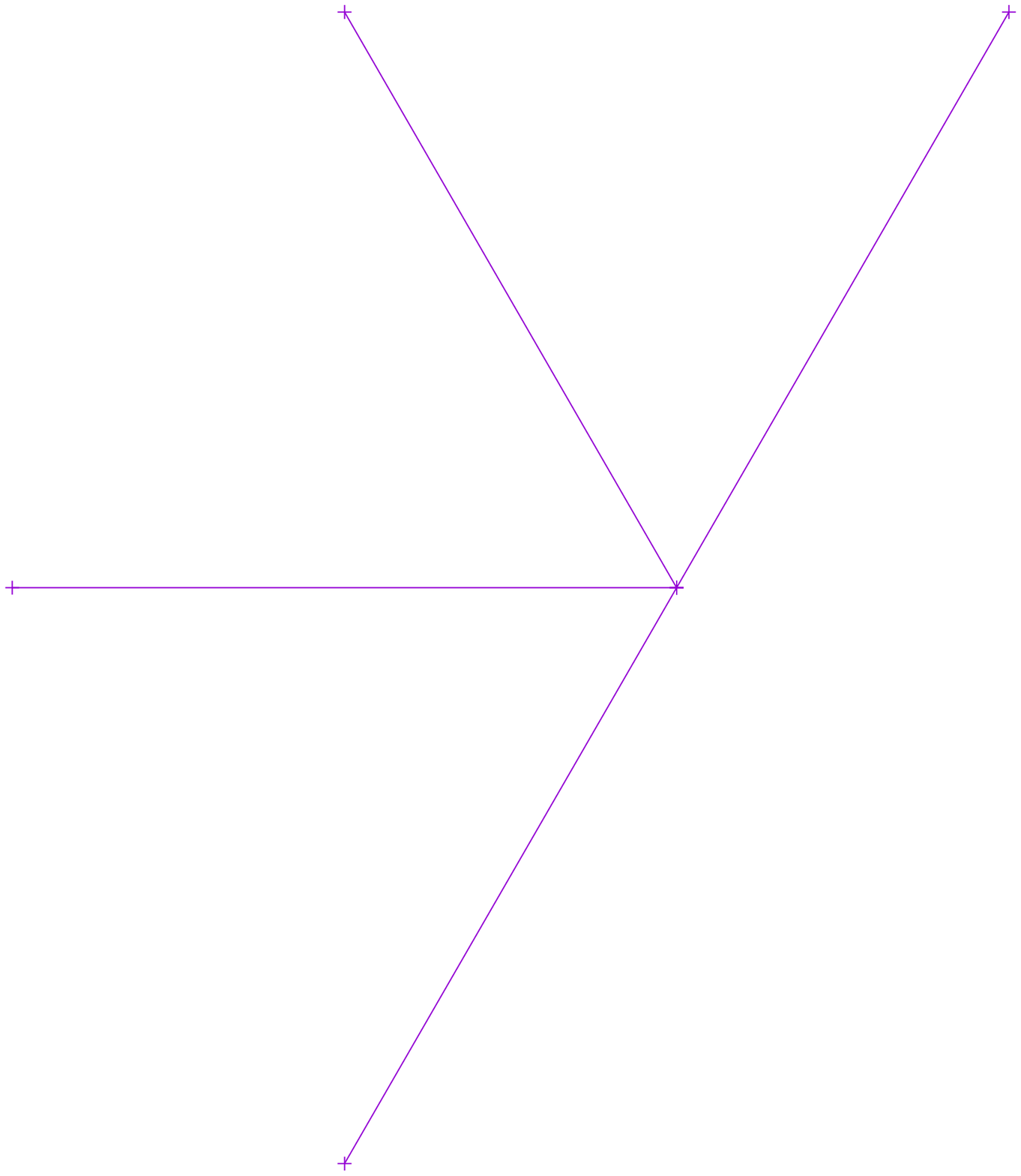})
=
V(\includegraphics[scale=0.02]{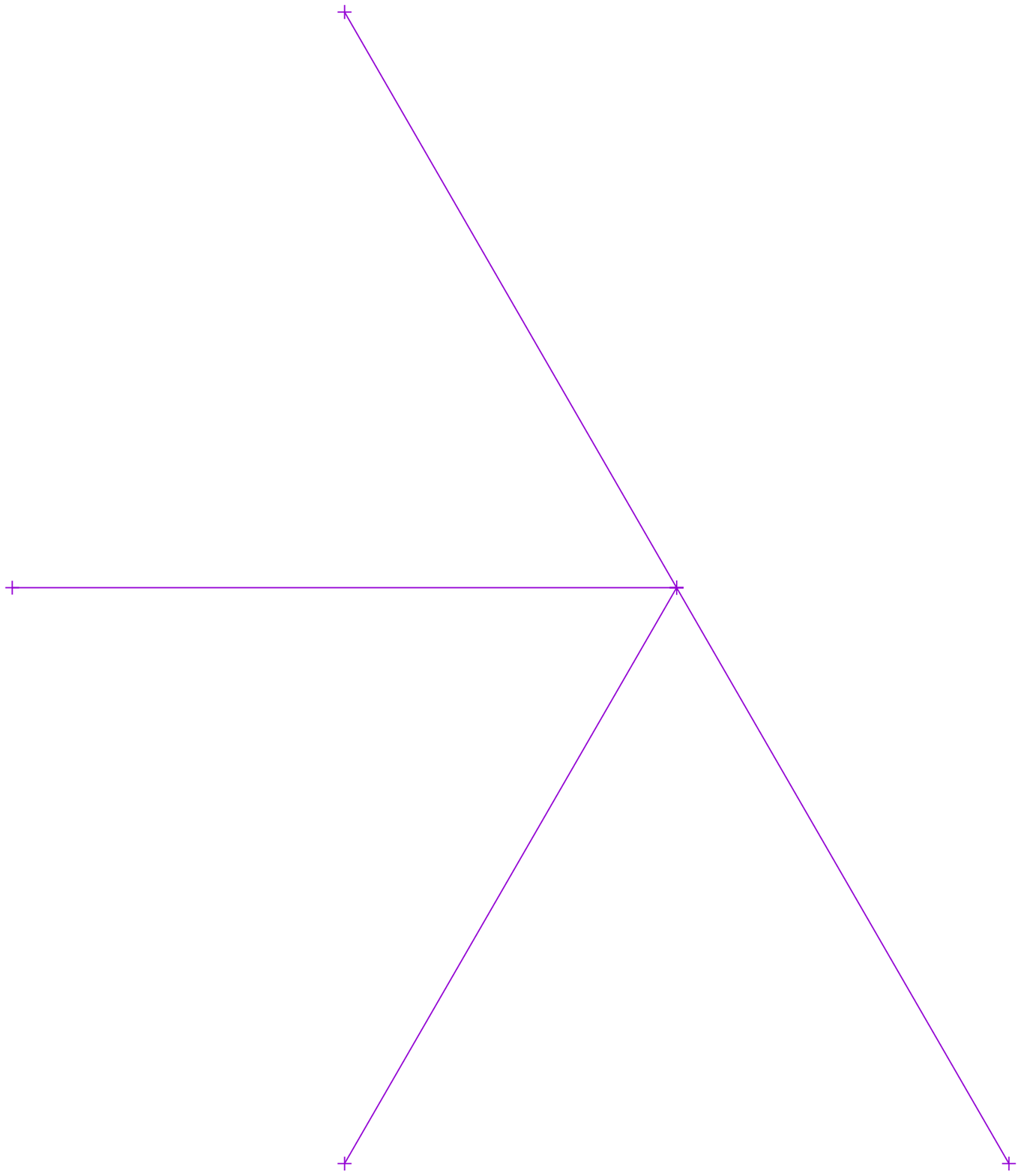})
=
V(\includegraphics[scale=0.02]{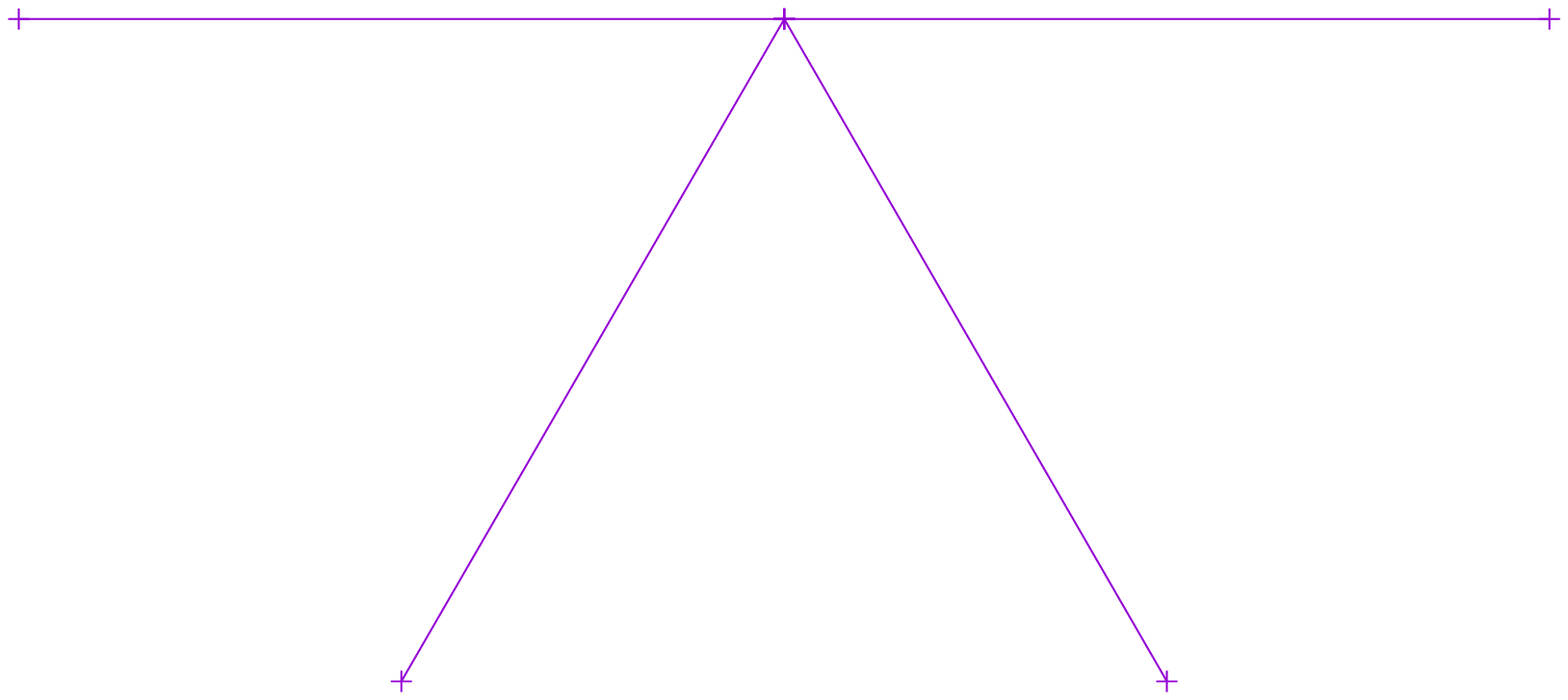})
=
V(\includegraphics[scale=0.02]{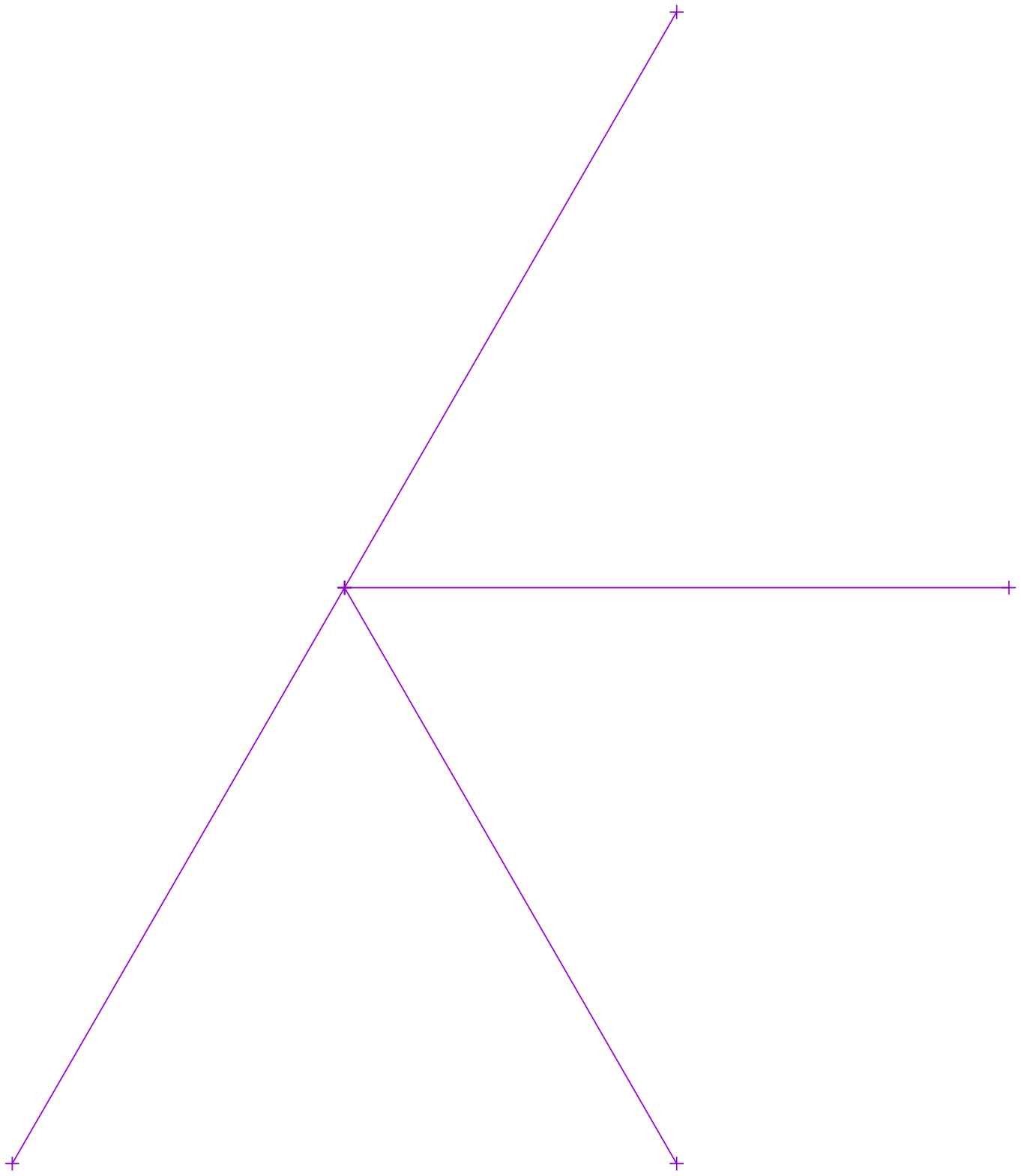})
=
V(\includegraphics[scale=0.02]{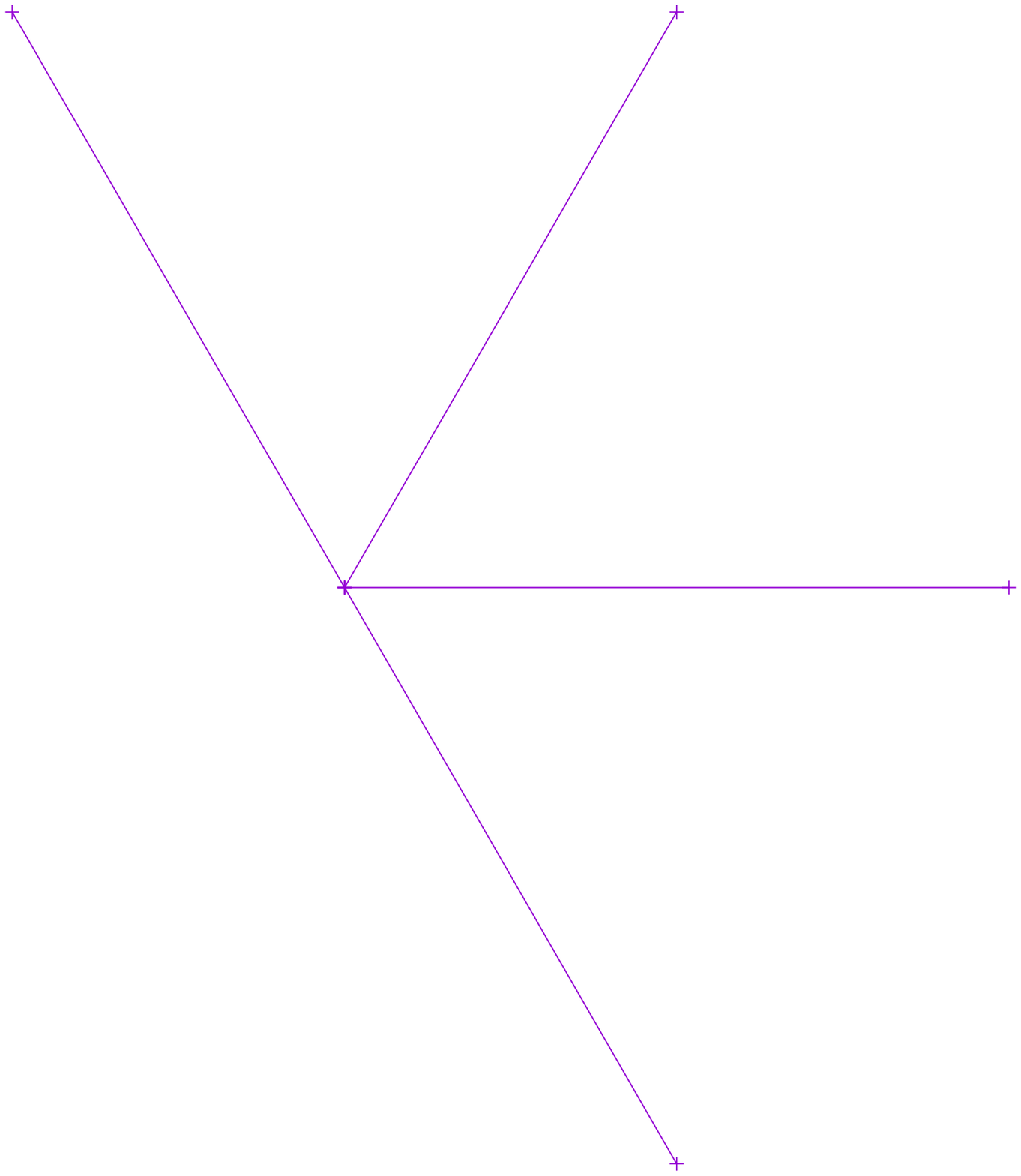})
=T_{n-2};
\end{equation}

\begin{equation}
V(\includegraphics[scale=0.02]{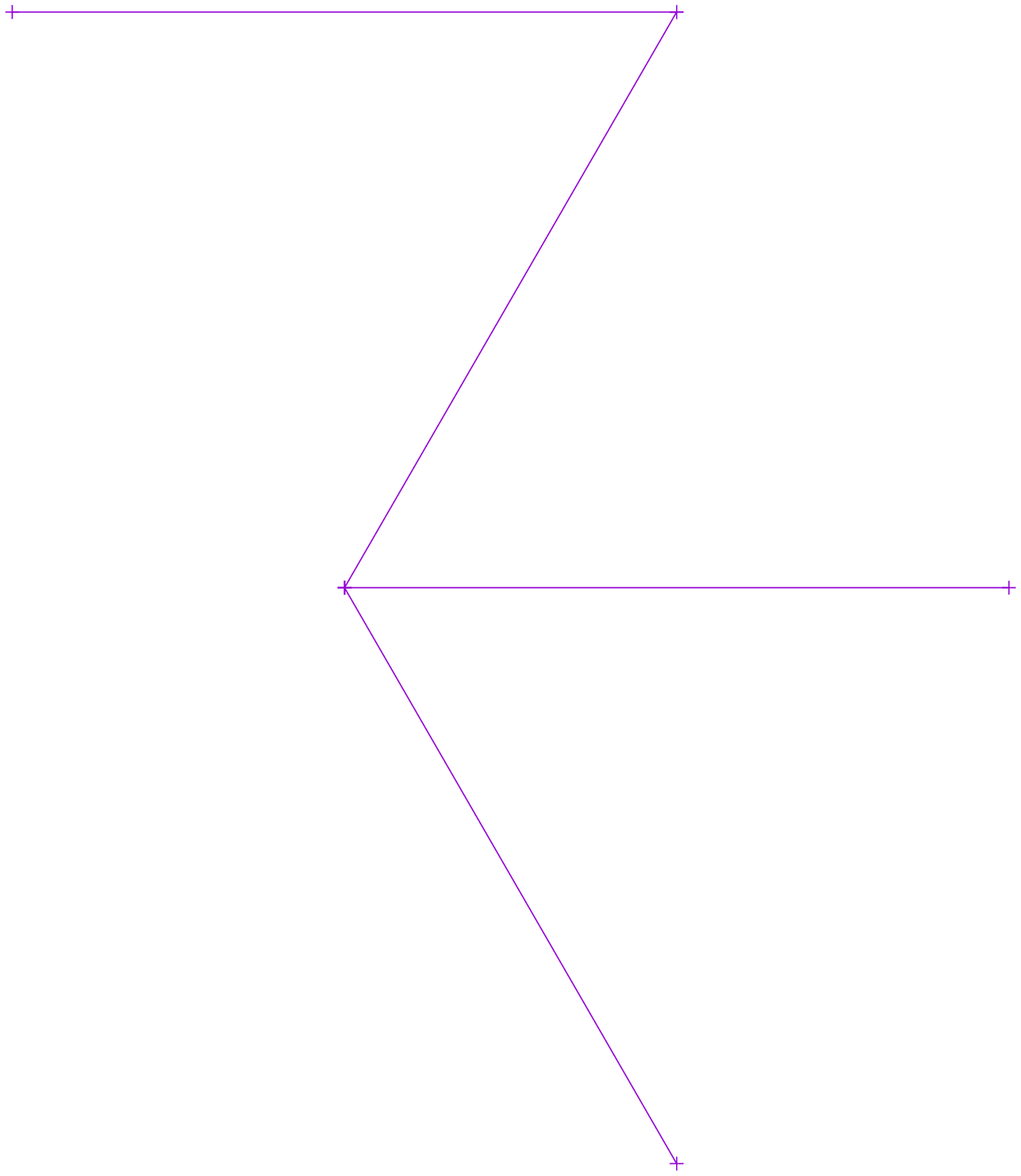})
=
V(\includegraphics[scale=0.02]{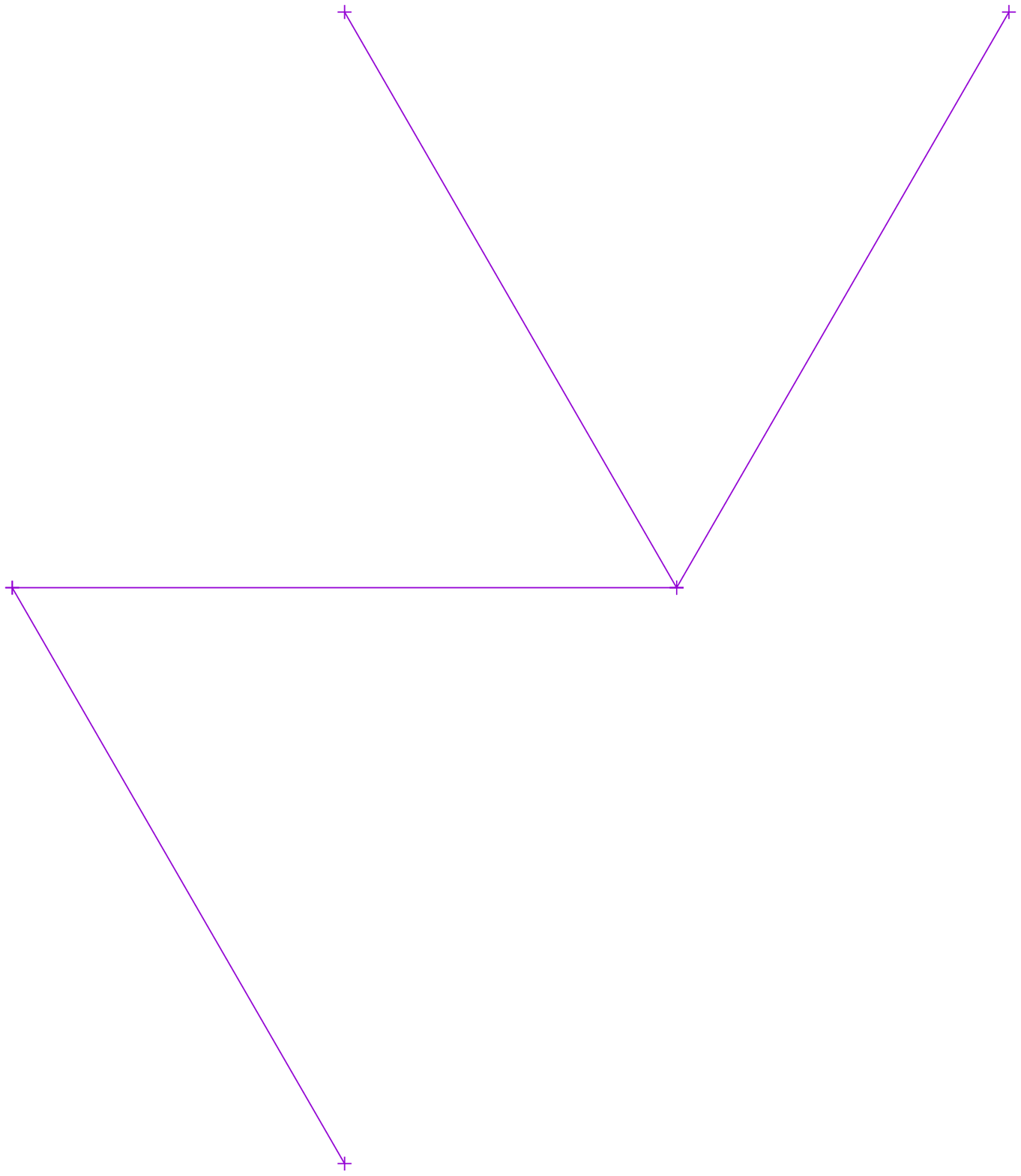})
=
V(\includegraphics[scale=0.02]{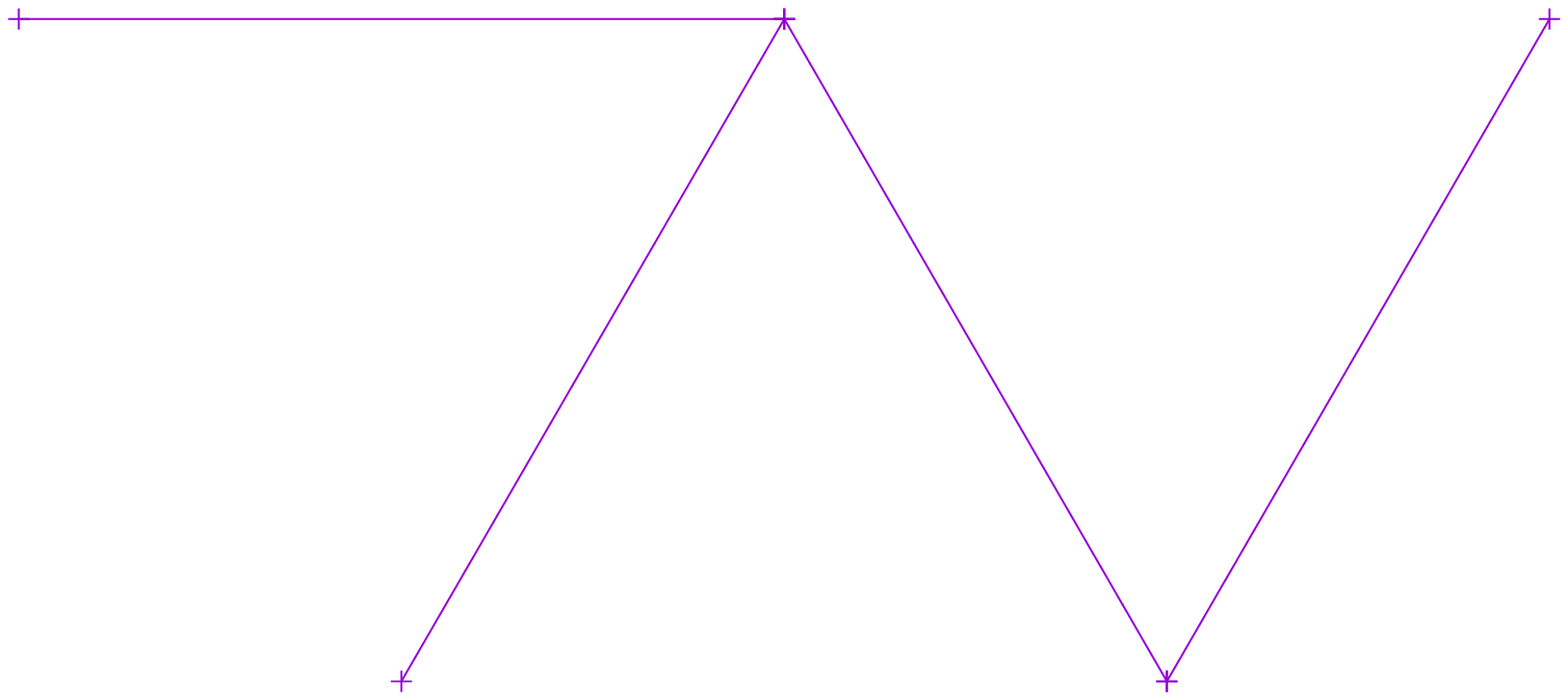})
=T_{n-2};\quad
V(\includegraphics[scale=0.02]{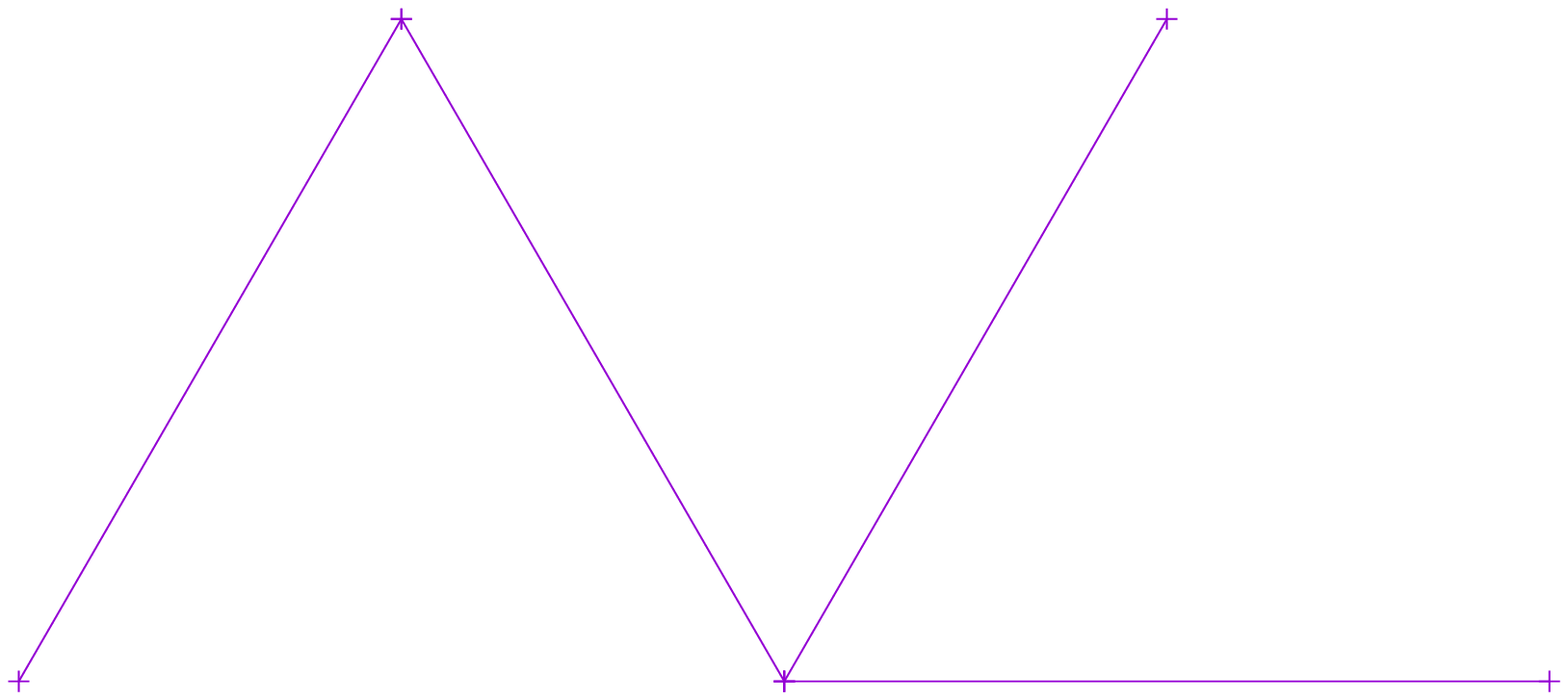})
=
V(\includegraphics[scale=0.02]{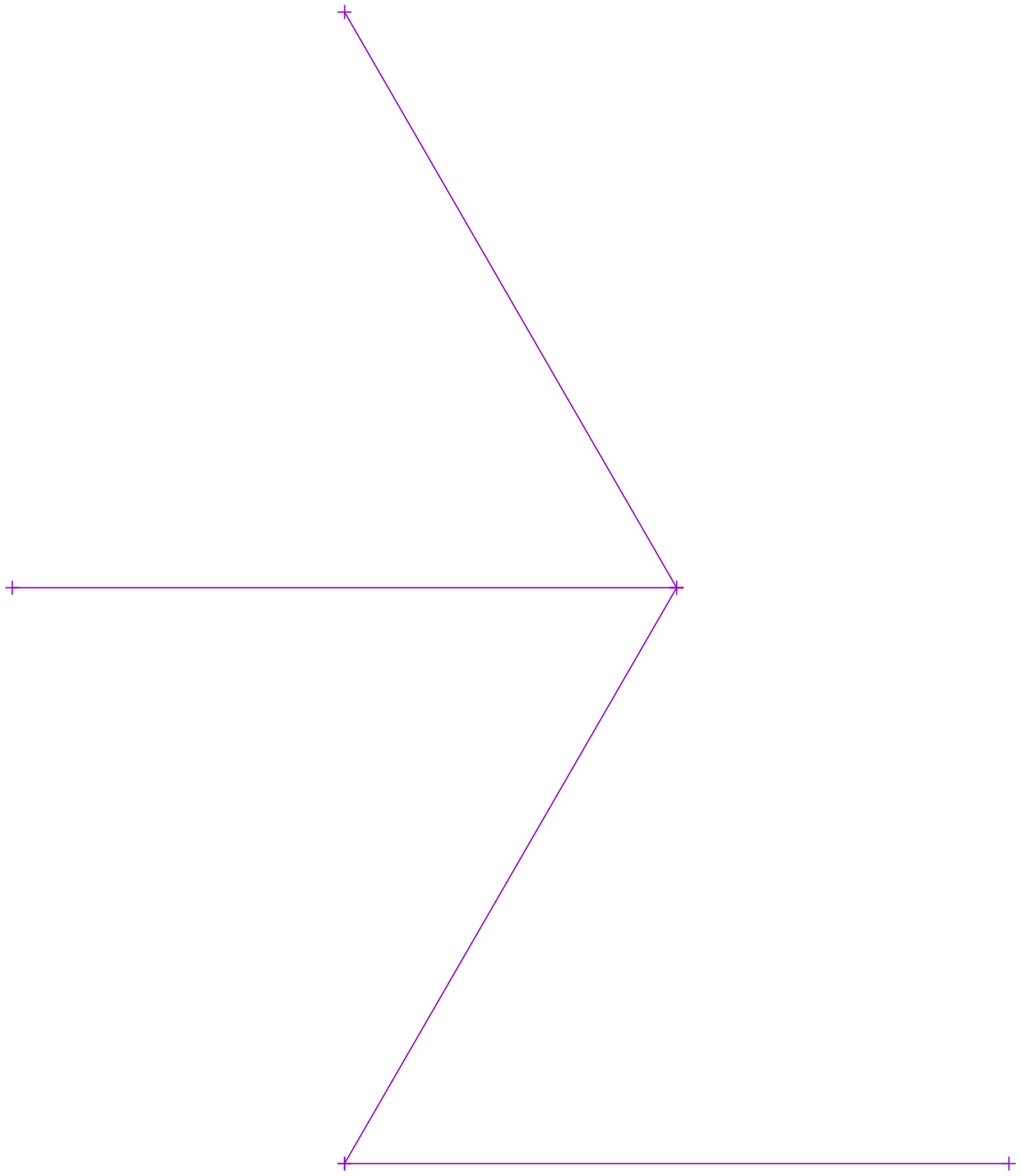})
=
V(\includegraphics[scale=0.02]{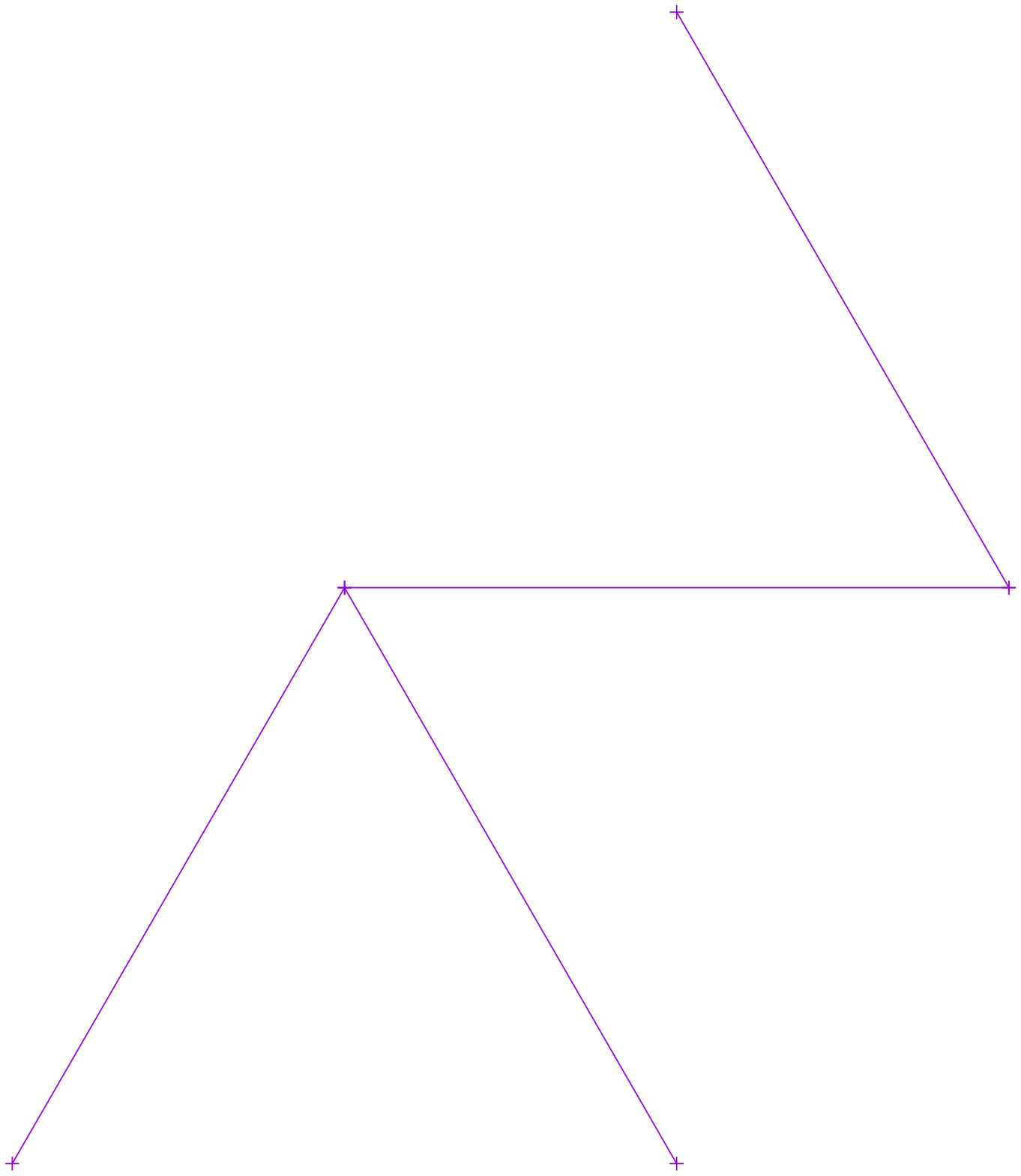})
=T_{n-3};
\end{equation}

\begin{equation}
V(\includegraphics[scale=0.02]{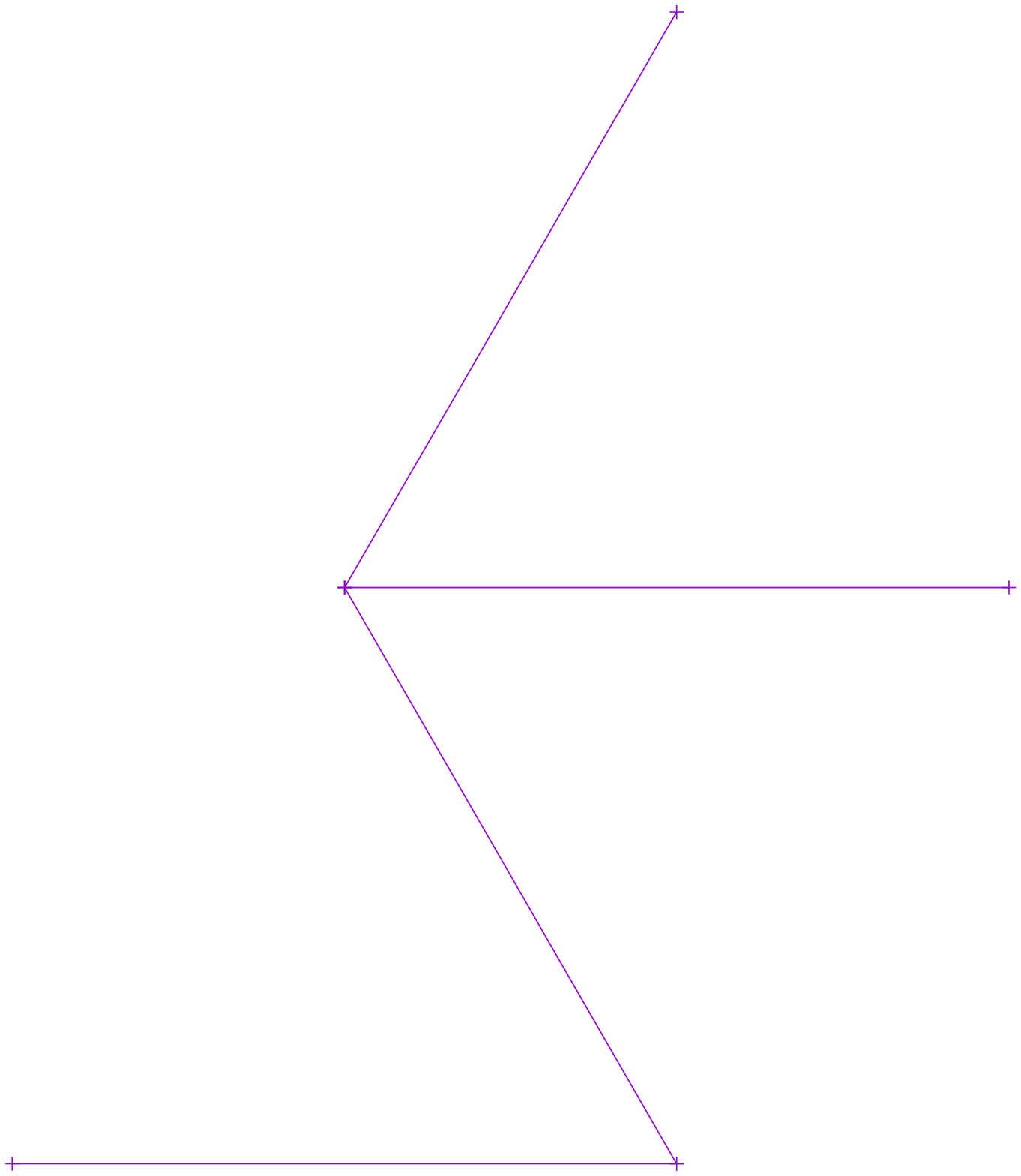})
=
V(\includegraphics[scale=0.02]{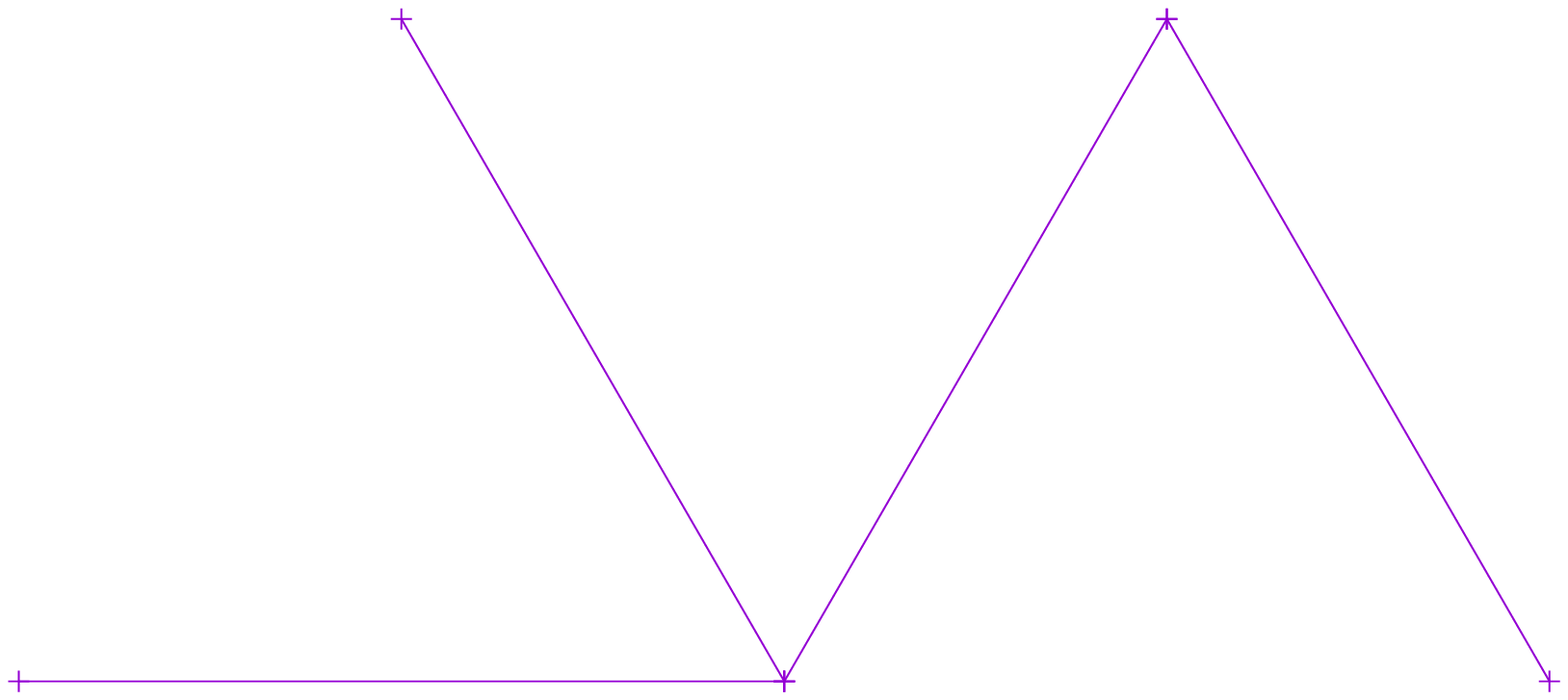})
=
V(\includegraphics[scale=0.02]{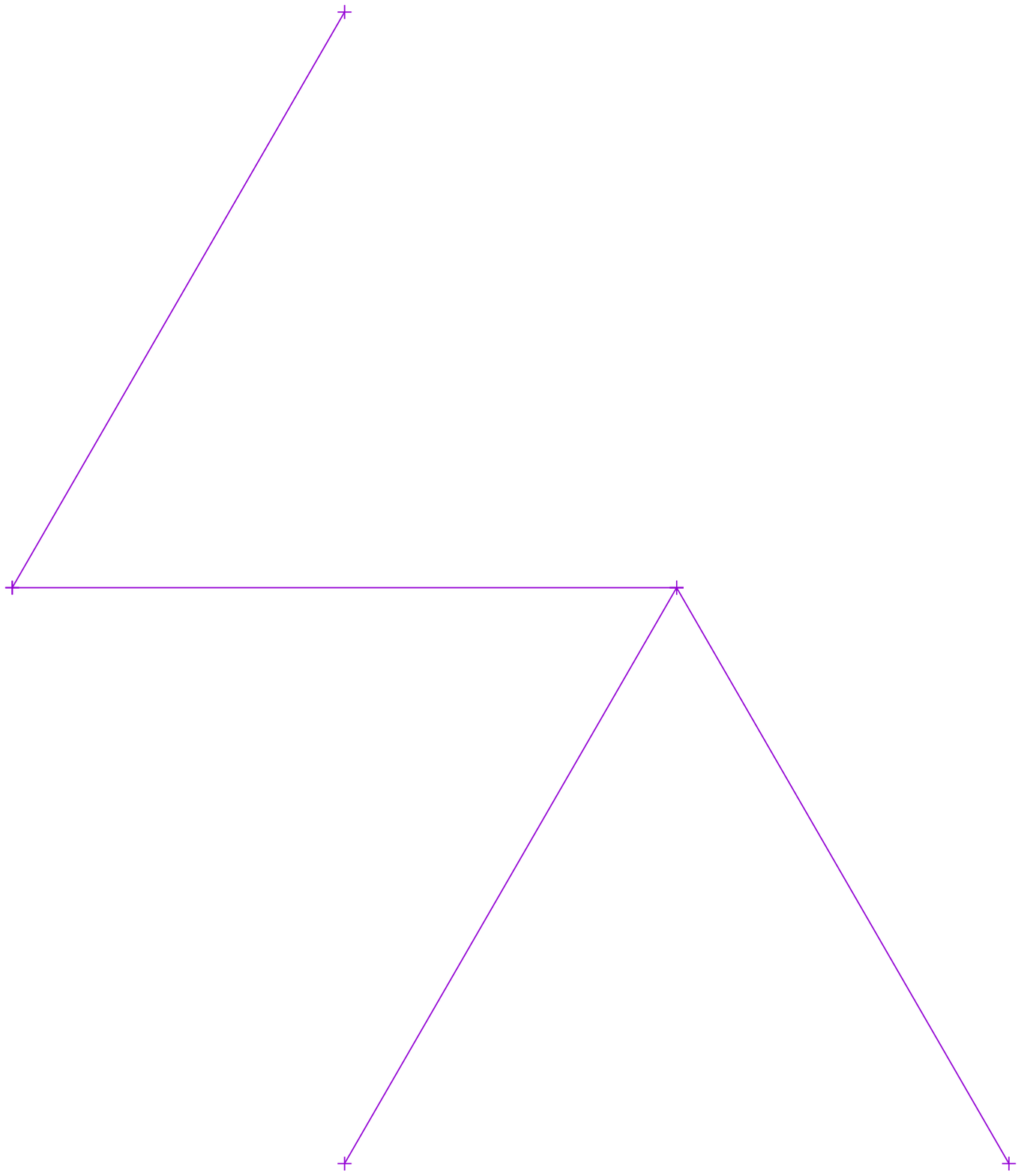})
=T_{n-3};\quad
V(\includegraphics[scale=0.02]{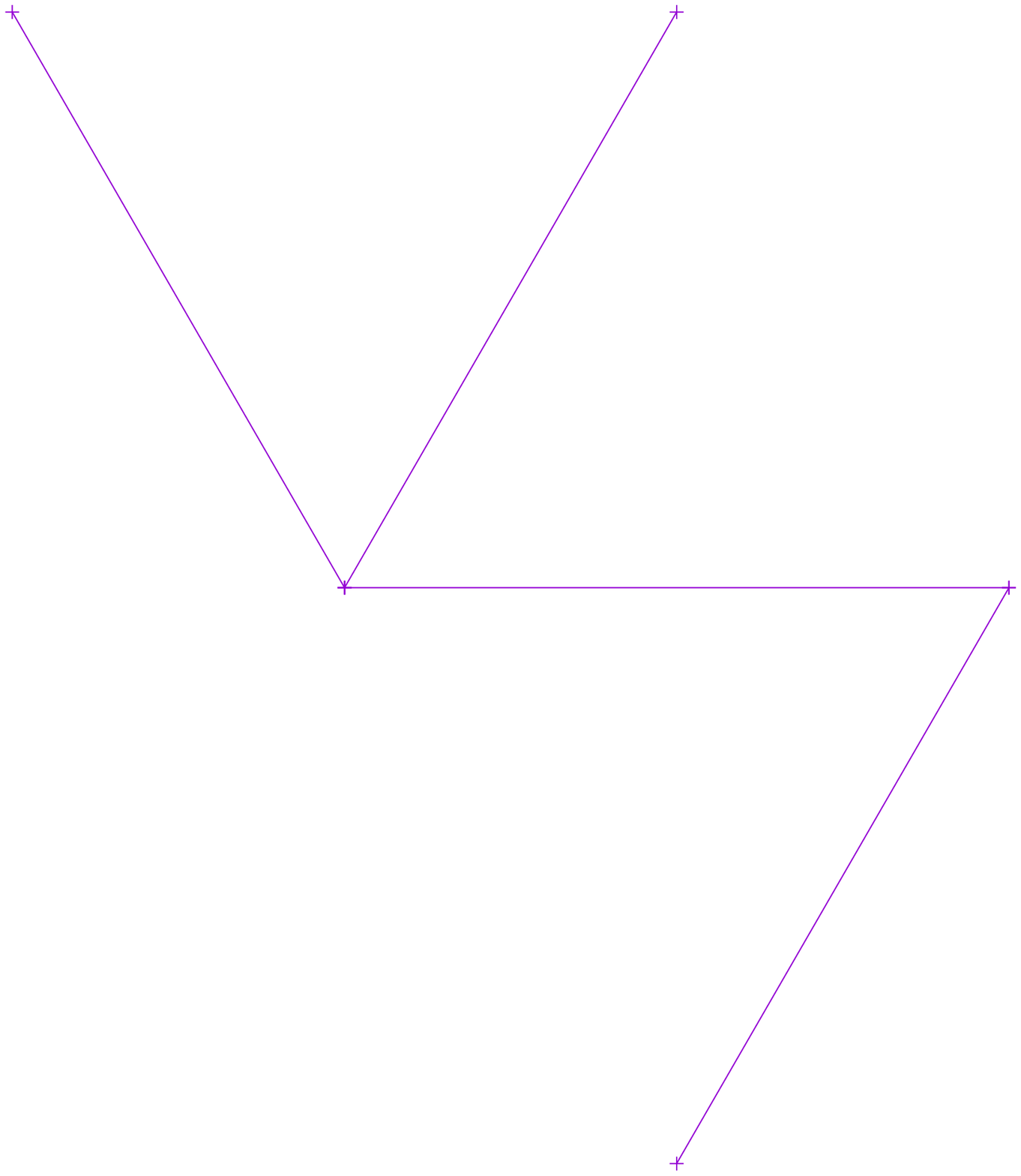})
=
V(\includegraphics[scale=0.02]{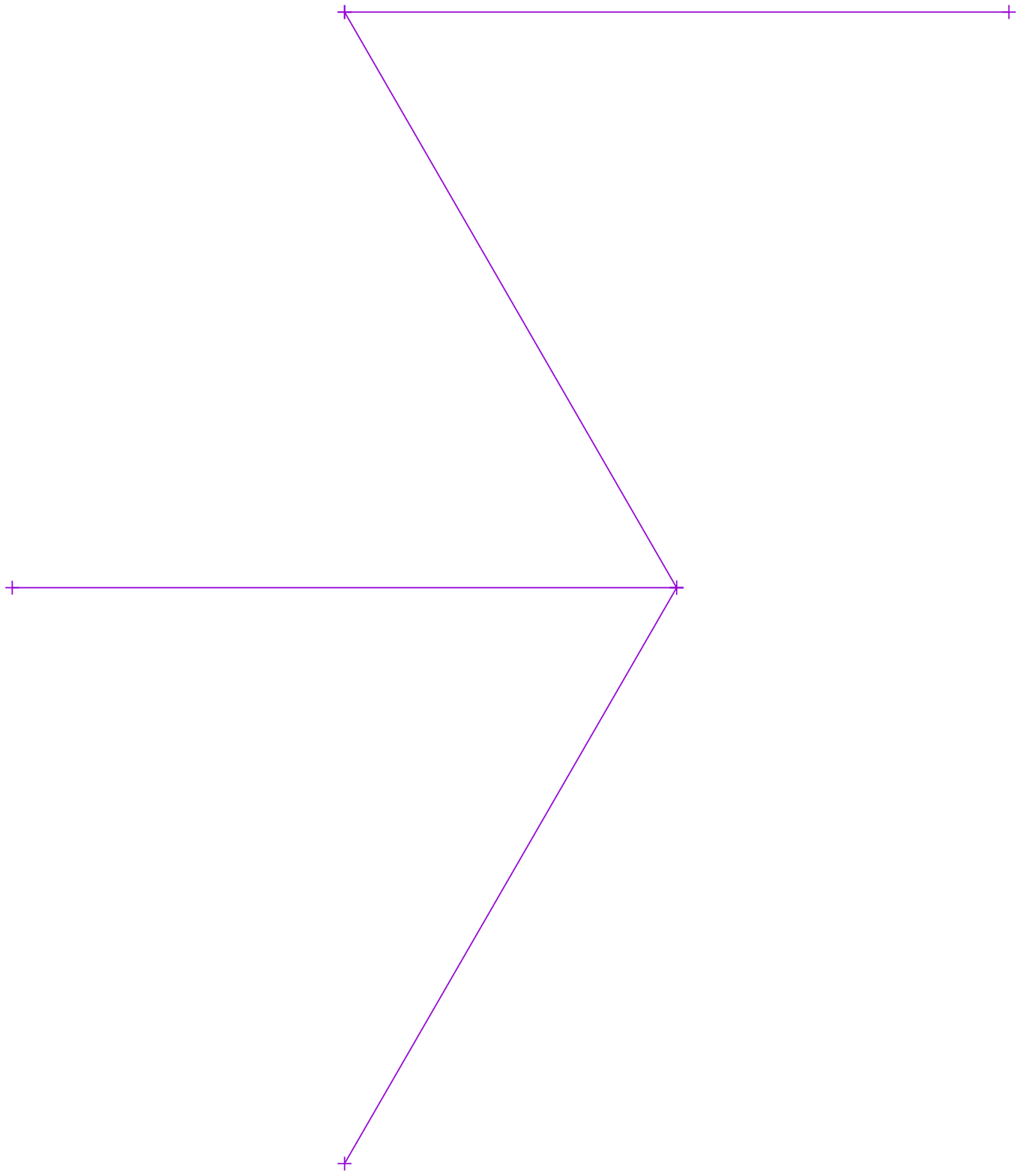})
=
V(\includegraphics[scale=0.02]{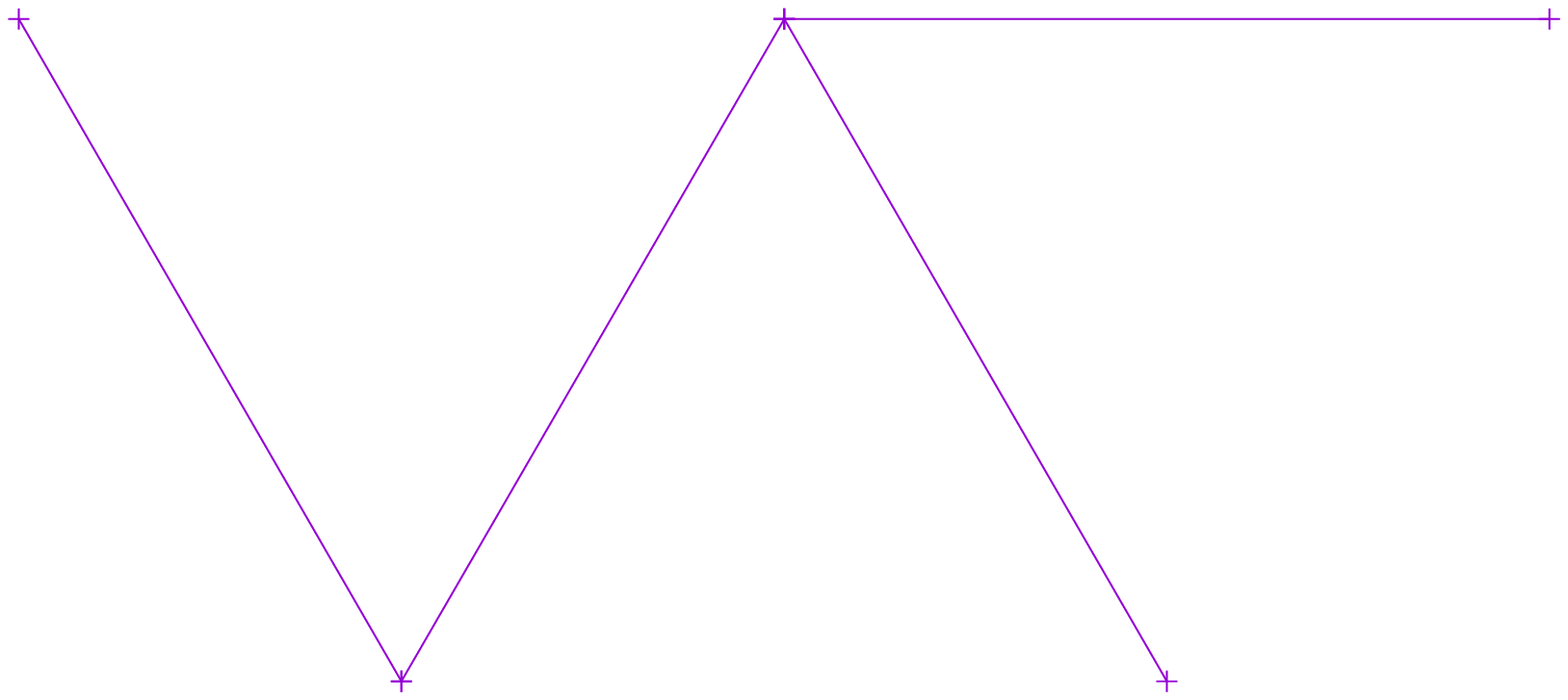})
=T_{n-2}.
\label{eq.pl4last}
\end{equation}
The sum of these 36 terms is
\begin{equation}
3(6n-11)(n-2).
\end{equation}
The Hasse diagram is sketched in Figure \ref{fig.pos4}. Groups of similar polyedges
related by symmetry are represented by a representative and a multiplier.
There is essentially
a copy of the diagram for 3 edges with a uncorrelated further edge added to each
polyedge. Then two more posets are added in rank 2 which are 
pairs of $V$'s with the same orientation (6 possible orientations)
or with different orientation ($\binom{6}{2}=15$ combinations). 
A new rank 3 with 36 
elements contains the polyedges enumerated
by eqs.\ \eqref{eq.pl4first}--\eqref{eq.pl4last}.

By selecting a $V$ and a pair of the remaining $M_{n-1}-2$ inner edges we find in rank 1
\begin{eqnarray}
\bar V(\includegraphics[scale=0.02]{pl2_1.eps}\;\includegraphics[scale=0.015]{pl1_2C.eps}\;\includegraphics[scale=0.015]{pl1_2C.eps})
&=&
V(\includegraphics[scale=0.02]{pl2_1.eps})\binom{M_{n-1}-2}{2}
-\binom{ V(\includegraphics[scale=0.02]{pl2_1.eps})}{2}
;\\
\bar V(\includegraphics[scale=0.02]{pl2_2.eps}\;\includegraphics[scale=0.015]{pl1_2C.eps}\;\includegraphics[scale=0.015]{pl1_2C.eps})
&=&
V(\includegraphics[scale=0.02]{pl2_2.eps})\binom{M_{n-1}-2}{2}
-\binom{ V(\includegraphics[scale=0.02]{pl2_2.eps})}{2}
;\\
\bar V(\includegraphics[scale=0.02]{pl2_3.eps}\;\includegraphics[scale=0.015]{pl1_2C.eps}\;\includegraphics[scale=0.015]{pl1_2C.eps})
&=&
V(\includegraphics[scale=0.02]{pl2_3.eps})\binom{M_{n-1}-2}{2}
-\binom{ V(\includegraphics[scale=0.02]{pl2_3.eps})}{2}
;\\
\bar V(\includegraphics[scale=0.02]{pl2_4.eps}\;\includegraphics[scale=0.015]{pl1_2C.eps}\;\includegraphics[scale=0.015]{pl1_2C.eps})
&=&
V(\includegraphics[scale=0.02]{pl2_4.eps})\binom{M_{n-1}-2}{2}
-\binom{ V(\includegraphics[scale=0.02]{pl2_4.eps})}{2}
;\\
\bar V(\includegraphics[scale=0.02]{pl2_5.eps}\;\includegraphics[scale=0.015]{pl1_2C.eps}\;\includegraphics[scale=0.015]{pl1_2C.eps})
&=&
V(\includegraphics[scale=0.02]{pl2_5.eps})\binom{M_{n-1}-2}{2}
-\binom{ V(\includegraphics[scale=0.02]{pl2_5.eps})}{2}
;\\
\bar V(\includegraphics[scale=0.02]{pl2_6.eps}\;\includegraphics[scale=0.015]{pl1_2C.eps}\;\includegraphics[scale=0.015]{pl1_2C.eps})
&=&
V(\includegraphics[scale=0.02]{pl2_6.eps})\binom{M_{n-1}-2}{2}
-\binom{ V(\includegraphics[scale=0.02]{pl2_6.eps})}{2}
.
\end{eqnarray}
The two green edges in the arguments of these functions indicate a pair of edges with 
any of the 6 mixed orientations.
The negative terms on the right hand sides avoid that geometries with two V's of
the same orientation are counted twice.

In a first type of rank 2 there is a polyedge with 3 edges plus one edge of the remaining $M_{n-1}-3$ inner edges:
\begin{eqnarray}
\bar V(\includegraphics[scale=0.02]{pl3_1.eps}\;\includegraphics[scale=0.015]{pl1_2C.eps})
&=&
V(\includegraphics[scale=0.02]{pl3_1.eps})(M_{n-1}-3);\\
\bar V(\includegraphics[scale=0.02]{pl3_2.eps}\;\includegraphics[scale=0.015]{pl1_2C.eps})
&=&
V(\includegraphics[scale=0.02]{pl3_2.eps})(M_{n-1}-3);\\
\bar V(\includegraphics[scale=0.02]{pl3_3.eps}\;\includegraphics[scale=0.015]{pl1_2C.eps})
&=&
V(\includegraphics[scale=0.02]{pl3_3.eps})(M_{n-1}-3);\\
\bar V(\includegraphics[scale=0.02]{pl3_4.eps}\;\includegraphics[scale=0.015]{pl1_2C.eps})
&=&
V(\includegraphics[scale=0.02]{pl3_4.eps})(M_{n-1}-3);\\
\bar V(\includegraphics[scale=0.02]{pl3_5.eps}\;\includegraphics[scale=0.015]{pl1_2C.eps})
&=&
V(\includegraphics[scale=0.02]{pl3_5.eps})(M_{n-1}-3);\\
\bar V(\includegraphics[scale=0.02]{pl3_6.eps}\;\includegraphics[scale=0.015]{pl1_2C.eps})
&=&
V(\includegraphics[scale=0.02]{pl3_6.eps})(M_{n-1}-3).
\end{eqnarray}
\begin{eqnarray}
\bar V(\includegraphics[scale=0.02]{pl3_7.eps}\;\includegraphics[scale=0.015]{pl1_2C.eps})
&=&
V(\includegraphics[scale=0.02]{pl3_7.eps})(M_{n-1}-3);\\
\bar V(\includegraphics[scale=0.02]{pl3_8.eps}\;\includegraphics[scale=0.015]{pl1_2C.eps})
&=&
V(\includegraphics[scale=0.02]{pl3_8.eps})(M_{n-1}-3);
\end{eqnarray}
\begin{eqnarray}
\bar V(\includegraphics[scale=0.02]{pl3_9.eps}\;\includegraphics[scale=0.015]{pl1_2C.eps})
&=&
V(\includegraphics[scale=0.02]{pl3_9.eps})(M_{n-1}-3);\\
\bar V(\includegraphics[scale=0.02]{pl3_10.eps}\;\includegraphics[scale=0.015]{pl1_2C.eps})
&=&
V(\includegraphics[scale=0.02]{pl3_10.eps})(M_{n-1}-3);\\
\bar V(\includegraphics[scale=0.02]{pl3_11.eps}\;\includegraphics[scale=0.015]{pl1_2C.eps})
&=&
V(\includegraphics[scale=0.02]{pl3_11.eps})(M_{n-1}-3);\\
\bar V(\includegraphics[scale=0.02]{pl3_12.eps}\;\includegraphics[scale=0.015]{pl1_2C.eps})
&=&
V(\includegraphics[scale=0.02]{pl3_12.eps})(M_{n-1}-3);\\
\bar V(\includegraphics[scale=0.02]{pl3_13.eps}\;\includegraphics[scale=0.015]{pl1_2C.eps})
&=&
V(\includegraphics[scale=0.02]{pl3_13.eps})(M_{n-1}-3);\\
\bar V(\includegraphics[scale=0.02]{pl3_14.eps}\;\includegraphics[scale=0.015]{pl1_2C.eps})
&=&
V(\includegraphics[scale=0.02]{pl3_14.eps})(M_{n-1}-3).
\end{eqnarray}
In another type of rank 2, pairs of V of the same orientation are selectible by the usual binomial count:
\begin{eqnarray}
\bar V(\includegraphics[scale=0.02]{pl2_1.eps}\;\includegraphics[scale=0.02]{pl2_1.eps})
&=&
\binom{V(\includegraphics[scale=0.02]{pl2_1.eps})}{2};\\
\bar V(\includegraphics[scale=0.02]{pl2_2.eps}\;\includegraphics[scale=0.02]{pl2_2.eps})
&=&
\binom{V(\includegraphics[scale=0.02]{pl2_2.eps})}{2};\\
\bar V(\includegraphics[scale=0.02]{pl2_3.eps}\;\includegraphics[scale=0.02]{pl2_3.eps})
&=&
\binom{V(\includegraphics[scale=0.02]{pl2_3.eps})}{2};\\
\bar V(\includegraphics[scale=0.02]{pl2_4.eps}\;\includegraphics[scale=0.02]{pl2_4.eps})
&=&
\binom{V(\includegraphics[scale=0.02]{pl2_4.eps})}{2};\\
\bar V(\includegraphics[scale=0.02]{pl2_5.eps}\;\includegraphics[scale=0.02]{pl2_5.eps})
&=&
\binom{V(\includegraphics[scale=0.02]{pl2_5.eps})}{2};\\
\bar V(\includegraphics[scale=0.02]{pl2_6.eps}\;\includegraphics[scale=0.02]{pl2_6.eps})
&=&
\binom{V(\includegraphics[scale=0.02]{pl2_6.eps})}{2}.
\end{eqnarray}
In a third type of rank 2, the 15 pairs of V with different orientations fall in two classes; 
The products on the right hand sides select two $V$'s independently,
and the one or two negative terms account for site selections
with common edges counted in Section \ref{sec.3l}:

\begin{eqnarray}
\bar V(\includegraphics[scale=0.02]{pl2_1.eps}\;\includegraphics[scale=0.02]{pl2_2.eps})
&=&
V(\includegraphics[scale=0.02]{pl2_1.eps}) V(\includegraphics[scale=0.02]{pl2_2.eps})-V(\includegraphics[scale=0.02]{pl3_14.eps})
;\\
\bar V(\includegraphics[scale=0.02]{pl2_1.eps}\;\includegraphics[scale=0.02]{pl2_3.eps})
&=&
V(\includegraphics[scale=0.02]{pl2_1.eps}) V(\includegraphics[scale=0.02]{pl2_3.eps})-V(\includegraphics[scale=0.02]{pl3_7.eps})
;\\
\bar V(\includegraphics[scale=0.02]{pl2_1.eps}\;\includegraphics[scale=0.02]{pl2_4.eps})
&=&
V(\includegraphics[scale=0.02]{pl2_1.eps}) V(\includegraphics[scale=0.02]{pl2_4.eps})-V(\includegraphics[scale=0.02]{pl3_3.eps})
-V(\includegraphics[scale=0.02]{pl3_5.eps})
;\\
\bar V(\includegraphics[scale=0.02]{pl2_1.eps}\;\includegraphics[scale=0.02]{pl2_5.eps})
&=&
V(\includegraphics[scale=0.02]{pl2_1.eps}) V(\includegraphics[scale=0.02]{pl2_5.eps})-V(\includegraphics[scale=0.02]{pl3_7.eps})
;\\
\bar V(\includegraphics[scale=0.02]{pl2_1.eps}\;\includegraphics[scale=0.02]{pl2_6.eps})
&=&
V(\includegraphics[scale=0.02]{pl2_1.eps}) V(\includegraphics[scale=0.02]{pl2_6.eps})-V(\includegraphics[scale=0.02]{pl3_13.eps})
;\\
\bar V(\includegraphics[scale=0.02]{pl2_2.eps}\;\includegraphics[scale=0.02]{pl2_3.eps})
&=&
V(\includegraphics[scale=0.02]{pl2_2.eps}) V(\includegraphics[scale=0.02]{pl2_3.eps})-V(\includegraphics[scale=0.02]{pl3_9.eps})
;\\
\bar V(\includegraphics[scale=0.02]{pl2_2.eps}\;\includegraphics[scale=0.02]{pl2_4.eps})
&=&
V(\includegraphics[scale=0.02]{pl2_2.eps}) V(\includegraphics[scale=0.02]{pl2_4.eps})-V(\includegraphics[scale=0.02]{pl3_8.eps})
;\\
\bar V(\includegraphics[scale=0.02]{pl2_2.eps}\;\includegraphics[scale=0.02]{pl2_5.eps})
&=&
V(\includegraphics[scale=0.02]{pl2_2.eps}) V(\includegraphics[scale=0.02]{pl2_5.eps})-V(\includegraphics[scale=0.02]{pl3_6.eps})
-V(\includegraphics[scale=0.02]{pl3_1.eps})
;\\
\bar V(\includegraphics[scale=0.02]{pl2_2.eps}\;\includegraphics[scale=0.02]{pl2_6.eps})
&=&
V(\includegraphics[scale=0.02]{pl2_2.eps}) V(\includegraphics[scale=0.02]{pl2_6.eps})-V(\includegraphics[scale=0.02]{pl3_8.eps})
;\\
\bar V(\includegraphics[scale=0.02]{pl2_3.eps}\;\includegraphics[scale=0.02]{pl2_4.eps})
&=&
V(\includegraphics[scale=0.02]{pl2_3.eps}) V(\includegraphics[scale=0.02]{pl2_4.eps})-V(\includegraphics[scale=0.02]{pl3_10.eps})
;\\
\bar V(\includegraphics[scale=0.02]{pl2_3.eps}\;\includegraphics[scale=0.02]{pl2_5.eps})
&=&
V(\includegraphics[scale=0.02]{pl2_3.eps}) V(\includegraphics[scale=0.02]{pl2_5.eps})-V(\includegraphics[scale=0.02]{pl3_7.eps})
;\\
\bar V(\includegraphics[scale=0.02]{pl2_3.eps}\;\includegraphics[scale=0.02]{pl2_6.eps})
&=&
V(\includegraphics[scale=0.02]{pl2_3.eps}) V(\includegraphics[scale=0.02]{pl2_6.eps})-V(\includegraphics[scale=0.02]{pl3_4.eps})
-V(\includegraphics[scale=0.02]{pl3_2.eps})
;\\
\bar V(\includegraphics[scale=0.02]{pl2_4.eps}\;\includegraphics[scale=0.02]{pl2_5.eps})
&=&
V(\includegraphics[scale=0.02]{pl2_4.eps}) V(\includegraphics[scale=0.02]{pl2_5.eps})-V(\includegraphics[scale=0.02]{pl3_11.eps})
;\\
\bar V(\includegraphics[scale=0.02]{pl2_4.eps}\;\includegraphics[scale=0.02]{pl2_6.eps})
&=&
V(\includegraphics[scale=0.02]{pl2_4.eps}) V(\includegraphics[scale=0.02]{pl2_6.eps})-V(\includegraphics[scale=0.02]{pl3_8.eps})
;\\
\bar V(\includegraphics[scale=0.02]{pl2_5.eps}\;\includegraphics[scale=0.02]{pl2_6.eps})
&=&
V(\includegraphics[scale=0.02]{pl2_5.eps}) V(\includegraphics[scale=0.02]{pl2_6.eps})-V(\includegraphics[scale=0.02]{pl3_12.eps})
.
\end{eqnarray}
Mobius inversion of the 35 graphs of rank 2 of the Hasse diagram yields:
\begin{eqnarray}
V(\includegraphics[scale=0.02]{pl3_1.eps}\;\includegraphics[scale=0.015]{pl1_2C.eps})
&=&
\bar V(\includegraphics[scale=0.02]{pl3_1.eps}\;\includegraphics[scale=0.015]{pl1_2C.eps})
-
V(\includegraphics[scale=0.02]{pl4_4.eps})
-
V(\includegraphics[scale=0.02]{pl4_10.eps})
-
V(\includegraphics[scale=0.02]{pl4_29.eps})
-
V(\includegraphics[scale=0.02]{pl4_15.eps})
-
V(\includegraphics[scale=0.02]{pl4_26.eps})
-
V(\includegraphics[scale=0.02]{pl4_2.eps})
\label{eq.mobr2first}
;\\
V(\includegraphics[scale=0.02]{pl3_2.eps}\;\includegraphics[scale=0.015]{pl1_2C.eps})
&=&
\bar V(\includegraphics[scale=0.02]{pl3_2.eps}\;\includegraphics[scale=0.015]{pl1_2C.eps})
-
V(\includegraphics[scale=0.02]{pl4_1.eps})
-
V(\includegraphics[scale=0.02]{pl4_17.eps})
-
V(\includegraphics[scale=0.02]{pl4_30.eps})
-
V(\includegraphics[scale=0.02]{pl4_11.eps})
-
V(\includegraphics[scale=0.02]{pl4_27.eps})
-
V(\includegraphics[scale=0.02]{pl4_5.eps})
;\\
V(\includegraphics[scale=0.02]{pl3_3.eps}\;\includegraphics[scale=0.015]{pl1_2C.eps})
&=&
\bar V(\includegraphics[scale=0.02]{pl3_3.eps}\;\includegraphics[scale=0.015]{pl1_2C.eps})
-
V(\includegraphics[scale=0.02]{pl4_6.eps})
-
V(\includegraphics[scale=0.02]{pl4_14.eps})
-
V(\includegraphics[scale=0.02]{pl4_28.eps})
-
V(\includegraphics[scale=0.02]{pl4_8.eps})
-
V(\includegraphics[scale=0.02]{pl4_25.eps})
-
V(\includegraphics[scale=0.02]{pl4_3.eps})
;\\
V(\includegraphics[scale=0.02]{pl3_4.eps}\;\includegraphics[scale=0.015]{pl1_2C.eps})
&=&
\bar V(\includegraphics[scale=0.02]{pl3_4.eps}\;\includegraphics[scale=0.015]{pl1_2C.eps})
-
V(\includegraphics[scale=0.02]{pl4_5.eps})
-
V(\includegraphics[scale=0.02]{pl4_7.eps})
-
V(\includegraphics[scale=0.02]{pl4_31.eps})
-
V(\includegraphics[scale=0.02]{pl4_34.eps})
-
V(\includegraphics[scale=0.02]{pl4_13.eps})
-
V(\includegraphics[scale=0.02]{pl4_1.eps})
;\\
 V(\includegraphics[scale=0.02]{pl3_5.eps}\;\includegraphics[scale=0.015]{pl1_2C.eps})
&=&
\bar V(\includegraphics[scale=0.02]{pl3_5.eps}\;\includegraphics[scale=0.015]{pl1_2C.eps})
-
V(\includegraphics[scale=0.02]{pl4_6.eps})
-
V(\includegraphics[scale=0.02]{pl4_16.eps})
-
V(\includegraphics[scale=0.02]{pl4_32.eps})
-
V(\includegraphics[scale=0.02]{pl4_35.eps})
-
V(\includegraphics[scale=0.02]{pl4_9.eps})
-
V(\includegraphics[scale=0.02]{pl4_3.eps})
;\\
V(\includegraphics[scale=0.02]{pl3_6.eps}\;\includegraphics[scale=0.015]{pl1_2C.eps})
&=&
\bar V(\includegraphics[scale=0.02]{pl3_6.eps}\;\includegraphics[scale=0.015]{pl1_2C.eps})
-
V(\includegraphics[scale=0.02]{pl4_2.eps})
-
V(\includegraphics[scale=0.02]{pl4_18.eps})
-
V(\includegraphics[scale=0.02]{pl4_33.eps})
-
V(\includegraphics[scale=0.02]{pl4_12.eps})
-
V(\includegraphics[scale=0.02]{pl4_36.eps})
-
V(\includegraphics[scale=0.02]{pl4_4.eps})
.
\end{eqnarray}
\begin{eqnarray}
V(\includegraphics[scale=0.02]{pl3_7.eps}\;\includegraphics[scale=0.015]{pl1_2C.eps})
&=&
\bar V(\includegraphics[scale=0.02]{pl3_7.eps}\;\includegraphics[scale=0.015]{pl1_2C.eps})
-
V(\includegraphics[scale=0.02]{pl4_7.eps})
-
V(\includegraphics[scale=0.02]{pl4_8.eps})
-
V(\includegraphics[scale=0.02]{pl4_9.eps})
-
V(\includegraphics[scale=0.02]{pl4_10.eps})
-
V(\includegraphics[scale=0.02]{pl4_11.eps})
-
V(\includegraphics[scale=0.02]{pl4_12.eps})
;\\
V(\includegraphics[scale=0.02]{pl3_8.eps}\;\includegraphics[scale=0.015]{pl1_2C.eps})
&=&
\bar V(\includegraphics[scale=0.02]{pl3_8.eps}\;\includegraphics[scale=0.015]{pl1_2C.eps})
-
V(\includegraphics[scale=0.02]{pl4_13.eps})
-
V(\includegraphics[scale=0.02]{pl4_14.eps})
-
V(\includegraphics[scale=0.02]{pl4_15.eps})
-
V(\includegraphics[scale=0.02]{pl4_16.eps})
-
V(\includegraphics[scale=0.02]{pl4_17.eps})
-
V(\includegraphics[scale=0.02]{pl4_18.eps})
;
\end{eqnarray}
\begin{eqnarray}
V(\includegraphics[scale=0.02]{pl3_9.eps}\;\includegraphics[scale=0.015]{pl1_2C.eps})
&=&
\bar V(\includegraphics[scale=0.02]{pl3_9.eps}\;\includegraphics[scale=0.015]{pl1_2C.eps})
-
V(\includegraphics[scale=0.02]{pl4_27.eps})
-
V(\includegraphics[scale=0.02]{pl4_10.eps})
-
V(\includegraphics[scale=0.02]{pl4_20.eps})
-
V(\includegraphics[scale=0.02]{pl4_19.eps})
-
V(\includegraphics[scale=0.02]{pl4_13.eps})
-
V(\includegraphics[scale=0.02]{pl4_33.eps})
;\\
V(\includegraphics[scale=0.02]{pl3_10.eps}\;\includegraphics[scale=0.015]{pl1_2C.eps})
&=&
\bar V(\includegraphics[scale=0.02]{pl3_10.eps}\;\includegraphics[scale=0.015]{pl1_2C.eps})
-
V(\includegraphics[scale=0.02]{pl4_28.eps})
-
V(\includegraphics[scale=0.02]{pl4_17.eps})
-
V(\includegraphics[scale=0.02]{pl4_21.eps})
-
V(\includegraphics[scale=0.02]{pl4_20.eps})
-
V(\includegraphics[scale=0.02]{pl4_9.eps})
-
V(\includegraphics[scale=0.02]{pl4_34.eps})
;\\
V(\includegraphics[scale=0.02]{pl3_11.eps}\;\includegraphics[scale=0.015]{pl1_2C.eps})
&=&
\bar V(\includegraphics[scale=0.02]{pl3_11.eps}\;\includegraphics[scale=0.015]{pl1_2C.eps})
-
V(\includegraphics[scale=0.02]{pl4_29.eps})
-
V(\includegraphics[scale=0.02]{pl4_35.eps})
-
V(\includegraphics[scale=0.02]{pl4_18.eps})
-
V(\includegraphics[scale=0.02]{pl4_8.eps})
-
V(\includegraphics[scale=0.02]{pl4_22.eps})
-
V(\includegraphics[scale=0.02]{pl4_21.eps})
;\\
V(\includegraphics[scale=0.02]{pl3_12.eps}\;\includegraphics[scale=0.015]{pl1_2C.eps})
&=&
\bar V(\includegraphics[scale=0.02]{pl3_12.eps}\;\includegraphics[scale=0.015]{pl1_2C.eps})
-
V(\includegraphics[scale=0.02]{pl4_30.eps})
-
V(\includegraphics[scale=0.02]{pl4_36.eps})
-
V(\includegraphics[scale=0.02]{pl4_7.eps})
-
V(\includegraphics[scale=0.02]{pl4_15.eps})
-
V(\includegraphics[scale=0.02]{pl4_22.eps})
-
V(\includegraphics[scale=0.02]{pl4_23.eps})
;\\
V(\includegraphics[scale=0.02]{pl3_13.eps}\;\includegraphics[scale=0.015]{pl1_2C.eps})
&=&
\bar V(\includegraphics[scale=0.02]{pl3_13.eps}\;\includegraphics[scale=0.015]{pl1_2C.eps})
-
V(\includegraphics[scale=0.02]{pl4_25.eps})
-
V(\includegraphics[scale=0.02]{pl4_31.eps})
-
V(\includegraphics[scale=0.02]{pl4_11.eps})
-
V(\includegraphics[scale=0.02]{pl4_16.eps})
-
V(\includegraphics[scale=0.02]{pl4_23.eps})
-
V(\includegraphics[scale=0.02]{pl4_24.eps})
;\\
V(\includegraphics[scale=0.02]{pl3_14.eps}\;\includegraphics[scale=0.015]{pl1_2C.eps})
&=&
\bar V(\includegraphics[scale=0.02]{pl3_14.eps}\;\includegraphics[scale=0.015]{pl1_2C.eps})
-
V(\includegraphics[scale=0.02]{pl4_26.eps})
-
V(\includegraphics[scale=0.02]{pl4_32.eps})
-
V(\includegraphics[scale=0.02]{pl4_12.eps})
-
V(\includegraphics[scale=0.02]{pl4_14.eps})
-
V(\includegraphics[scale=0.02]{pl4_19.eps})
-
V(\includegraphics[scale=0.02]{pl4_24.eps})
;
\end{eqnarray}
\begin{eqnarray}
V(\includegraphics[scale=0.02]{pl2_1.eps}\;\includegraphics[scale=0.02]{pl2_1.eps})
&=&
\bar V(\includegraphics[scale=0.02]{pl2_1.eps}\;\includegraphics[scale=0.02]{pl2_1.eps})
-
V(\includegraphics[scale=0.02]{pl4_6.eps})
;\\
V(\includegraphics[scale=0.02]{pl2_2.eps}\;\includegraphics[scale=0.02]{pl2_2.eps})
&=&
\bar V(\includegraphics[scale=0.02]{pl2_2.eps}\;\includegraphics[scale=0.02]{pl2_2.eps})
-
V(\includegraphics[scale=0.02]{pl4_4.eps})
;\\
V(\includegraphics[scale=0.02]{pl2_3.eps}\;\includegraphics[scale=0.02]{pl2_3.eps})
&=&
\bar V(\includegraphics[scale=0.02]{pl2_3.eps}\;\includegraphics[scale=0.02]{pl2_3.eps})
-
V(\includegraphics[scale=0.02]{pl4_1.eps})
;\\
V(\includegraphics[scale=0.02]{pl2_4.eps}\;\includegraphics[scale=0.02]{pl2_4.eps})
&=&
\bar V(\includegraphics[scale=0.02]{pl2_4.eps}\;\includegraphics[scale=0.02]{pl2_4.eps})
-
V(\includegraphics[scale=0.02]{pl4_3.eps})
;\\
V(\includegraphics[scale=0.02]{pl2_5.eps}\;\includegraphics[scale=0.02]{pl2_5.eps})
&=&
\bar V(\includegraphics[scale=0.02]{pl2_5.eps}\;\includegraphics[scale=0.02]{pl2_5.eps})
-
V(\includegraphics[scale=0.02]{pl4_2.eps}) 
;\\
V(\includegraphics[scale=0.02]{pl2_6.eps}\;\includegraphics[scale=0.02]{pl2_6.eps})
&=&
\bar V(\includegraphics[scale=0.02]{pl2_6.eps}\;\includegraphics[scale=0.02]{pl2_6.eps})
-
V(\includegraphics[scale=0.02]{pl4_5.eps})
;
\end{eqnarray}
\begin{eqnarray}
V(\includegraphics[scale=0.02]{pl2_1.eps}\;\includegraphics[scale=0.02]{pl2_2.eps})
&=&
\bar V(\includegraphics[scale=0.02]{pl2_1.eps}\;\includegraphics[scale=0.02]{pl2_2.eps})
-
V(\includegraphics[scale=0.02]{pl4_10.eps})
-
V(\includegraphics[scale=0.02]{pl4_16.eps})
;\\
V(\includegraphics[scale=0.02]{pl2_1.eps}\;\includegraphics[scale=0.02]{pl2_3.eps})
&=&
\bar V(\includegraphics[scale=0.02]{pl2_1.eps}\;\includegraphics[scale=0.02]{pl2_3.eps})
-
V(\includegraphics[scale=0.02]{pl4_31.eps})
-
V(\includegraphics[scale=0.02]{pl4_19.eps})
-
V(\includegraphics[scale=0.02]{pl4_28.eps})
;\\
V(\includegraphics[scale=0.02]{pl2_1.eps}\;\includegraphics[scale=0.02]{pl2_4.eps})
&=&
\bar V(\includegraphics[scale=0.02]{pl2_1.eps}\;\includegraphics[scale=0.02]{pl2_4.eps})
;\\
V(\includegraphics[scale=0.02]{pl2_1.eps}\;\includegraphics[scale=0.02]{pl2_5.eps})
&=&
\bar V(\includegraphics[scale=0.02]{pl2_1.eps}\;\includegraphics[scale=0.02]{pl2_5.eps})
-
V(\includegraphics[scale=0.02]{pl4_23.eps})
-
V(\includegraphics[scale=0.02]{pl4_35.eps})
-
V(\includegraphics[scale=0.02]{pl4_26.eps})
;\\
V(\includegraphics[scale=0.02]{pl2_1.eps}\;\includegraphics[scale=0.02]{pl2_6.eps})
&=&
\bar V(\includegraphics[scale=0.02]{pl2_1.eps}\;\includegraphics[scale=0.02]{pl2_6.eps})
-
V(\includegraphics[scale=0.02]{pl4_7.eps})
-
V(\includegraphics[scale=0.02]{pl4_14.eps})
;\\
V(\includegraphics[scale=0.02]{pl2_2.eps}\;\includegraphics[scale=0.02]{pl2_3.eps})
&=&
\bar V(\includegraphics[scale=0.02]{pl2_2.eps}\;\includegraphics[scale=0.02]{pl2_3.eps})
-
V(\includegraphics[scale=0.02]{pl4_12.eps})
-
V(\includegraphics[scale=0.02]{pl4_17.eps})
;\\
V(\includegraphics[scale=0.02]{pl2_2.eps}\;\includegraphics[scale=0.02]{pl2_4.eps})
&=&
\bar V(\includegraphics[scale=0.02]{pl2_2.eps}\;\includegraphics[scale=0.02]{pl2_4.eps})
-
V(\includegraphics[scale=0.02]{pl4_20.eps})
-
V(\includegraphics[scale=0.02]{pl4_32.eps})
-
V(\includegraphics[scale=0.02]{pl4_29.eps}) 
;\\
V(\includegraphics[scale=0.02]{pl2_2.eps}\;\includegraphics[scale=0.02]{pl2_5.eps})
&=&
\bar V(\includegraphics[scale=0.02]{pl2_2.eps}\;\includegraphics[scale=0.02]{pl2_5.eps})
;\\
V(\includegraphics[scale=0.02]{pl2_2.eps}\;\includegraphics[scale=0.02]{pl2_6.eps})
&=&
\bar V(\includegraphics[scale=0.02]{pl2_2.eps}\;\includegraphics[scale=0.02]{pl2_6.eps})
-
V(\includegraphics[scale=0.02]{pl4_24.eps})
-
V(\includegraphics[scale=0.02]{pl4_36.eps})
-
V(\includegraphics[scale=0.02]{pl4_27.eps})  
;\\
V(\includegraphics[scale=0.02]{pl2_3.eps}\;\includegraphics[scale=0.02]{pl2_4.eps})
&=&
\bar V(\includegraphics[scale=0.02]{pl2_3.eps}\;\includegraphics[scale=0.02]{pl2_4.eps})
-
V(\includegraphics[scale=0.02]{pl4_8.eps})
-
V(\includegraphics[scale=0.02]{pl4_13.eps})
;\\
V(\includegraphics[scale=0.02]{pl2_3.eps}\;\includegraphics[scale=0.02]{pl2_5.eps})
&=&
\bar V(\includegraphics[scale=0.02]{pl2_3.eps}\;\includegraphics[scale=0.02]{pl2_5.eps})
-
V(\includegraphics[scale=0.02]{pl4_21.eps})
-
V(\includegraphics[scale=0.02]{pl4_30.eps})
-
V(\includegraphics[scale=0.02]{pl4_33.eps})
;\\
V(\includegraphics[scale=0.02]{pl2_3.eps}\;\includegraphics[scale=0.02]{pl2_6.eps})
&=&
\bar V(\includegraphics[scale=0.02]{pl2_3.eps}\;\includegraphics[scale=0.02]{pl2_6.eps})
;\\
V(\includegraphics[scale=0.02]{pl2_4.eps}\;\includegraphics[scale=0.02]{pl2_5.eps})
&=&
\bar V(\includegraphics[scale=0.02]{pl2_4.eps}\;\includegraphics[scale=0.02]{pl2_5.eps})
-
V(\includegraphics[scale=0.02]{pl4_9.eps})
-
V(\includegraphics[scale=0.02]{pl4_15.eps})
;\\
V(\includegraphics[scale=0.02]{pl2_4.eps}\;\includegraphics[scale=0.02]{pl2_6.eps})
&=&
\bar V(\includegraphics[scale=0.02]{pl2_4.eps}\;\includegraphics[scale=0.02]{pl2_6.eps})
-
V(\includegraphics[scale=0.02]{pl4_22.eps})
-
V(\includegraphics[scale=0.02]{pl4_25.eps})
-
V(\includegraphics[scale=0.02]{pl4_34.eps})
;\\
V(\includegraphics[scale=0.02]{pl2_5.eps}\;\includegraphics[scale=0.02]{pl2_6.eps})
&=&
\bar V(\includegraphics[scale=0.02]{pl2_5.eps}\;\includegraphics[scale=0.02]{pl2_6.eps})
-
V(\includegraphics[scale=0.02]{pl4_11.eps})
-
V(\includegraphics[scale=0.02]{pl4_18.eps})
.
\label{eq.mobr2last}
\end{eqnarray}
Each $V$-term on the right hand side of equations \eqref{eq.mobr2first}--\eqref{eq.mobr2last} is
equivalent to an arrow between rank 2 and rank 3 diagrams (120 poset covers).
The sum of these 35 terms is
\begin{equation}
\frac{3}{2}(n-2)(10n^3-20n^2-47n+95).
\end{equation}

Mobius inversion of the 6 sets at rank 1:
\begin{eqnarray*}
V(\includegraphics[scale=0.02]{pl2_1.eps}\;\includegraphics[scale=0.015]{pl1_2C.eps}\;\includegraphics[scale=0.015]{pl1_2C.eps})
&=&
\bar V(\includegraphics[scale=0.02]{pl2_1.eps}\;\includegraphics[scale=0.015]{pl1_2C.eps}\;\includegraphics[scale=0.015]{pl1_2C.eps})
-
\bar V(\includegraphics[scale=0.02]{pl3_3.eps}\;\includegraphics[scale=0.015]{pl1_2C.eps})
-
\bar V(\includegraphics[scale=0.02]{pl3_5.eps}\;\includegraphics[scale=0.015]{pl1_2C.eps})
-
\bar V(\includegraphics[scale=0.02]{pl3_7.eps}\;\includegraphics[scale=0.015]{pl1_2C.eps})
-
\bar V(\includegraphics[scale=0.02]{pl3_13.eps}\;\includegraphics[scale=0.015]{pl1_2C.eps})
-
\bar V(\includegraphics[scale=0.02]{pl3_14.eps}\;\includegraphics[scale=0.015]{pl1_2C.eps})
\\
&&
-
\bar V(\includegraphics[scale=0.02]{pl2_1.eps}\;\includegraphics[scale=0.02]{pl2_1.eps})
-
\bar V(\includegraphics[scale=0.02]{pl2_1.eps}\;\includegraphics[scale=0.02]{pl2_2.eps})
-
\bar V(\includegraphics[scale=0.02]{pl2_1.eps}\;\includegraphics[scale=0.02]{pl2_3.eps})
-
\bar V(\includegraphics[scale=0.02]{pl2_1.eps}\;\includegraphics[scale=0.02]{pl2_4.eps})
-
\bar V(\includegraphics[scale=0.02]{pl2_1.eps}\;\includegraphics[scale=0.02]{pl2_5.eps})
-
\bar V(\includegraphics[scale=0.02]{pl2_1.eps}\;\includegraphics[scale=0.02]{pl2_6.eps})
\\&&
+ V(\includegraphics[scale=0.02]{pl4_3.eps})
+2 V(\includegraphics[scale=0.02]{pl4_6.eps})
+ V(\includegraphics[scale=0.02]{pl4_7.eps})
+ V(\includegraphics[scale=0.02]{pl4_8.eps})
+ V(\includegraphics[scale=0.02]{pl4_9.eps})
+ V(\includegraphics[scale=0.02]{pl4_10.eps})
\\&&
+ V(\includegraphics[scale=0.02]{pl4_11.eps})
+ V(\includegraphics[scale=0.02]{pl4_12.eps})
+2 V(\includegraphics[scale=0.02]{pl4_14.eps})
+2 V(\includegraphics[scale=0.02]{pl4_16.eps})
+ V(\includegraphics[scale=0.02]{pl4_19.eps})
+ V(\includegraphics[scale=0.02]{pl4_23.eps})
\\&&
+ V(\includegraphics[scale=0.02]{pl4_24.eps})
+ V(\includegraphics[scale=0.02]{pl4_25.eps})
+ V(\includegraphics[scale=0.02]{pl4_26.eps})
+ V(\includegraphics[scale=0.02]{pl4_28.eps})
+ V(\includegraphics[scale=0.02]{pl4_31.eps})
+ V(\includegraphics[scale=0.02]{pl4_32.eps})
+ V(\includegraphics[scale=0.02]{pl4_35.eps})
;\\
V(\includegraphics[scale=0.02]{pl2_2.eps}\;\includegraphics[scale=0.015]{pl1_2C.eps}\;\includegraphics[scale=0.015]{pl1_2C.eps})
&=&
\bar V(\includegraphics[scale=0.02]{pl2_2.eps}\;\includegraphics[scale=0.015]{pl1_2C.eps}\;\includegraphics[scale=0.015]{pl1_2C.eps})
-
\bar V(\includegraphics[scale=0.02]{pl3_1.eps}\;\includegraphics[scale=0.015]{pl1_2C.eps})
-
\bar V(\includegraphics[scale=0.02]{pl3_6.eps}\;\includegraphics[scale=0.015]{pl1_2C.eps})
-
\bar V(\includegraphics[scale=0.02]{pl3_8.eps}\;\includegraphics[scale=0.015]{pl1_2C.eps})
-
\bar V(\includegraphics[scale=0.02]{pl3_9.eps}\;\includegraphics[scale=0.015]{pl1_2C.eps})
-
\bar V(\includegraphics[scale=0.02]{pl3_14.eps}\;\includegraphics[scale=0.015]{pl1_2C.eps})
\\
&&
-
\bar V(\includegraphics[scale=0.02]{pl2_1.eps}\;\includegraphics[scale=0.02]{pl2_2.eps})
-
\bar V(\includegraphics[scale=0.02]{pl2_2.eps}\;\includegraphics[scale=0.02]{pl2_2.eps})
-
\bar V(\includegraphics[scale=0.02]{pl2_2.eps}\;\includegraphics[scale=0.02]{pl2_3.eps})
-
\bar V(\includegraphics[scale=0.02]{pl2_2.eps}\;\includegraphics[scale=0.02]{pl2_4.eps})
-
\bar V(\includegraphics[scale=0.02]{pl2_2.eps}\;\includegraphics[scale=0.02]{pl2_5.eps})
-
\bar V(\includegraphics[scale=0.02]{pl2_2.eps}\;\includegraphics[scale=0.02]{pl2_6.eps})
\\&&
+ V(\includegraphics[scale=0.02]{pl4_2.eps})
+2 V(\includegraphics[scale=0.02]{pl4_4.eps})
+2 V(\includegraphics[scale=0.02]{pl4_10.eps})
+2 V(\includegraphics[scale=0.02]{pl4_12.eps})
+ V(\includegraphics[scale=0.02]{pl4_13.eps})
+ V(\includegraphics[scale=0.02]{pl4_14.eps})
\\&&
+ V(\includegraphics[scale=0.02]{pl4_15.eps})
+ V(\includegraphics[scale=0.02]{pl4_16.eps})
+ V(\includegraphics[scale=0.02]{pl4_17.eps})
+ V(\includegraphics[scale=0.02]{pl4_18.eps})
+ V(\includegraphics[scale=0.02]{pl4_19.eps})
+ V(\includegraphics[scale=0.02]{pl4_20.eps})
\\&&
+ V(\includegraphics[scale=0.02]{pl4_24.eps})
+ V(\includegraphics[scale=0.02]{pl4_26.eps})
+ V(\includegraphics[scale=0.02]{pl4_27.eps})
+ V(\includegraphics[scale=0.02]{pl4_29.eps})
+ V(\includegraphics[scale=0.02]{pl4_32.eps})
+ V(\includegraphics[scale=0.02]{pl4_33.eps})
+ V(\includegraphics[scale=0.02]{pl4_36.eps})
;\\
V(\includegraphics[scale=0.02]{pl2_3.eps}\;\includegraphics[scale=0.015]{pl1_2C.eps}\;\includegraphics[scale=0.015]{pl1_2C.eps})
&=&
\bar V(\includegraphics[scale=0.02]{pl2_3.eps}\;\includegraphics[scale=0.015]{pl1_2C.eps}\;\includegraphics[scale=0.015]{pl1_2C.eps})
-
\bar V(\includegraphics[scale=0.02]{pl3_2.eps}\;\includegraphics[scale=0.015]{pl1_2C.eps})
-
\bar V(\includegraphics[scale=0.02]{pl3_4.eps}\;\includegraphics[scale=0.015]{pl1_2C.eps})
-
\bar V(\includegraphics[scale=0.02]{pl3_7.eps}\;\includegraphics[scale=0.015]{pl1_2C.eps})
-
\bar V(\includegraphics[scale=0.02]{pl3_9.eps}\;\includegraphics[scale=0.015]{pl1_2C.eps})
-
\bar V(\includegraphics[scale=0.02]{pl3_10.eps}\;\includegraphics[scale=0.015]{pl1_2C.eps})
\\
&&
-
\bar V(\includegraphics[scale=0.02]{pl2_1.eps}\;\includegraphics[scale=0.02]{pl2_3.eps})
-
\bar V(\includegraphics[scale=0.02]{pl2_2.eps}\;\includegraphics[scale=0.02]{pl2_3.eps})
-
\bar V(\includegraphics[scale=0.02]{pl2_3.eps}\;\includegraphics[scale=0.02]{pl2_3.eps})
-
\bar V(\includegraphics[scale=0.02]{pl2_3.eps}\;\includegraphics[scale=0.02]{pl2_4.eps})
-
\bar V(\includegraphics[scale=0.02]{pl2_3.eps}\;\includegraphics[scale=0.02]{pl2_5.eps})
-
\bar V(\includegraphics[scale=0.02]{pl2_3.eps}\;\includegraphics[scale=0.02]{pl2_6.eps})
\\&&
+2 V(\includegraphics[scale=0.02]{pl4_1.eps})
+ V(\includegraphics[scale=0.02]{pl4_5.eps})
+ V(\includegraphics[scale=0.02]{pl4_7.eps})
+ V(\includegraphics[scale=0.02]{pl4_8.eps})
+ V(\includegraphics[scale=0.02]{pl4_9.eps})
+ V(\includegraphics[scale=0.02]{pl4_10.eps})
\\&&
+ V(\includegraphics[scale=0.02]{pl4_11.eps})
+ V(\includegraphics[scale=0.02]{pl4_12.eps})
+2 V(\includegraphics[scale=0.02]{pl4_13.eps})
+2 V(\includegraphics[scale=0.02]{pl4_17.eps})
+ V(\includegraphics[scale=0.02]{pl4_19.eps})
+ V(\includegraphics[scale=0.02]{pl4_20.eps})
\\&&
+ V(\includegraphics[scale=0.02]{pl4_21.eps})
+ V(\includegraphics[scale=0.02]{pl4_27.eps})
+ V(\includegraphics[scale=0.02]{pl4_28.eps})
+ V(\includegraphics[scale=0.02]{pl4_30.eps})
+ V(\includegraphics[scale=0.02]{pl4_31.eps})
+ V(\includegraphics[scale=0.02]{pl4_33.eps})
+ V(\includegraphics[scale=0.02]{pl4_34.eps})
;\\
V(\includegraphics[scale=0.02]{pl2_4.eps}\;\includegraphics[scale=0.015]{pl1_2C.eps}\;\includegraphics[scale=0.015]{pl1_2C.eps})
&=&
\bar V(\includegraphics[scale=0.02]{pl2_4.eps}\;\includegraphics[scale=0.015]{pl1_2C.eps}\;\includegraphics[scale=0.015]{pl1_2C.eps})
-
\bar V(\includegraphics[scale=0.02]{pl3_3.eps}\;\includegraphics[scale=0.015]{pl1_2C.eps})
-
\bar V(\includegraphics[scale=0.02]{pl3_5.eps}\;\includegraphics[scale=0.015]{pl1_2C.eps})
-
\bar V(\includegraphics[scale=0.02]{pl3_8.eps}\;\includegraphics[scale=0.015]{pl1_2C.eps})
-
\bar V(\includegraphics[scale=0.02]{pl3_10.eps}\;\includegraphics[scale=0.015]{pl1_2C.eps})
-
\bar V(\includegraphics[scale=0.02]{pl3_11.eps}\;\includegraphics[scale=0.015]{pl1_2C.eps})
\\
&&
-
\bar V(\includegraphics[scale=0.02]{pl2_1.eps}\;\includegraphics[scale=0.02]{pl2_4.eps})
-
\bar V(\includegraphics[scale=0.02]{pl2_2.eps}\;\includegraphics[scale=0.02]{pl2_4.eps})
-
\bar V(\includegraphics[scale=0.02]{pl2_3.eps}\;\includegraphics[scale=0.02]{pl2_4.eps})
-
\bar V(\includegraphics[scale=0.02]{pl2_4.eps}\;\includegraphics[scale=0.02]{pl2_4.eps})
-
\bar V(\includegraphics[scale=0.02]{pl2_4.eps}\;\includegraphics[scale=0.02]{pl2_5.eps})
-
\bar V(\includegraphics[scale=0.02]{pl2_4.eps}\;\includegraphics[scale=0.02]{pl2_6.eps})
\\&&
+2 V(\includegraphics[scale=0.02]{pl4_3.eps})
+ V(\includegraphics[scale=0.02]{pl4_6.eps})
+2 V(\includegraphics[scale=0.02]{pl4_8.eps})
+2 V(\includegraphics[scale=0.02]{pl4_9.eps})
+ V(\includegraphics[scale=0.02]{pl4_13.eps})
+ V(\includegraphics[scale=0.02]{pl4_14.eps})
\\&&
+ V(\includegraphics[scale=0.02]{pl4_15.eps})
+ V(\includegraphics[scale=0.02]{pl4_16.eps})
+ V(\includegraphics[scale=0.02]{pl4_17.eps})
+ V(\includegraphics[scale=0.02]{pl4_18.eps})
+ V(\includegraphics[scale=0.02]{pl4_20.eps})
+ V(\includegraphics[scale=0.02]{pl4_21.eps})
\\&&
+ V(\includegraphics[scale=0.02]{pl4_22.eps})
+ V(\includegraphics[scale=0.02]{pl4_25.eps})
+ V(\includegraphics[scale=0.02]{pl4_28.eps})
+ V(\includegraphics[scale=0.02]{pl4_29.eps})
+ V(\includegraphics[scale=0.02]{pl4_32.eps})
+ V(\includegraphics[scale=0.02]{pl4_34.eps})
+ V(\includegraphics[scale=0.02]{pl4_35.eps})
;\\
V(\includegraphics[scale=0.02]{pl2_5.eps}\;\includegraphics[scale=0.015]{pl1_2C.eps}\;\includegraphics[scale=0.015]{pl1_2C.eps})
&=&
\bar V(\includegraphics[scale=0.02]{pl2_5.eps}\;\includegraphics[scale=0.015]{pl1_2C.eps}\;\includegraphics[scale=0.015]{pl1_2C.eps})
-
\bar V(\includegraphics[scale=0.02]{pl3_1.eps}\;\includegraphics[scale=0.015]{pl1_2C.eps})
-
\bar V(\includegraphics[scale=0.02]{pl3_6.eps}\;\includegraphics[scale=0.015]{pl1_2C.eps})
-
\bar V(\includegraphics[scale=0.02]{pl3_7.eps}\;\includegraphics[scale=0.015]{pl1_2C.eps})
-
\bar V(\includegraphics[scale=0.02]{pl3_11.eps}\;\includegraphics[scale=0.015]{pl1_2C.eps})
-
\bar V(\includegraphics[scale=0.02]{pl3_12.eps}\;\includegraphics[scale=0.015]{pl1_2C.eps})
\\
&&
-
\bar V(\includegraphics[scale=0.02]{pl2_1.eps}\;\includegraphics[scale=0.02]{pl2_5.eps})
-
\bar V(\includegraphics[scale=0.02]{pl2_2.eps}\;\includegraphics[scale=0.02]{pl2_5.eps})
-
\bar V(\includegraphics[scale=0.02]{pl2_3.eps}\;\includegraphics[scale=0.02]{pl2_5.eps})
-
\bar V(\includegraphics[scale=0.02]{pl2_4.eps}\;\includegraphics[scale=0.02]{pl2_5.eps})
-
\bar V(\includegraphics[scale=0.02]{pl2_5.eps}\;\includegraphics[scale=0.02]{pl2_5.eps})
-
\bar V(\includegraphics[scale=0.02]{pl2_5.eps}\;\includegraphics[scale=0.02]{pl2_6.eps})
\\&&
+2 V(\includegraphics[scale=0.02]{pl4_2.eps})
+ V(\includegraphics[scale=0.02]{pl4_4.eps})
+ V(\includegraphics[scale=0.02]{pl4_7.eps})
+ V(\includegraphics[scale=0.02]{pl4_8.eps})
+ V(\includegraphics[scale=0.02]{pl4_9.eps})
+ V(\includegraphics[scale=0.02]{pl4_10.eps})
\\&&
+ V(\includegraphics[scale=0.02]{pl4_11.eps})
+ V(\includegraphics[scale=0.02]{pl4_12.eps})
+2 V(\includegraphics[scale=0.02]{pl4_15.eps})
+2 V(\includegraphics[scale=0.02]{pl4_18.eps})
+ V(\includegraphics[scale=0.02]{pl4_21.eps})
+ V(\includegraphics[scale=0.02]{pl4_22.eps})
\\&&
+ V(\includegraphics[scale=0.02]{pl4_23.eps})
+ V(\includegraphics[scale=0.02]{pl4_26.eps})
+ V(\includegraphics[scale=0.02]{pl4_29.eps})
+ V(\includegraphics[scale=0.02]{pl4_30.eps})
+ V(\includegraphics[scale=0.02]{pl4_33.eps})
+ V(\includegraphics[scale=0.02]{pl4_35.eps})
+ V(\includegraphics[scale=0.02]{pl4_36.eps})
;\\
V(\includegraphics[scale=0.02]{pl2_6.eps}\;\includegraphics[scale=0.015]{pl1_2C.eps}\;\includegraphics[scale=0.015]{pl1_2C.eps})
&=&
\bar V(\includegraphics[scale=0.02]{pl2_6.eps}\;\includegraphics[scale=0.015]{pl1_2C.eps}\;\includegraphics[scale=0.015]{pl1_2C.eps})
-
\bar V(\includegraphics[scale=0.02]{pl3_2.eps}\;\includegraphics[scale=0.015]{pl1_2C.eps})
-
\bar V(\includegraphics[scale=0.02]{pl3_4.eps}\;\includegraphics[scale=0.015]{pl1_2C.eps})
-
\bar V(\includegraphics[scale=0.02]{pl3_8.eps}\;\includegraphics[scale=0.015]{pl1_2C.eps})
-
\bar V(\includegraphics[scale=0.02]{pl3_12.eps}\;\includegraphics[scale=0.015]{pl1_2C.eps})
-
\bar V(\includegraphics[scale=0.02]{pl3_13.eps}\;\includegraphics[scale=0.015]{pl1_2C.eps})
\\
&&
-
\bar V(\includegraphics[scale=0.02]{pl2_1.eps}\;\includegraphics[scale=0.02]{pl2_6.eps})
-
\bar V(\includegraphics[scale=0.02]{pl2_2.eps}\;\includegraphics[scale=0.02]{pl2_6.eps})
-
\bar V(\includegraphics[scale=0.02]{pl2_3.eps}\;\includegraphics[scale=0.02]{pl2_6.eps})
-
\bar V(\includegraphics[scale=0.02]{pl2_4.eps}\;\includegraphics[scale=0.02]{pl2_6.eps})
-
\bar V(\includegraphics[scale=0.02]{pl2_5.eps}\;\includegraphics[scale=0.02]{pl2_6.eps})
-
\bar V(\includegraphics[scale=0.02]{pl2_6.eps}\;\includegraphics[scale=0.02]{pl2_6.eps})
\\&&
+ V(\includegraphics[scale=0.02]{pl4_1.eps})
+2 V(\includegraphics[scale=0.02]{pl4_5.eps})
+2 V(\includegraphics[scale=0.02]{pl4_7.eps})
+2 V(\includegraphics[scale=0.02]{pl4_11.eps})
+ V(\includegraphics[scale=0.02]{pl4_13.eps})
+ V(\includegraphics[scale=0.02]{pl4_14.eps})
\\&&
+ V(\includegraphics[scale=0.02]{pl4_15.eps})
+ V(\includegraphics[scale=0.02]{pl4_16.eps})
+ V(\includegraphics[scale=0.02]{pl4_17.eps})
+ V(\includegraphics[scale=0.02]{pl4_18.eps})
+ V(\includegraphics[scale=0.02]{pl4_22.eps})
+ V(\includegraphics[scale=0.02]{pl4_23.eps})
\\&&
+ V(\includegraphics[scale=0.02]{pl4_24.eps})
+ V(\includegraphics[scale=0.02]{pl4_25.eps})
+ V(\includegraphics[scale=0.02]{pl4_27.eps})
+ V(\includegraphics[scale=0.02]{pl4_30.eps})
+ V(\includegraphics[scale=0.02]{pl4_31.eps})
+ V(\includegraphics[scale=0.02]{pl4_34.eps})
+ V(\includegraphics[scale=0.02]{pl4_36.eps}).
\end{eqnarray*}
The sum of these 6 terms is
\begin{equation}
\frac{3}{8}(n-2)(9n^5-18n^4-96n^3+198n^2+235n-520).
\end{equation}

Subtracting the 6 $V$-terms of rank 1, the
35 $V$-terms of rank 2, and the 36 $V$-terms of rank 3
from $\binom{M_{n-1}}{4}$ yields
\cite{DresdenLozenge}\cite[A326369]{sloane}
\begin{multline}
L_{n,4} = 
\binom{M_{n-1}}{4}-\bigg[3(6n-11)(n-2)+\frac{3}{2}(n-2)(10n^3-20n^2-47n+95)\\
+\frac38(n-2)(9n^5-18n^4-96n^3+198n^2+235n-520)\bigg]
\\
=
\frac{3}{128}(n-2)(n-3)(9n^6+9n^5-135n^4-81n^3+670n^2+104n-1216),\quad n\ge 3.
\end{multline}

\subsection{Five or Six Lozenges}

The common shape of $L_{n,l}$ is that $M_{n-1}$ is a polynomial of degree 2,
and that $\binom{M_{n-1}}{l}$, the upper bound \eqref{eq.bino}, is a polynomial of degree $2l$. The first-order
corrections of the leading term $3(n-1)^2\binom{M_{n-1}-2}{l-2}$, counting
uncorrelated $V$ subgraphs, are of lesser
degree $2l-2$, because at larger $n$ the spatial correlations of the (deleted) internal edges
play a lesser role.
Conjectures for polynomials ensue assuming that the second-order corrections
are of degree $2l-4$, once a sufficiently large set of $L_{n,l}$ for small $n$ is known.

\begin{con}
\begin{multline}
L_{n,5} = 
\binom{M_{n-1}}{5}-\bigg[3(n-1)^2\binom{M_{n-1}-2}{3}
\\
-\frac{1}{4}(4704-3102n+1845n^3-2031n^2+60n^4-315n^5+63n^6)\bigg]\\
=
\frac{3}{1280}(n-3)(n+3)(27n^8-135n^7-387n^6+2835n^5-168n^4-18732n^3+19568n^2+36992n-56320)
,\quad n\ge 3.
\end{multline}
\end{con}

\begin{con}
\begin{multline}
L_{n,6} = 
\binom{M_{n-1}}{6}-\bigg[3(n-1)^2\binom{M_{n-1}-2}{4}
\\
-\frac{1}{16}(-131088+61472n-41206n^3+69420n^2-90n^6-918n^7+153n^8-10851n^4+9828n^5)\bigg]\\
=
\frac{1}{5120}(
81 n^{12}
-486 n^{11}
-2835 n^{10}
+21870 n^9
+26775 n^8
-384786 n^7
+131751 n^6
+3275730 n^5
\\
-3798716 n^4
-13254088 n^3
+22481984 n^2
19678080 n
-42024960
),\quad n\ge 4.
\end{multline}
\end{con}

\section{Summary}
Table \ref{tab.main} summarizes the numerical results which were calculated by the
Java program listed in the ancillary directory.

\begin{table}[htb]
\begin{tabular}{r|rrrrrrrr}
$n\backslash l$ & 0 & 1 & 2 & 3 & 4& 5 & 6 &7 \\
\hline
1 & 1 \\
2 & 1 & 3 \\
3 & 1 & 9 & 24 & 18\\
4 & 1 & 18 & 126& 434& 762& 630& 187\\
5 & 1 & 30 & 387 & 2814 & 12699 & 36894& 69242 & 81936 \\
6 & 1 & 45 & 915 & 11127 & 90270 & 515970& 2139120 & 6523428 \\
7 & 1 & 63 & 1845 & 33365 & 417435 & 3836439 & 26841853 & 146208393 \\
8 & 1 & 84 & 3339 & 83568 & 1478160 & 19662060 & 204334715 & 1701554868 \\
9 & 1 & 108 & 5586 & 184254 & 4354497 & 78536358 & 1124301411 & 13119112488 \\
10 & 1 & 135 & 8802 & 369254 & 11203269 & 261985815 & 4914087052 & 75970268748 \\
11 & 1 & 165 & 13230 & 686952 & 25970895 & 762098799 & 18070041680 & 355864850838 \\
12 & 1 & 198 & 19140 & 1203930 & 55414395 & 1990014156 & 58055896449 & 1414611219018 \\
13 & 1 & 234 & 26829 & 2009018 & 110505120 & 4761037260 & 167316709165 & 4931688363498 \\
14 &  1 & 273 & 36621 & 3217749 & 208300257 & 10594451901 & 440911546295 & 15439933756251 \\
15 &  1 & 315 & 48867 & 4977219 & 374375664 & 22178743326 & 1077784772922 & 44182928710470 \\
\end{tabular}
\begin{tabular}{r|rrrrrrrrrrrrrrrrrr}
$n\backslash l$ & 8& 9& 10 & 11 \\
\hline
5 & 57672 & 21432 & 3135\\
6 & 14683401 & 24256853 & 28975770 & 24383838 \\
7 & 628823088  & 2153224090 & 5892984618 & 12892017948 \\
8 & 11554013295 & 64766667704 & 302315092020 & 1181998895448 \\
9 & 127156871457 & 1038068322606 & 7212713283360 & 42993319234518 \\
10 & 987147811836 & 10940096605816 & 104581114754595 & 869988063985737 \\
11 & 5938169156829 & 85230974965513 & 1064629166358066 & 11681266282861098 \\
12 & 29375579984238 & 527873999198830 & 8307168403048731 & 115585010198220444 \\
13 & 124419130905960 & 2728420121843584 & 52640100670770348 & 902231390539173210 \\
14 & 464317587238419 & 12178604171344167 & 282021772415608164 & 5822744874311864316 \\
15 & 1559497806005040 & 48137813623437500 & 1315457502665712336 & 32139701729335767774 \\
\end{tabular}

\caption{The number $L_{n,l}$ of lozenge tilings for sides circumscribing
regular triangles of edge length $n\ge 1$
with $l\ge 0$ unit lozenges. Sequence \cite[A273464]{sloane} contains an extended 
version of this array.
}
\label{tab.main}
\end{table}
Row sums $\sum_{l\ge 0}L_{n,l}$ are 1, 4, 52, 2158, 286242, 121479420,\ldots
Following the conjectured (\ref{eq.lmax}), 
the maximum $l$ for nonzero entries is $T_{n-2}$.

The values $L_{n,T_{n-2}}$ = 1, 3, 18, 187, 3135, 81462, 3198404, 186498819,
15952438877, 1983341709785 ,\ldots of these configurations with the maximum
number of lozenges
have already been computed by Santos \cite[Table 1]{SantosCMath374}\cite[Table 9.2]{SantosACM25}.

\appendix
\section{JAVA Program}
\subsection{Algorithm}
The ancillary directory contains a Java program that generates
Table \ref{tab.main}.
The main function in \texttt{LozeTil2.java} uses an edge-growing
recursive algorithm which computes a lozenge statistics $L_{n,l}$
refined by the set of lozenges that have one of their 4 sides
on one of the three sides of the big triangle.

\bibliographystyle{amsplain}
\bibliography{all}

\end{document}